\documentclass[a4paper,11pt]{amsart}

\usepackage{amsmath,amsthm,amssymb}
\usepackage[mathscr]{eucal}
 \usepackage{cite}
\usepackage{upgreek}
\usepackage[bookmarks=false]{hyperref}
\usepackage{enumerate}
\usepackage{color}

\numberwithin{equation}{section}

\theoremstyle{plain}
\newtheorem{main}{Theorem}
\newtheorem{mcor}[main]{Corollary}

\newtheorem{theorem}{Theorem}[section]

\newtheorem{lemma}[theorem]{Lemma}
\newtheorem{proposition}[theorem]{Proposition}
\newtheorem{corollary}[theorem]{Corollary}
\theoremstyle{definition}
\newtheorem{definition}[theorem]{Definition}

\newtheorem{notation}[theorem]{Notation}
\newtheorem{remark}[theorem]{Remark}

\newcommand{\R}{\mathbb{R}}
\newcommand{\C}{\mathbb{C}}
\newcommand{\Z}{\mathbb{Z}}

\newcommand{\G}{\mathbb{G}}

\renewcommand{\P}{\mathbb{P}}

\newcommand{\cH}{\mathcal{H}}

\newcommand{\cK}{\mathcal{K}}

\newcommand{\cS}{\mathcal{S}}

\newcommand{\fg}{\mathfrak{g}}
\newcommand{\fh}{\mathfrak{h}}
\renewcommand{\d}{\text{d}}

\newcommand{\Ad}{\operatorname{Ad}}

\newcommand{\SL}{\operatorname{SL}}
\newcommand{\GL}{\operatorname{GL}}

\newcommand{\Prob}{\operatorname{Prob}}
\newcommand{\supp}{\operatorname{supp}}
\newcommand{\rk}{\operatorname{rk}}

\newcommand{\Ker}{\operatorname{Ker}}
\newcommand{\Range}{\operatorname{Range}}
\newcommand{\End}{\operatorname{End}}
\newcommand{\Char}{\operatorname{Char}}
\newcommand{\actson}{\curvearrowright}

\newcommand{\oneG}{\overline{\Gamma}^1}

\def\bbq{\mathbb{Q}}

\def\bbr{\mathbb{R}}

\def\bbc{\mathbb{C}}

\def\bbh{\mathbb{H}}

\def\bbg{\mathbb{G}}

\def\ocal{\mathcal{O}}
\def\ocal{\mathcal{O}}

\def\hfr{\mathfrak{h}}

\def\Bfr{\mathfrak{B}}

\def\pfr{\mathfrak{p}}

\def\gfr{\mathfrak{g}}

\def\vbf{\mathbf{v}}

\def\lbf{\mathbf{l}}

\DeclareMathOperator\Lie{Lie}

\def\h{\hspace{1mm}}

\def\be{\begin{equation}}
\def\ee{\end{equation}}

\begin{document}

\title[Local spectral gap]
{Local spectral gap in simple Lie groups and applications}
\author{R\'{e}mi Boutonnet, Adrian Ioana and Alireza Salehi Golsefidy}
\thanks {R.B. was partially supported by NSF Grant DMS \#1161047, NSF Career Grant DMS \#1253402 and ANR Grant NEUMANN}
\thanks{A.I. was partially supported by NSF  Grant DMS \#1161047,  NSF Career Grant DMS \#1253402, and a Sloan Foundation Fellowship}
\thanks{A.S.G. was partially supported by NSF Grant DMS \#1303121, and a Sloan Foundation Fellowship}
\address{Mathematics Department; University of California, San Diego, CA 90095-1555 (United States).}
\email{rboutonnet@ucsd.edu}
\email{aioana@ucsd.edu}
\email{asalehigolsefidy@ucsd.edu}

\begin{abstract}  We introduce a novel notion of {\it local spectral gap} for general, possibly infinite, measure preserving actions. We establish local spectral gap for the left translation action $\Gamma\curvearrowright G$, whenever $\Gamma$ is a dense subgroup generated by algebraic elements of an arbitrary connected simple Lie group $G$. 
This extends to the non-compact setting works of Bourgain and Gamburd \cite{BG06,BG10}, and Benoist and de Saxc\'{e} \cite{BdS14}. We present several applications to the Banach-Ruziewicz problem, orbit equivalence rigidity,  continuous and monotone expanders, and bounded random walks on $G$. In particular, we prove that, up to a multiplicative constant, the Haar measure is the unique $\Gamma$-invariant finitely additive measure defined on all bounded measurable subsets of $G$. \end{abstract}

\maketitle


\tableofcontents

\section{Introduction}

\subsection{Background and motivation} Spectral gap for probability measure preserving actions is a fundamental notion in mathematics with a wide range of applications.  
The goal of this paper is to introduce and study a notion of spectral gap for general measure preserving actions. 

We begin our discussion by recalling the following:

\begin{definition}
A measure preserving action $\Gamma\curvearrowright (X,\mu)$ of a countable group $\Gamma$ on a standard probability space $(X,\mu)$ is said to have {\it spectral gap}  if there exist $S\subset\Gamma$ finite and  $\kappa>0$ such that
$$\|F\|_{2}\leqslant\kappa\sum_{g\in S}\|g\cdot F-F\|_{2}\;\;\;\text{for any $F\in L^2(X,\mu)$ with $\int_{X}F\;\text{d}\mu=0$}.$$
Here,  $g\cdot F$ denotes the function given by $(g\cdot F)(x)=F(g^{-1}x)$, for every $g\in\Gamma$ and $x\in X$.
To justify the terminology, consider the self-adjoint averaging operator $P_S(\xi)=\frac{1}{2|S|}\sum_{g\in S}(g\cdot F+g^{-1}\cdot F)$. 
Then the constant function {\bf 1} is an eigenfunction of $P_S$ with eigenvalue $1$, and the existence of $\kappa>0$ as above is equivalent to the presence of a gap right below $1$ in the spectrum of $P_S$. 
\end{definition}

Let $G$ be a compact Lie group and denote by $m_G$ its Haar measure.
An important question, which has been investigated intensively over the last three decades, is whether the left translation action $\Gamma\curvearrowright (G,m_G)$, associated to a countable dense subgroup $\Gamma<G$, has spectral gap. 
Interest in this question first arose in the early 1980s, in connection with  Ruziewicz's problem for the $n$-sphere $S^n$ (also known as the Banach-Ruziewicz problem). The latter asks if the Lebesgue  measure on $S^n$ is the unique finitely additive, rotation-invariant measure defined on all Lebesgue measurable subsets. For $n=1$,  Banach used the amenability of $SO(2)$ (as a discrete group) to show that the answer is negative \cite{Ba23}.  For $n\geqslant 2$, however, the problem remained open for a long time.
First, it was realized that  the existence of a countable dense subgroup of $SO(n+1)$ with the spectral gap property implies an affirmative answer \cite{dJR79,Ro81}.  By using Kazhdan's property (T), Margulis \cite{Ma80} and Sullivan \cite{Su81} then obtained an affirmative answer for every  $n\geqslant 4$. The remaining cases $n=2, 3$ were finally settled in the affirmative by Drinfeld \cite{Dr84} via the construction of a countable dense subgroup of $SU(2)$ with the spectral gap property.
An optimal such construction was achieved soon after by Lubotzky, Phillips, and Sarnak \cite{LPS86,LPS87} (see \cite{Oh05} for a generalization to compact simple Lie groups not locally isomorphic to $SO(3)$). For all of this, see the excellent survey \cite{Lu94}.
Later on, a new robust method for proving the spectral gap property for subgroups of $SU(2)$ was developed by Gamburd, Jakobson, and Sarnak \cite{GJS99}. It is worth pointing out that in all of these results, the subgroups involved are generated by matrices with algebraic entries. 

In 2006, a breakthrough was made by Bourgain and Gamburd who established the spectral gap property for any dense subgroup of $SU(2)$ generated by matrices with algebraic entries \cite{BG06}. Their approach followed two earlier major works: the authors' work on expansion for Cayley graphs of $SL_2(\mathbb F_p)$ \cite{BG05}, and Helfgott's product theorem for subsets of $SL_2(\mathbb F_p)$ \cite{He05}. 
Subsequently,  Bourgain and Gamburd established the spectral gap property for dense subgroups of $G=SU(d)$ generated by matrices with algebraic entries, for any $d\geqslant 2$ \cite{BG10}. Recently, this was generalized further by Benoist and de Saxc\'{e} to cover arbitrary connected compact simple Lie groups $G$ \cite{BdS14}.

If $G$ is a compact group and $\Gamma$ is a countable dense subgroup with the spectral gap property, then the Haar measure $m_G$ is the unique finitely additive $\Gamma$-invariant measure defined on all measurable subsets of $G$.
One of the main motivations for this paper is to formulate and prove analogues of the results from \cite{BG06,BG10,BdS14} that apply to general simple Lie groups $G$.  By analogy with the compact case, it would thus be desirable to find a notion of spectral gap for infinite measure preserving actions, which in the case of left translation actions on locally compact groups $G$, implies a uniqueness property for its left Haar measures as finitely additive measures.
\subsection{Local spectral gap} 
As we explain in Corollary \ref{BR}, the following new notion of spectral gap satisfies the desired property.

\begin{definition}
Let $\Gamma\curvearrowright (X,\mu)$ be a measure preserving action of a countable group $\Gamma$ on a standard measure space $(X,\mu)$.  We say that $\Gamma\curvearrowright (X,\mu)$ has {\it local spectral gap} with respect to a measurable set $B\subset X$ of finite measure if there exist $S\subset\Gamma$ finite and $\kappa>0$ such that
$$\|F\|_{2,B}\leqslant\kappa\sum_{g\in S}\|g\cdot F-F\|_{2,B}\;\;\;\text{for any $F\in L^2(X,\mu)$ with $\int_{B}F\;\text{d}\mu=0$}.$$
Here, $\|F\|_{2,B}:=\displaystyle{\big(\int_{B}|F|^2\;\text{d}\mu\big)^{\frac{1}{2}}}$ denotes the $L^2$-norm of the restriction of $F$ to $B$. 
\end{definition}
\begin{remark}\label{rem1} We continue with a few remarks on this definition:
\begin{enumerate}
\item Although the action $\Gamma\curvearrowright (X,\mu)$ is not required to be ergodic, this is automatic if the action has local spectral gap with respect to a set $B$ such that $\cup_{g\in\Gamma}\;g\cdot B$ is co-null in $X$.

\item If $(X,\mu)$ is a probability space and $B=X$,  then local spectral gap coincides with  spectral gap. Assume that $(X,\mu)$ is an infinite measure space with $X$ being a locally compact space and $\mu$ a Radon measure. Then the notion of local spectral gap aims to capture the intuitive idea that functions on $X$ which are locally almost $\Gamma$-invariant, must be locally almost constant (see also Proposition \ref{kazhdan}). This is different from the ``global"  notion of spectral gap requiring that there is no sequence of unit almost $\Gamma$-invariant functions in $L^2(X,\mu)$. Indeed, since any sequence of unit almost $\Gamma$-invariant functions in $L^2(X,\mu)$ converges weakly to $0$ on compact subsets of $X$, the latter reflects only the dynamics of the action at infinity.

\item The notion of local spectral gap appears implicitly in Margulis'   positive resolution of the Banach-Ruziewicz problem for $\mathbb R^n$ ($n\geqslant 3$). 
More precisely, with the above terminology,
he first shows the existence of a  subgroup $\Gamma<\mathbb R^n\rtimes SO(n)$ such that the action $\Gamma\curvearrowright (\mathbb R^n,\lambda^n)$ has local spectral gap, and then concludes that the Lebesgue measure $\lambda^n$ is indeed the unique finitely additive isometry-invariant measure defined on all bounded measurable subsets of $\mathbb R^n$ \cite{Ma82}. 
  
\item  While local spectral gap might depend on the choice of $B$, the following independence result can be easily shown: assume that $B_1,B_2$ are measurable subsets of $X$ such that $B_1\subset K\cdot B_2$ and $B_2\subset K\cdot B_1$, for some finite set $K\subset\Gamma$. Then local spectral gap with respect to $B_1$ is equivalent to local spectral gap with respect to $B_2$ (see Proposition \ref{indep}).

\end{enumerate}
\end{remark}
{\bf Notation.} 
Let $G$ be locally compact second countable group and $H<G$ be a closed subgroup. 
 Here and after, we assume that the locally compact topology on $G$ is Hausdorff. 
We denote by $m_G$ a fixed left Haar measure of $G$. We also denote by $m_{G/H}$ a fixed quasi-invariant Borel regular measure on $G/H$ which  is ``nice", in the sense that it arises from a rho-function for the pair $(G,H)$ (see  \cite[Theorem B.1.4.]{BdHV08}). 

The following is our main result. 

\begin{main}[local spectral gap]\label{main} 
Let $G$ be a connected simple Lie group. Denote by $\frak g$ the Lie algebra of $G$ and by ${\Ad}:G\rightarrow\GL(\frak g)$ its adjoint representation.
Let $\Gamma<G$ be a dense subgroup. Assume that there is a basis $\frak B$ of $\frak g$ such that the matrix of $\Ad(g)$ in the basis $\frak B$ has algebraic entries, for any $g\in\Gamma$.  Let $B\subset G$ be a measurable set with compact closure and non-empty interior.

Then the left translation action $\Gamma\curvearrowright (G,m_G)$ has local spectral gap with respect to $B$. \end{main}

In the case $G$ is compact, Theorem \ref{main} recovers the main results of \cite{BG06,BG10,BdS14}. On the other hand, if $G$ is not compact, Theorem \ref{main} reveals an entirely new type of phenomenon for locally compact groups.


\begin{remark}
The assumption on $\Gamma<G$ is in particular satisfied if $G=\SL_n(\mathbb R)$, for some $n\geqslant 2$, and $\Gamma$ is a dense subgroup of $G$ such that every matrix $g\in\Gamma$ has algebraic entries. 
\end{remark}

\begin{remark} In view of Remark \ref{rem1} (4), the conclusion of Theorem \ref{main} does not depend on the choice of the set $B$. Indeed, if $B_1\subset G$ has compact closure and $B_2\subset G$ has non-empty interior, then there exists a finite set $K\subset\Gamma$ such that $B_1\subset K\cdot B_2$.
\end{remark}

 Theorem \ref{main} is a consequence of our main technical result proving a restricted spectral gap estimate in the spirit of Bourgain and Yehudayoff's pioneering work  \cite{Bo09,BY11}.

\begin{main}[restricted spectral gap]\label{restricted}
Assume that  $\Gamma<G$ are as in Theorem \ref{main}. 
Let $B\subset G$ be a measurable set with compact closure. Let $U$ be a neighborhood of the identity element in $G$.

Then there exist a finite set $T\subset\Gamma\cap U$ and a finite dimensional subspace $V\subset L^2(B)$ such that the probability measure $\mu=\frac{1}{2|T|}\sum_{g\in T}(\delta_g+\delta_{g^{-1}})$ satisfies  $\|\mu*F\|_2<\frac{1}{2}\|F\|_2$, for every $F\in L^2(B)\ominus V$.
\end{main}

Note that unlike Theorem \ref{main}, this result is new even in the case of compact groups, where it leads to some unexpected consequences (see Remark \ref{spectra}).

Theorem \ref{restricted} concerns the left regular representation of $G$. The proof of Theorem \ref{restricted} moreover shows that
  for any $0<r<1$ there exists a finite set $T\subset\Gamma\cap U$ such that the conclusion holds with $r$ in place of $\frac{1}{2}$. As a consequence, it follows that a more general statement, addressing all quasi-regular representations of $G$, holds true.

 \begin{mcor}\label{by} 
 Assume that  $\Gamma<G$ are as in Theorem \ref{main}.
 Let $H<G$ be a closed subgroup  and denote by $\pi:G\rightarrow \mathcal U(L^2(G/H,m_{G/H}))$ be the associated quasi-regular unitary representation.  Let $B\subset G/H$ be a measurable set with compact closure. Let   $U$ be a neighborhood of the identity in $G$.

Then there exist a finite set $T\subset\Gamma\cap U$ and a finite dimensional subspace $V\subset L^2(B)$ such that  the probability measure $\mu=\frac{1}{2|T|}\sum_{g\in T}(\delta_g+\delta_{g^{-1}})$ satisfies  $\|\pi(\mu)(F)\|_2<\frac{1}{2}\|F\|_2$, for any $F\in L^2(B)\ominus V$.
\end{mcor}

Here, for a probability measure $\mu$, we denote by $\pi(\mu)$ the averaging operator $\sum_{g\in G}\mu(\{g\})\pi(g)$.

Corollary \ref{by} generalizes \cite[Theorem 5]{BY11} which deals with the case when $G$ is $SL_2(\mathbb R)$ and $H$ is the subgroup of upper triangular matrices. Then $G/H$ can be identified with the real projective line, $\mathbb P^1(\mathbb R)$.
The proof of \cite[Theorem 5]{BY11}  is specific to this situation, as it relies on the fact that the action of $SL_2(\mathbb R)$ on $\mathbb P^1(\mathbb R)$ is $2$-transitive to show a certain mixing property. Corollary \ref{by} provides an alternative approach to the mixing property in this case. Corollary \ref{by} is new in all other cases with $G$ non-compact, including  the simplest one when $G=SL_2(\mathbb R)$ and $H$ is trivial.

\begin{remark}\label{quant}
In the case $G$ has trivial center, the proof of Theorem \ref{restricted} yields a more quantitative statement (see Theorem \ref{restricted2}). 
 To explain this, identify $G$ with a subgroup of $\GL_n(\mathbb R)$, for some $n$, and endow it with the metric induced by the Hilbert-Schmidt norm $\|.\|_2$. For $\varepsilon>0$,  denote $B_{\varepsilon}(1)=\{g\in G|\;\|g-1\|_2<\varepsilon\}$.
  
 Then the proof of Theorem \ref{restricted} shows that there is a constant $C>1$ (depending on $\Gamma$) such that for any small enough $\varepsilon>0$, there exist a finite set $T\subset\Gamma\cap B_{\varepsilon}(1)$ which freely generates a group, and a finite dimensional subspace  $V\subset L^2(B)$ such that denoting $\mu=\frac{1}{2|T|}\sum_{g\in T}(\delta_g+\delta_{g^{-1}})$ we have
 \begin{itemize}

 \item $|T|<\frac{1}{\varepsilon^{C}}$, and
 \item   $\|\mu*F\|_2<\varepsilon\|F\|_2$, for every $F\in L^2(B)\ominus V$.
\end{itemize}
\end{remark}

\begin{remark}\label{spectra}
Theorem \ref{restricted} (and its quantitative version) sheds some new light on the spectra of averaging operators on compact groups. 
In order to briefly recall known results along these lines, assume for simplicity that $G=SU(2)$. 
Then the irreducible representations of $G$ can be listed as $\pi_n:G\rightarrow\mathcal U(\mathcal H_n)$, where dim$(\mathcal H_n)=n+1$, for every $n\geqslant 0$, and by the Peter-Weyl theorem we have that $L^2(G)=\bigoplus_{n\geqslant 0}\mathcal H_n^{\oplus {(n+1)}}$.
Let $T\subset G$ be a finite set which freely generates a subgroup, consider the probability measure $\mu=\frac{1}{2|T|}\sum_{g\in T}(\delta_g+\delta_{g^{-1}})$, and denote by $P_{\mu}$ the operator $F\mapsto\mu*F$. 
Then $P_{\mu}$ is self-adjoint and since $\|P_{\mu}\|\leqslant 1$, its spectrum is contained in $[-1,1]$.
Moreover, since $P_{\mu}$ can be identified with $\bigoplus_{n\geqslant 0}\pi_n(\mu)^{\oplus {n+1}}$, it is also diagonalizable. 

 The  asymptotic distribution of the eigenvalues of $P_{\mu}$ has been studied in  \cite{LPS86, GJS99}, where it is shown that most of them lie in the interval $\big[-\frac{\sqrt{2|T|-1}}{|T|},\frac{\sqrt{2|T|-1}}{|T|} \big]$.
More precisely, if $d_n$ denotes the number of  eigenvalues of $\pi_n(\mu)$ that lie outside this interval (so-called ``exceptional" eigenvalues), then $\frac{d_n}{n}\rightarrow 0$ (see \cite[Theorem 1.1]{LPS86}). 
Assume from now on that the  elements of $T$ have algebraic entries. Then the following sharper estimate holds: $\frac{d_n}{n}\ll\frac{1}{\log n}$, for large $n$ (see \cite[Theorem 1.3]{GJS99}). A remarkable fact, discovered by Lubotzky, Phillips and Sarnak  is that for certain  sets  $T$, the operator $P_{\mu}$ has no exceptional eigenvalues, i.e. $d_n=0$, for all $n\geqslant 1$ (see \cite{LPS86,LPS87}).
As already mentioned above, the more recent work \cite{BG06} implies that $P_{\mu}$ has a spectral gap, i.e. the supremum $\kappa_{\mu}$ of the spectrum of $P_{\mu}$ but $1$ satisfies $\kappa_{\mu}<1$.
However, besides these facts, not much is known about the exceptional eigenvalues of $P_{\mu}$. In particular,  to the best of our knowledge, it is unknown whether $\kappa_{\mu}$ is ever an eigenvalue of $P_{\mu}$.

 Theorem \ref{restricted} implies that $k_{\mu}$ can be an isolated eigenvalue of $P_{\mu}$, and thus $P_{\mu}$ can have a second spectral gap.  
Moreover, it shows that operators of the form $P_{\mu}$ may have arbitrarily many gaps at the top of their spectrum.
To make this precise, let $\varepsilon>0$ small enough, and let $T\subset\Gamma\cap B_{\varepsilon}(1)$ and $\mu$ as given  by Remark \ref{quant}.
Then $P_{\mu}$ has only finitely many eigenvalues outside the interval $[-\varepsilon,\varepsilon]$. 
On the other hand,  since $T\subset B_{\varepsilon}(1)$, the number of eigenvalues of $P_{\mu}$ belonging to the interval $(\frac{1}{2},1)$ gets arbitrarily large, as $\varepsilon\rightarrow 0$. In fact, it is easy to see that this number is $\gg\frac{1}{\varepsilon^2}$. 
In combination with \cite[Theorem 1.1]{LPS86} the following picture emerges: 
the spectrum of $P_{\mu}$ contains 
\begin{itemize}
\item the whole interval $\big[-\frac{\sqrt{2|T|-1}}{|T|},\frac{\sqrt{2|T|-1}}{|T|} \big]$ 
\item  only finitely many points, all of which are isolated eigenvalues, outside $\big[-\frac{1}{|T|^{\frac{1}{C}}},\frac{1}{|T|^{\frac{1}{C}}}\big]$.
\item  $\gg |T|^{\frac{2}{C}}$ points in the interval $(\frac{1}{2},1)$.

\end{itemize}
\end{remark}

\subsection{Applications} We now turn to discussing several applications of our main results.

\subsubsection{\bf{The Banach-Ruziewicz problem}} The original Banach-Ruziewicz problem asks whether the Lebesgue measure on $S^n$ (resp. $\mathbb R^n$) is the unique rotation-invariant (resp. isometry-invariant) finitely additive measure defined on all bounded Lebesgue measurable sets.  This problem is an illustration of a general question: 
let $\Gamma$ be a locally compact group acting isometrically on a locally compact metric space $X$ with an invariant Radon measure $\mu$. Is $\mu$ the unique $\Gamma$-invariant finitely additive measure defined on all $\mu$-measurable subsets of $X$ with compact closure? Here and after, uniqueness is of course meant up to a multiplicative constant.  

If the space $X$ is compact and the group $\Gamma$ is countable discrete, then a positive answer to this question is closely connected to the spectral gap of the action. The connection stems from the well-known fact that the action $\Gamma\curvearrowright (X,\mu)$ has spectral gap if and only if  integration against $\mu$ is the unique $\Gamma$-invariant mean on $L^{\infty}(X,\mu)$ (see \cite{Ro81,Sc81}). On the other hand, if $\mu$ is unique among invariant finitely additive measures, then  integration against $\mu$ is unique among invariant means. The converse of this statement is also true for certain classes of actions, including left translation actions on compact groups (see Remark \ref{cominv}).   

In Section \ref{6}, we  generalize these results to the case when  $X$ is locally compact. 
Assume that every orbit $\Gamma\cdot x$ is dense in $X$, and denote by $L^{\infty}_{\text{c}}(X,\mu)$ the algebra of $L^{\infty}$-functions with compact support. Firstly, we prove that the action $\Gamma\curvearrowright (X,\mu)$ has local spectral gap with respect to a measurable set with compact closure and non-empty interior if and only if integration against $\mu$ is the unique $\Gamma$-invariant positive linear functional on $L^{\infty}_{\text{c}}(X,\mu)$ (see Theorem \ref{mean}). This result is partially inspired by Margulis' work \cite{Ma82}, which we follow closely in the proof of the only if assertion. 
Secondly, in the case of left translation actions $\Gamma\curvearrowright (G,m_G)$ on locally compact groups $G$, we show that the Haar measure $m_G$ is unique among invariant finitely additive measures if and only if integration against $m_G$ is unique among invariant positive linear functionals on $L^{\infty}_{\text{c}}(G,m_G)$ (see Theorem \ref{BRtext}).
Altogether, by combining these two results we derive the following:

\begin{main}\label{BR}

Let $G$ be a locally compact second countable group and  $\Gamma<G$ be a countable dense subgroup. Denote by $\mathcal C(G)$ the family of measurable subsets $A\subset G$ with compact closure.

 Then the following conditions are equivalent:

\begin{enumerate}
\item If $\nu:\mathcal C(G)\rightarrow [0,\infty)$ is a $\Gamma$-invariant, finitely additive measure, then there exists $\alpha\geqslant 0$ such that $\nu(A)=\alpha\;m_G(A)$, for all $A\in\mathcal C(G)$.
\item The left translation action $\Gamma\curvearrowright (G,m_G)$ has local spectral gap.

\end{enumerate}
\end{main}

Note that in order to treat arbitrary locally compact groups, we use the structure theory of locally compact groups \cite{MZ55} as well as Breuillard and Gelander's  topological Tits alternative \cite{BG04}. 

As an immediate consequence of Theorems \ref{main} and \ref{BR} we deduce the following uniqueness characterization of Haar measures on simple Lie groups, in the spirit of the Banach-Ruziewicz problem:

\begin{mcor}
Assume that $\Gamma<G$ are as in Theorem \ref{main}. 

Then, up to a multiplicative constant, the Haar measure $m_G$ of $G$ is the unique finitely additive $\Gamma$-invariant  measure defined on $\mathcal C(G)$.
\end{mcor}

\subsubsection{\bf{Orbit equivalence rigidity}} Next, we apply our results to the theory of orbit equivalence of actions.  This area has flourished in the last 15 years, with many new exciting developments (see the surveys \cite{Po07,Fu09,Ga10}). To recall the notion of orbit equivalence, consider two ergodic measure preserving actions $\Gamma\curvearrowright (X,\mu)$ and $\Lambda\curvearrowright (Y,\nu)$  of countable groups $\Gamma$, $\Lambda$ on standard  measure spaces $(X,\mu)$, $(Y,\nu)$. The actions are called {\it orbit equivalent} if there exists a measure class preserving Borel isomorphism $\theta:X\rightarrow Y$ such that $\theta(\Gamma\cdot x)=\Lambda\cdot\theta(x)$, for $\mu$-almost every $x\in X$. The simplest instance of when the actions are orbit equivalent is when they are {\it conjugate}, i.e. there exists a measure class preserving Borel isomorphism $\theta:X\rightarrow Y$ and a group isomorphism $\delta:\Gamma\rightarrow\Lambda$ such that $\theta(g\cdot x)=\delta(g)\cdot\theta(x)$, for all $g\in\Gamma$ and $\mu$-almost every $x\in X$. 

In general, however, orbit equivalence is a much weaker notion of equivalence than conjugacy. This is best illustrated by the striking theorem of Ornstein-Weiss and Connes-Feldman-Weiss  showing that if the groups $\Gamma,\Lambda$ are both infinite amenable and the measure spaces $(X,\mu),(Y,\nu)$ are either both finite or both infinite, then the actions are orbit equivalent (see \cite{OW80,CFW81}).
In sharp contrast, there exist ``rigid" situations when for certain classes of actions of non-amenable groups one can deduce conjugacy  from orbit equivalence. 

It was recently discovered in \cite{Io13} that such a rigidity phenomenon occurs for left translation actions on compact groups in the presence of spectral gap. More precisely, let $\Gamma<G$ and $\Lambda<H$ be countable dense subgroups of compact connected Lie groups with trivial centers. 
Assuming that $\Gamma\curvearrowright (G,m_G)$ has spectral gap, it follows from \cite[Corollary 6.3]{Io13} that the actions $\Gamma\curvearrowright (G,m_G)$ and $\Lambda\curvearrowright (H,m_H)$ are orbit equivalent if and only if they are conjugate. Most recently, this result has been generalized to the case when $G$ and $H$ are arbitrary, not necessarily compact, connected Lie groups with trivial centers (see \cite[Theorems A and 4.1]{Io14}). The only difference is that in the locally compact setting, the spectral gap assumption has to be replaced with the assumption that the action $\Gamma\curvearrowright (G,m_G)$ is strongly ergodic.

To recall the latter notion, let $\Gamma\curvearrowright (X,\mu)$ be  an ergodic measure preserving action. Then, loosely speaking, strong ergodicity requires that any sequence of asymptotically invariant subsets of $X$ must be asymptotically trivial. In order to make this precise, since the measure $\mu$ can be infinite, we first choose a probability measure $\mu_0$ on $X$ with the same null sets as $\mu$. 
The action  is said to be {\it strongly ergodic} 
if any sequence $\{A_n\}$ of measurable subsets of $X$ satisfying $\mu_0(g\cdot A_n\;\Delta\; A_n)\rightarrow 0$, for all $g\in\Gamma$, must satisfy $\mu_0(A_n)(1-\mu_0(A_n))\rightarrow 0$ \cite{CW80,Sc80}. It is easy to see that this definition does not depend on the choice of $\mu_0$.

For translation actions on compact groups, strong ergodicity is implied by the spectral gap property, which is now known to hold  in considerably large generality by  \cite{BG06,BG10,BdS14}. On the other hand, in the case of translation actions on locally compact non-compact groups, strong ergodicity seems much harder to work with, and so far could only be checked in two rather specific situations (see \cite[Propositions G and H]{Io14}).

Nevertheless, strong ergodicity is implied by local spectral gap, for arbitrary ergodic measure preserving actions. Moreover, for translation actions on locally compact groups, we prove that  local spectral gap and strong ergodicity are equivalent (see Theorem \ref{BRtext}). This generalizes \cite[Theorem 4]{AE10}, which dealt with the compact case.
Consequently, all actions covered by Theorem \ref{main} are strongly ergodic, which in combination with \cite{Io14} allows us to conclude the following:

\begin{mcor}\label{OErig}
Assume that $\Gamma<G$ are as in Theorem \ref{main}. Suppose that $G$ has trivial center. Let $H$ be any connected Lie group with trivial center and $\Lambda<H$ be any countable dense subgroup.

Then the left translation actions $\Gamma\curvearrowright (G,m_G)$ and $\Lambda\curvearrowright (H,m_H)$ are orbit equivalent if and only if there is a topological isomorphism $\delta:G\rightarrow H$ such that $\delta(\Gamma)=\Lambda$.
\end{mcor}

\begin{remark} If $\widetilde\Gamma<G$ is a countable subgroup that contains $\Gamma$, then Theorem \ref{main} implies that the action $\widetilde\Gamma\curvearrowright (G,m_G)$ has local spectral gap. 
Hence, Corollary \ref{OErig} remains valid if $\Gamma$ is replaced by $\widetilde\Gamma$.
\end{remark}

\begin{remark}
In the context of Corollary \ref{OErig}, assume moreover that $\Gamma$ is a free group. Since the left translation action $\Gamma\curvearrowright (G,m_G)$ is strongly ergodic, it is not amenable in the sense of \cite{Zi78}. Then \cite[Theorem A]{HV12} (which builds on \cite{OP07,PV11}) implies that $L^{\infty}(G)$ is the unique Cartan subalgebra of the  $L^{\infty}(G)\rtimes\Gamma$, up to unitary conjugacy. In combination with Corollary \ref{OErig}, we deduce that the crossed product von Neumann algebras $L^{\infty}(G)\rtimes\Gamma$ and $L^{\infty}(H)\rtimes\Lambda$ are isomorphic if and only if there is a topological isomorphism $\delta:G\rightarrow H$ such that $\delta(\Gamma)=\Lambda$.

\end{remark}

\subsubsection{\bf{Continuous and monotone expanders}} Our main results also lead to a general construction of continuous and monotone expanders, extending the main result of \cite{BY11}. Expander graphs are infinite families of highly connected sparse finite graphs. It is sometimes desirable to find expander graphs within certain classes of graphs. A finite graph is called {\it monotone} if it is defined by monotone functions. This means that the vertex set of the graph  can be identified with  $[n]=\{1,2,...,n\}$ in such a way that there exist partially defined monotone maps $\varphi_i:[n]\rightarrow [n], 1\leqslant i\leqslant d,$ such that two  vertices $a,b$ are connected iff $b=\varphi_i(a)$, for some $i$. 

Bourgain and Yehudayoff recently found the first explicit construction of constant degree monotone expander graphs \cite{Bo09,BY11}. Their approach is to first build a continuous monotone expander and then discretize it to obtain monotone  expanders. 
In their terminology, a {\it continuous expander} consists of a family of smooth partially defined maps $\varphi_i:B\rightarrow B$, $1\leqslant i\leqslant d$,  where $B$ is a compact subset of a manifold endowed with a finite measure $A\mapsto |A|$, such that the following holds: there is $\kappa>0$ such that for every measurable set $A\subset B$ with $|A|\leqslant\frac{|B|}{2}$, we have $|\cup_{i=1}^d\varphi_i(A)|\geqslant (1+\kappa)|A|.$

As a consequence of Theorem \ref{restricted}, we obtain the following result.

\begin{mcor}\label{monexp}

Assume that $\Gamma<G$ are as in Theorem \ref{main}.
Let $H<G$ be a closed subgroup and $B\subset G/H$ be a  measurable set with compact closure and non-empty interior. For a measurable subset $A\subset G/H$, denote $|A|:=m_{G/H}(A)$.

Then  there exists a finite set $S\subset\Gamma$ finite for which there is a constant $\kappa>0$ such that  for any measurable set $A\subset B$ with $|A|\leqslant\frac{|B|}{2}$ we have
 $$|\big(\cup_{g\in S}g\cdot A\big)\cap B|\geqslant (1+\kappa)|A|.$$
Moreover, if $B$ is open and connected, and  $\varepsilon>0$ is given, then  $S\subset\Gamma$ can be taken inside $B_{\varepsilon}(1)$.
\end{mcor}

Assume that $G$ is equal to $SL_2(\mathbb R)$, $H$ is the subgroup of upper triangular matrices, identify $G/H$ with the real projective line $\mathbb P^1(\mathbb R)=\mathbb R\cup\{\infty\}$, and let $B=[0,1]$.
With this notation,  \cite[Theorem 4]{BY11} provides a finite set $S\subset SL_2(\mathbb Q)$ which satisfies the conclusion of Corollary \ref{monexp}. Moreover, $S$ can  be taken close enough to the identity so that the restriction $\tilde g$ of every $g\in S$ to $B\cap g^{-1}B$ is monotonically increasing. Therefore, the family $\{\tilde g\}_{g\in S}$ is a continuous monotone expander. 

Corollary \ref{monexp} generalizes \cite[Theorem 4]{BY11} by showing the existence of such a set $S$ inside any dense subgroup of $G$ generated by algebraic elements. Note that, as opposed to \cite{BY11}, our construction of $S$ is not explicit. On the other hand, unlike \cite{BY11}, our construction does not rely on the strong Tits alternative from \cite{Br08}.

\subsubsection{\bf{Spectral gap for delayed bounded random walks}} Our last application concerns random walks on Lie groups that are bounded and delayed, in a sense made precise below.  Let $G$ be a connected simple Lie group and $S\subset G$ be a finite symmetric set.  Denote $k=|S|$ and enumerate $S=\{g_1,...,g_k\}$. Let  $B\subset G$ be  a measurable set which is bounded (i.e. has compact closure). 
We define a random walk on $B$ as follows: a given point $x\in B$ moves with probability $\frac{1}{k}$ to each of the points $h_1x, h_2x,...,h_kx$, where $h_i=g_i$, if $g_ix\in B$, and $h_i=e$, if $g_ix\notin B$. In other words, with probability $\frac{1}{k}$, $x$ either moves to $g_ix$ or stays put, depending on whether $g_ix$ belongs to $B$ or not.

The associated transition operator $P_S:L^2(B)\rightarrow L^2(B)$ is given by $$P_S(F)=\frac{1}{k}\sum_{i=1}^k\big( {\bf 1}_{B\cap g_iB}\;g_i\cdot F+{\bf 1}_{B\setminus g_iB}\;F\big), \;\;\;\;\text{for every $F\in L^2(B)$}.$$

Then $P_S$ is symmetric, $\|P_S\|\leqslant 1$, and $P_S({\bf 1}_B)={\bf 1}_B$, where ${\bf 1}_B$ denotes the characteristic function of $B$. 
 Theorem \ref{restricted} allows us to deduce the existence of many sets $S$ such that $P_S$ has a spectral gap.

\begin{mcor}\label{rwalk}

Assume that $\Gamma<G$ are as in Theorem \ref{main}.

Then there exists a finite symmetric set $S\subset\Gamma$ such that the operator $P_S:L^2(B)\rightarrow L^2(B)$ satisfies 
$$\|{P_S}_{|L^2(B)\ominus\mathbb C{\bf 1}_B}\|<1.$$
\end{mcor}

When $G$ is compact and $B=G$, this result is a consequence of \cite{BG06,BG10,BdS14}. Corollary \ref{rwalk} is new in all other cases, including the case when $G$ is compact and $B$ is a proper subset. 

\subsection{On the proof of restricted spectral gap}
Our approach to proving restricted spectral gap is a combination of general results from \cite{dS14,BdS14}, refinements of techniques from  \cite{BG10,SGV11}, and ideas from \cite{BY11} on how to treat non-compact situations.
It relies on the remarkable strategy invented by Bourgain and Gamburd \cite{BG05,BG06} to prove spectral gap in the compact setting.

 To briefly recall this strategy, consider a symmetric probability measure $\mu$ on a compact group $G$, for which we want to establish the spectral gap property.
A first step is to show that the convolution powers of $\mu$ become ``flat''  rather quickly. Then one uses a mixing inequality to deduce spectral gap for the corresponding operator $P_{\mu}:L^2(G)\rightarrow L^2(G)$ given by $P_{\mu}(F)=\mu*F$.

{\bf Flattening.} The term flat roughly means that after ``discretizing'' the group $G$, the measure has a small $2$-norm, compared to the scale at which we discretize $G$. In \cite{BdS14}, Benoist and de Saxc\'e proved a general flattening lemma for connected compact simple Lie groups $G$. They showed that if a measure $\nu$ on  $G$ is not already flat and does not concentrate on any proper closed subgroup of $G$, then its convolution square $\nu \ast \nu$ will be significantly flatter. A repeated application of this result shows that a measure on $G$ with small mass on closed subgroups will flatten rather quickly.

{\bf Escaping subgroups.} Thus, in order to show that $P_{\mu}$ has spectral gap, it is necessary to show that, quite quickly, convolution powers of $\mu$ have small mass on closed subgroups. To guarantee this, one needs (due to the currently available techniques) to impose a diophantine assumption on the support of $\mu$. Typically, one assumes that $\mu$ is supported on finitely many elements with algebraic entries (when viewed as matrices via the adjoint representation).

{\bf Mixing inequality.} The concluding part, deducing spectral gap out of flatness of some small power of $\mu$, relies on a mixing inequality. If $G$ is a finite group, the mixing inequality bounds the norm of the operator $P_{\mu}$ in terms of the $2$-norm of $\mu$ (see \cite{BNP08} and \cite[Proposition 1.3.7]{Ta15}). This step relies on the representation theory of the ambient group $G$. Specifically, one usually uses the idea, due to Sarnak and Xue \cite{SX91}, of exploiting ``high multiplicity" of eigenvalues.    

In the non-compact setting, we will prove restricted spectral gap with a similar strategy. Recall, however, that our aim is somewhat different from the spectral gap property for compact groups. Indeed, we are given a  connected simple Lie group $G$, a dense subgroup $\Gamma<G$ and an open ball $B \subset G$. Our goal is to produce a measure $\mu$ supported on $\Gamma$ and on an arbitrarily small neighborhood of $1$, such that the averaging operator $P_\mu: L^2(B) \to L^2(G)$ has norm less than $\frac{1}{2}$ (after discarding a finite dimensional subspace $V \subset L^2(B)$). Let us emphasize the main differences that occur in the proof.

Firstly, we  show that the mixing inequality still holds in our setting, leading to a result that might be of independent interest (see Theorem \ref{rho}). Our proof is inspired by the ``geometric approach" introduced in \cite{BY11,BG10}, but here we address a far greater level of generality. Also, our proof is elementary, in that it only relies on basic results from the representation theory of $G$, and essentially self-contained. 
Using this inequality, we reduce to the task of producing a measure $\mu$ with support contained in $\Gamma$ and arbitrarily close to $1$, whose convolution powers flatten rather quickly.

The flattening lemma of \cite{BdS14} relies on two main tools: a {\it product theorem} due to de Saxc\'e \cite{dS14}, and the {\it non-commutative Balog-Szemer\'edi-Gowers Lemma} due to Tao \cite{Ta06}. It turns out that these two tools actually hold for general (not necessarily compact) connected simple Lie groups. So, by reproducing the proof of \cite[Lemma 2.5]{BdS14}, we get a similar flattening lemma in the locally compact setting (see Corollary \ref{flattening}). An important aspect is that our lemma only applies to measures whose support is {\it controlled} (relative to the scale at which we discretize $G$).

Next, refining techniques from \cite[Section 3]{SGV11} we construct a measure $\mu$, supported on $\Gamma$ and on an arbitrarily small neighborhood of $1$, that will escape proper subgroups quickly when taking convolution powers (Propositions \ref{ping-pong} and \ref{neigh}). Therefore, we are almost in position to apply the flattening lemma to some convolution powers of $\mu$. However, we need to make sure that these convolution powers still have a controlled support. This amounts to bounding the speed of escape of subgroups in terms of the size of the support of $\mu$. A priori, the measure $\mu$ that we construct does not admit such a nice bound. As in \cite{BY11}, an application of the pigeonhole principle allows us to construct a new measure $\mu'$ with an improved bound.
Then $\mu'$ satisfies all the required assumptions to ensure that it will become flat quickly enough. Finally,  our mixing inequality will allow us to show restricted spectral gap for this new measure $\mu'$. 

We will provide more quantitative statements of the main steps of the proof in Section \ref{outline}.

\subsection{Organization of the paper} Besides the introduction, this paper has seven other sections and an appendix. In Section 2, we establish  some basic properties of local spectral gap, explain how Theorem \ref{main} follows from Theorem \ref{restricted}, and provide a detailed outline of the proof of Theorem \ref{restricted}. Sections 3, 4, and 5 are each devoted to one of the three main parts of the proof of Theorem \ref{restricted}. In Section 6 we conclude the proof of Theorem \ref{restricted} and derive Corollary \ref{by}.
In Sections 7 and 8, we prove Theorem \ref{BR} and Corollaries \ref{monexp}, \ref{rwalk}, respectively. Finally, the Appendix deals with the proof of Lemma \ref{BdS}. 

\subsection{Acknowledgements} We are grateful to Cyril Houdayer, Hee Oh and Peter Sarnak for helpful comments. 


\tableofcontents


\section{Preliminaries}

\subsection{Terminology} We begin by introducing various terminology concerning analysis on groups.
Let $G$ be a 
 locally compact second countable ({\it l.c.s.c.}) group and fix a left Haar measure $m_G$.

 Given a measurable set $A\subset G$  and  a measurable function $f:G\rightarrow\mathbb C$, we denote
 
$$|A|:=m_G(A),\;\;\;\;\int_{G}f(x)\;\text{d}x:=\int_{G}f\;\text{d}m_G,\;\;\;\;\text{and}$$ $$\|f\|_{p,A}:=\|{\bf 1}_Af\|_p=\Big(\int_{A}|f(x)|^p\;\text{d}x\Big)^{\frac{1}{p}}.$$ 

 We denote by
 $\mathcal M(G)$ the family of Borel probability measures on $G$.
Let $f,g:G\rightarrow\mathbb C$ be measurable functions and $\mu,\nu\in\mathcal M(G)$.
 Then the convolution functions $f*g$, $\mu*f:G\rightarrow\mathbb C$ and probability measure $\mu*\nu$ are defined (when the integrals make sense) by the formulae $$(f*g)(x)=\int_{G}f(y)g(y^{-1}x)\;\text{d}y,\;\;\;\;(\mu*f)(x)=\int_{G}f(y^{-1}x)\;\text{d}\mu(y)\;\;\;\;\text{and}$$ $$\int_{G}F\;\text{d}(\mu*\nu)=\int_{G}\int_{G}F(xy)\;\text{d}\mu(x)\text{d}\nu(y)$$
for any continuous $F:G\rightarrow\mathbb C$.
We will often use the following inequalities $$\|f*g\|_2\leqslant \|f\|_1\|g\|_2,\;\;\;\;\|f*g\|_{\infty}\leqslant\|f\|_2\|g\|_2\;\;\;\;\;\text{and}\;\;\;\;\|\mu*f\|_2\leqslant\|f\|_2.$$

Further, we denote by $\check{f}:G\rightarrow\mathbb C$ the function given by $\check{f}(x)=\overline{f(x^{-1})}$.
Similarly,  $\check{\mu}$ is the Borel probability measure given by ${\int_{G}F\;\text{d}\check{\mu}=\int_{G}\check{F}\;\text{d}\mu}$, for any continuous $F:G\rightarrow\mathbb R$. We say that $\mu$ is {\it symmetric} if  $\check{\mu}=\mu$. For $n\geqslant 1$, we denote by $\mu^{*n}$ the $n$-fold convolution product of $\mu$ with itself. 
We also denote by supp$(\mu)$ the {\it support} of $\mu$.
If $\mu$ and $\nu$ have finite support, then  $\check{\mu}(\{x\})=\overline{\mu(\{x^{-1}\})}$ and $(\mu*\nu)(\{x\})=\sum_{y\in G}\mu(\{y\})\nu(\{y^{-1}x\})$, for any $x,y\in G$.

If $G$ is unimodular, we denote by $\lambda,\rho:G\rightarrow\mathcal U(L^2(G))$ the {\it left} and {\it right} regular representations of $G$ given by $\lambda_g(f)(x)=f(g^{-1}x)$, $\rho_g(f)(x)=f(xg)$, for every $f\in L^2(G)$ and any $g,x\in G$. Notice that $\lambda_g(f)=\delta_g*f$ and $\rho_g(f)=f*\delta_{g^{-1}}$, where $\delta_g$ denotes the Dirac measure at $g\in G$. 

Next, we establish a useful result that we will need later on.

\begin{lemma}\label{powers}\label{A^{-1}A} Let $\mu$ be a symmetric Borel probability measure on  $G$ and $n\geqslant 1$. 
Then
\begin{enumerate}
\item $\|\mu^{*n}*f\|_2\geqslant \|\mu*f\|_2^{2n}$, for every $f\in L^2(G)$ with $\|f\|_2=1$. 
\item  $\mu^{*n}(A)^2\leqslant\mu^{*(2n)}(A^{-1}A)$, for every measurable set $A\subset G$.
\end{enumerate}
\end{lemma}

{\it Proof.} (1) 
Since $\mu$ is symmetric,  we have  $\|\mu^{*m}*f\|_2^2=\langle\mu^{*m}*f,\mu^{*m}*f\rangle=\langle\mu^{*2m}*f,f\rangle\leqslant\|\mu^{*2m}*f\|_2$, for any $m\geqslant 0$.
By induction, it follows that $\|\mu^{*2^m}*f\|_2\geqslant \|\mu*f\|^{2^m}$, for all $m\geqslant 0$. 
Choose  $m\geqslant 0$  such that $2^m\leqslant n<2^{m+1}$. Then $\|\mu^{*n}*f\|_2\geqslant\|\mu^{*2^{m+1}}*f\|_2\geqslant \|\mu*f\|^{2^{m+1}}\geqslant\|\mu*f\|^{2n},$ as claimed.

(2) Indeed, we have $\mu^{*n}(A)^2=\mu^{*n}(A^{-1})\mu^{*n}(A)\leqslant\mu^{*(2n)}(A^{-1}A)$.
\hfill$\blacksquare$

\subsection{Basic properties of local spectral gap}
We continue with several elementary properties of local spectral gap, starting with an easy, but useful, equivalent formulation of local spectral gap. 

\begin{proposition}\label{kazhdan}
Let $\Gamma\curvearrowright (X,\mu)$ be a measure preserving action of a countable group $\Gamma$, and $B\subset X$ a measurable set of finite measure. Then $\Gamma\curvearrowright (X,\mu)$ has local spectral gap with respect to $B$ if and only if
there exist a finite set
$F\subset\Gamma$ and a constant $\kappa>0$ such that the following holds:
$$\|\xi-\frac{1}{\mu(B)}{\int_{B}\xi\;\text{d}\mu}\|_{2,B}\leqslant\kappa\sum_{g\in F}\|g\cdot\xi-\xi\|_{2,B}\;\;\;\text{for any $\xi\in L^2(G)$}.$$
\end{proposition}

{\it Proof.} The {\it if} implication is clear. To prove the {\it only if} implication, suppose that $\Gamma\curvearrowright (X,\mu)$ has local spectral gap with respect to $B$. Then there are a finite set
$F\subset\Gamma$ and $\kappa>0$ such that 
$\|\eta\|_{2,B}\leqslant\kappa\sum_{g\in F}\|g\cdot\eta-\eta\|_{2,B}$, for any $\eta\in L^2(X)$ with ${\int_{B}\eta\;\text{d}\mu=0}$. We may assume that $e\in F$.
Let $\xi\in L^2(X)$ and put $\alpha={\frac{1}{\mu(B)}\int_{B}\xi\;\text{d}\mu}$. Let $C=\cup_{g\in F}g^{-1}B$ and define $\eta=\xi-\alpha {\bf 1}_C\in L^2(G)$. Then ${\int_B\eta\;\text{d}\mu=0}$ and $\|g\cdot\eta-\eta\|_{2,B}=\|g\cdot\xi-\xi\|_{2,B}$, for all $g\in F$. The conclusion now follows.
\hfill$\blacksquare$

\begin{proposition}\label{indep}
Let $\Gamma\curvearrowright (X,\mu)$ be a measure preserving action of a countable group $\Gamma$, and $B_1,B_2\subset X$ measurable sets of finite measure. Assume there is a finite set $K\subset\Gamma$ such that $B_1\subset\cup_{h\in K}hB_2$ and $B_2\subset\cup_{h\in K}hB_1$. 

Then  $\Gamma\curvearrowright (X,\mu)$ has local spectral gap with respect to $B_1$ if and only if it does with respect to $B_2$.
\end{proposition}

{\it Proof.} Assume that local spectral gap  holds with respect to $B_1$, but not $B_2$. Let $\xi_n\in L^2(X)$ be a sequence satisfying $\|\xi_n\|_{2,B_2}=1$, ${\int_{B_2}\xi_n\;\text{d}\mu=0}$, for all $n$, and $\|g\cdot\xi_n-\xi_n\|_{2,B_2}\rightarrow 0$, for every $g\in\Gamma$. 

If $g\in\Gamma$, then we have \begin{align*}\|g\cdot\xi_n-\xi_n\|_{2,B_1}&\leqslant\sum_{h\in K}\|g\cdot\xi_n-\xi_n\|_{2,hB_2}=\sum_{h\in K}\|(h^{-1}g)\cdot\xi_n-h^{-1}\cdot\xi_n\|_{2,B_2}\\&\leqslant\sum_{h\in K}\Big(\|(h^{-1}g)\cdot\xi_n-\xi_n\|_{2,B_2}+\|h^{-1}\cdot\xi_n-\xi_n\|_{2,B_2}\Big).\end{align*}

This implies that $\|g\cdot\xi_n-\xi_n\|_{2,B_1}\rightarrow 0$, for every $g\in\Gamma$.  Since  we have local spectral gap with respect to $B_1$, Proposition \ref{kazhdan} provides scalars $\alpha_n\in\mathbb C$ such that $\|\xi_n-\alpha_n\|_{2,B_1}\rightarrow 0$. By reasoning as above, it follows that $\|\xi_n-\alpha_n\|_{2,B_2}\rightarrow 0$. Since ${\int_{B_2}\xi_n\;\text{d}\mu=0}$, for all $n$, we get that $\alpha_n\rightarrow 0$. Hence, $\|\xi_n\|_{2,B_2}\rightarrow 0$, which gives the desired contradiction. \hfill$\blacksquare$

Next, we establish that local spectral gap passes to direct product actions.

\begin{proposition}\label{products} For $i\in\{1,2\}$, let $\Gamma_i\curvearrowright (X_i,\mu_i)$ be a measure preserving action which has local spectral gap with respect to a measurable set $B_i\subset X_i$ of finite measure.

Then the product action $\Gamma_1\times\Gamma_2\curvearrowright (X_1\times X_2,\mu_1\times\mu_2)$ has local spectral gap with respect to $B_1\times B_2$.
\end{proposition}

{\it Proof.} 
By Lemma \ref{kazhdan}, for  $i\in\{1,2\}$, we can find $F_i\subset\Gamma_i$ finite and $\kappa_i>0$ such that
\begin{equation}\label{kazhd}\|\xi-\frac{1}{\mu_i(B_i)}{\int_{B_i}\xi\;\text{d}\mu_i}\|_{2,B_i}^2\leqslant\kappa_i\sum_{g\in F_i}\|g\cdot\xi-\xi\|_{2,B_i}^2\;\;\;\text{for any $\xi\in L^2(X_i)$}.\end{equation}

Denote $(X,\mu)=(X_1\times X_2,\mu_1\times\mu_2)$ and $B=B_1\times B_2$. Let $\xi\in L^2(X,\mu)$ and put $\alpha={\frac{1}{\mu(B)}\int_B\xi\;\text{d}\mu}$. For $y\in X_2$, define  $\xi^y(x)=\xi(x,y)$ and $f(y)={\frac{1}{\mu_1(B_1)}\int_{B_1}\xi^y\;\text{d}\mu_1}$. 
Then it is easy to see that $f\in L^2(X_2)$  and $\|g\cdot f-f\|_{2,B_2}^2\leqslant{\frac{1}{\mu_1(B_1)}}\|g\cdot\xi-\xi\|_{2,B_1\times B_2}^2$, for all $g\in\Gamma_2$.
Since ${\frac{1}{\mu_2(B_2)}\int_{B_2}f\;\text{d}\mu_2=\alpha}$, by using the last inequality and applying \eqref{kazhd}  to $f$ we get that 

\begin{equation}\label{kazhd2}\|f-\alpha\|_{2,B_2}^2\leqslant\kappa_2\sum_{g\in F_2}\|g\cdot f-f\|_{2,B_2}^2\leqslant\frac{\kappa_2}{\mu_1(B_1)}\sum_{g\in F_2}\|g\cdot\xi-\xi\|_{2,B_1\times B_2}^2 \end{equation}

On the other hand, by applying \eqref{kazhd} to $\xi^y$ we get that ${\|\xi^y-f(y)\|_{2,B_1}^2\leqslant\kappa_1\sum_{g\in F_1}\|g\cdot\xi^y-\xi^y\|_{2,B_1}^2}$. By integrating over $y\in B_2$, we derive that 
\begin{equation}\label{kazhd3} \int_{B}|\xi(x,y)-f(y)|^2\;\text{d}\mu(x,y)\leqslant\kappa_1\sum_{g\in F_1}\|g\cdot\xi-\xi\|_{2,B_1\times B_2}^2.
\end{equation}

It is now clear that the combination of \eqref{kazhd2} and \eqref{kazhd3} implies the conclusion.
\hfill$\blacksquare$

Finally, we record a result asserting that local spectral gap passes through certain quotients. Since its proof is very similar to that of Corollary \ref{by}, we leave its details to the reader.

\begin{proposition}
Let $G$ be a l.c.s.c. group, $H<G$ a closed subgroup, and $\Gamma<G$ a countable dense subgroup. Assume that $G/H$ admits a $G$-invariant Borel regular measure $m_{G/H}$. Suppose that the left translation action $\Gamma\curvearrowright (G,m_G)$ has local spectral gap.

Then the left translation action $\Gamma\curvearrowright (G/H,m_{G/H})$ has local spectral gap.
\end{proposition}

\subsection{Deduction of Theorem \ref{main} from Theorem \ref{restricted}} The aim of this subsection is to show that Theorem \ref{restricted} implies Theorem \ref{main}. This relies on the following result.

\begin{proposition}\label{AtoB}
Let $G$ be a l.c.s.c. group, $\Gamma<G$ a countable dense subgroup, and $B\subset G$ a measurable set with non-empty interior and compact closure. 
Assume that there exists a constant $c > 0$ satisfying the following property: for  any neighborhood $U$ of the identity, there are a finite set $S \subset \Gamma \cap U$ and a finite dimensional vector space $V \subset L^2(G)$ such that for all $\xi\in L^2(B)\ominus V$ we have \[\max_{g \in S} \Vert g\cdot\xi - \xi \Vert_{2} \geqslant c\Vert \xi \Vert_{2}.\]

Then the left translation action $\Gamma\curvearrowright (G,m_G)$ has local spectral gap with respect to $B$.

\end{proposition}
{\it Proof.}  Assume by contradiction that the conclusion is false. Then there is a sequence $\xi_n\in L^2(G)$ satisfying
$\Vert \xi_n \Vert_{2,B} = 1$, ${\int_B\xi_n\;\text{d}\mu = 0}$, for all $n$, and  $\lim\limits_{n\rightarrow\infty} \Vert g \cdot \xi_n - \xi_n \Vert_{2,B} = 0$, for all $g \in \Gamma$.

 If $C\subset G$ is a compact set, then $C$ can be covered with finitely many of the sets $\{gB\}_{g\in\Gamma}$. It follows that $\sup_{n}\|\xi_n\|_{2,C}<\infty$ and $\lim\limits_{n\rightarrow\infty}\|g\cdot\xi_n-\xi_n\|_{2,C}=0$, for all $g\in\Gamma$. Since $G$ is second countable, we can find a subsequence $\{\xi_{n_k}\}$ of $\{\xi_n\}$ and $\xi\in L^2_{\text{loc}}(G)$  (i.e. a locally $L^2$-integrable function) such that ${\bf 1}_C\xi_{n_k}\rightarrow {\bf 1}_C\xi$, weakly, for every compact  set $C\subset G$. But then $\xi$ must be $\Gamma$-invariant, and hence constant by ergodicity. Since $\xi_n$ has mean zero on $B$, for all $n$, we derive that $\xi=0$, almost everywhere. This argument implies that ${\bf 1}_C\xi_n\rightarrow 0$, weakly, for any compact set $C\subset G$.
 
Let  $\lim\limits_n$  be a bounded linear functional on $\ell^{\infty}(\mathbb N)$ which extends the  limit.
Then $\nu(C)=\lim\limits_n \Vert {\bf 1}_C\xi_n \Vert_2^2$ defines a $\Gamma$-invariant finitely additive measure on bounded Haar measurable subsets of $G$. Since $\nu(B) \neq 0$, we get that $\nu\not=0$. 
Since $B$ has non-empty interior, by using finite additivity, we can find two open sets $B_1  \subsetneq B_2\subset B$ such that $\nu(B_1) \neq 0$, $\nu(B_2 \setminus B_1) \leqslant (c^2 \nu(B_1))/4$, and there exists a closed intermediate subset $B_1 \subset F \subset B_2$.

By local compactness, we can find an intermediate open set $B_0$ between $B_1$ and $B_2$ and a neighborhood $U$ of the identity small enough so that $B_1 \subset gB_0 \subset B_2$, for $g\in U$. Put $p = \mathbf{1}_{B_0}$.
As the sequence $\{p\xi_n\}$ converges weakly to $0$ and is supported on $B$, the following claim contradicts our assumption on $c$.

{\bf Claim.} For all $g \in \Gamma \cap U$, we have $\lim\limits_n \Vert g \cdot (p\xi_n) - (p\xi_n)\Vert_2 \leqslant \frac{c}{2} \lim\limits_n \Vert p\xi_n\Vert_2$. 

Indeed, for $g \in \Gamma \cap U$ we can estimate
\begin{align*}
\lim_n \Vert g \cdot (p\xi_n) - (p\xi_n)\Vert_2 & = \lim_n \Vert (g \cdot p)(g \cdot \xi_n) - (p\xi_n)\Vert_2\\
& \leqslant \lim_n \Vert (g \cdot p)\xi_n - (p\xi_n)\Vert_2 + \lim_n \Vert (g \cdot p)(g \cdot \xi_n) - (g \cdot p)\xi_n\Vert_2\\
& = \lim_n \Vert (g \cdot p - p)\xi_n\Vert_2 + \lim_n \Vert g \cdot \xi_n - \xi_n\Vert_{2,gB_0}
\end{align*}
But by the above, $\lim\limits_n \Vert g \cdot \xi_n - \xi_n\Vert_{2,gB_0} = 0$. Moreover, since $B_1 \subset gB_0 \subset B_2$, we get that
\[\lim_n \Vert (g \cdot p - p)\xi_n\Vert_2^2 \leqslant \nu(B_2 \setminus B_1) \leqslant \frac{c^2}{4}\nu(B_1) \leqslant \frac{c^2}{4} \nu(B_r(1))  = \frac{c^2}{4} \lim_n \Vert p\xi_n\Vert_2^2 \qedhere\]
\hfill$\blacksquare$
 
 {\bf Proof of Theorem \ref{main}}. Assume that Theorem \ref{restricted} holds and let us explain how Theorem \ref{main} follows. Let $S\subset\Gamma$ be a finite set and denote $\mu={\frac{1}{|S|}\sum_{g\in S}\delta_g}$.  If $\xi\in L^2(G)$, then  $$\sum_{g\in S}\|g\cdot\xi-\xi\|_2^2=2|S|\Big(\|\xi\|_2^2-\Re\langle\mu*\xi,\xi\rangle\Big).$$  Thus, if we have that $\|\mu*\xi\|<\frac{1}{2}\|\xi\|_2$, then $\max_{g\in S}\|g\cdot\xi-\xi\|_2>\|\xi\|_2$.
 By combining  Theorem \ref{restricted} and Proposition \ref{AtoB}, we conclude that $\Gamma\curvearrowright (G,m_G)$ has local spectral gap. \hfill$\blacksquare$
 
\subsection{Reduction to groups with trivial center} Next, we will argue that in order to prove Theorem \ref{restricted}, we may reduce to the case when $G$ has trivial center. 

Assume that Theorem \ref{restricted} holds for connected simple Lie groups with trivial center. 
Let $G$ be a connected simple Lie group, $B\subset G$ a measurable set with compact closure and non-empty interior,  and $c>0$.
Let $\pi:G\rightarrow\GL(\frak g)$ be the adjoint representation of $G$.  Put $G_0=\pi(G)$ and $\Gamma_0=\pi(\Gamma)$.

Since $\pi$ has discrete kernel, we can find  a small enough  compact neighborhood of the identity $C\subset G$ and $\varepsilon_0>0$ such that  $\pi$ is 1-1 on $\cup_{g\in B_{\varepsilon_0}(1)}gC$. Let $K\subset G$ be a finite set such that $B\subset\cup_{h\in K}Ch$. Write $B$ as a disjoint union $B=\sqcup_{h\in K}C_h$, where $C_h$ is a subset of $Ch$, for every $h\in K$.

Since $G_0$ has trivial center, the conclusion of Theorem \ref{restricted} holds for $(G_0,\Gamma_0,\pi(C))$ by our assumption. It is then easy to see that Theorem \ref{restricted} also holds for $(G,\Gamma, C)$. 

Thus, given $\varepsilon>0$, there are a finite set $T\subset\Gamma\cap B_{\varepsilon}(1)$ and a finite dimensional subspace $W\subset L^2(C)$ such that $\mu:=\frac{1}{2|T|}\sum_{g\in T}(\delta_g+\delta_{g^{-1}})$ satisfies  $\|\mu*F\|_2<\frac{c}{|K|}\|F\|_2$, for every $F\in L^2(C)\ominus W.$
 
 Let $V\subset L^2(B)$ be the linear span of $\{{\bf 1}_{C_h}\rho_h(W)|h\in K\}$, where $\{\rho_g\}_{g\in G}$ denotes the right regular representation of $G$. 
 Let $F\in L^2(B)\ominus V$. If $h\in K$, then $\rho_h^{-1}({\bf 1}_{C_h}F)\in L^2(C)\ominus W$, and therefore $$\|\mu*({\bf 1}_{C_h}F)\|_2=\|\mu*(\rho_h^{-1}({\bf 1}_{C_h}F))\|_2<\frac{c}{|K|}\|\rho_h^{-1}({\bf 1}_{C_h}F)\|_2=\frac{c}{|K|}\|{\bf 1}_{C_h}F\|_2.$$
 
 Since $F=\sum_{h\in K}{\bf 1}_{C_h}F$, we deduce that $\|\mu*F\|_2<c\|F\|_2$. Since $c>0$ is arbitrary, and $V\subset L^2(B)$ is finite dimensional, this implies that $G$ satisfies the conclusion of Theorem \ref{restricted}.

\subsection{Outline of the proof of Theorem \ref{restricted}}  
\label{outline}

The previous subsection allows us to work hereafter with connected simple groups with trivial center. 
Note, however, that although our results will be stated only for groups with trivial center, they have analogues
  for general connected simple groups.

 In order to outline the proof of Theorem \ref{restricted}, let us introduce some more notation.
Let $G$ be a connected simple Lie group with trivial center, and $\Gamma<G$ be a dense subgroup. Suppose there is a basis $\frak B$ of $\frak g$ such that the matrix of Ad$(g)$ in the basis $\frak B$ has algebraic entries, for any $g\in\Gamma$.  
 Let $n$ be the dimension of $G$, $\frak g$  its Lie algebra, and Ad$:G\rightarrow GL(\frak g)$ its adjoint representation. We identify $G\cong$ Ad$(G)$, $\frak g\cong\mathbb R^n$ via the basis $\frak B$, and $GL(\frak g)\cong GL_n(\mathbb R)\subset\mathbb M_n(\mathbb R)$. In particular, in this identification we have that $\Gamma<GL_n(\bar{\mathbb Q})$.

For  $\alpha=(\alpha_{i,j})_{i,j=1}^n\in\mathbb M_n(\mathbb R)$, we denote by $\|\alpha\|_2={(\sum_{i,j=1}^n|\alpha_{i,j}|^2)^{1/2}}$ its Hilbert-Schmidt norm. We endow $G$ with the metric given by $(g,h)\mapsto\|\text{Ad}(g)-\text{Ad}(h)\|_2$.  Abusing notation, we write $\|g-h\|_2:=\|\text{Ad}(g)-\text{Ad}(h)\|_2$ and $\|g\|_2:=\|\text{Ad}(g)\|_2$. 
Note that $$\|gh-gk\|_2\leqslant \|g\|_2\|h-k\|_2,\;\;\;\;\;\text{for all $g,h,k\in G$}.$$  

For $x\in G$ and $\delta>0$, we denote $B_{\delta}(x):=\{y\in G|\|x-y\|_2\leqslant\delta\}$.  For $\delta> 0$, we let $A^{(\delta)}=\cup_{x\in A}B_{\delta}(x)$ be the {\it $\delta$-neighborhood} of $A\subset G$, and denote $$P_{\delta}:=\frac{{\bf 1}_{B_{\delta}(1)}}{|B_{\delta}(1)|}\in L^1(G)_{+,1}.$$

As explained at the end of the Introduction, the proof of Theorem \ref{restricted} splits into three parts, dealt with in the following three sections.

\begin{itemize}
\item In Section 3 we produce measures with small support that {\bf Escape subgroups} quickly.
\end{itemize}
There are two steps for this. First, we produce for all $\varepsilon > 0$ a finite set $S\subset\Gamma\cap B_{\varepsilon}(1)$ and constants $d,C > 0$ such that for $\delta>0$ small enough, the measure $\mu_S :=\frac{1}{2|S|}\sum_{g\in S}(\delta_g+\delta_{g^{-1}})$ satisfies
$\mu_S^{*n}(H^{(\delta)})\leqslant\delta^{d}$, for all proper closed subgroups $H$ and $n\approx C\log{\frac{1}{\delta}}.$
This step is obtained by combining Propositions \ref{ping-pong} and \ref{neigh}. The set $S$ that we obtain freely generates a free group.

One can of course get a better constant $C$ by modifying accordingly the value of $d$. But for a fixed $d$, the value of $C$ depends on $\varepsilon$. Namely, it could happen that $C \to \infty$ as $\varepsilon \to 0$.
As explained in the introduction, we want to  control the speed of escape in terms of $\varepsilon$.
So the second step  is to upgrade the set $S$ to a set $T$, also contained  in $B_{\varepsilon}(1)$, such that the following holds (Theorem \ref{escape}).

{\it There are constants $d_1,d_2>0$ not depending on $T$ such that the probability measure $\mu_T$ satisfies
$\mu_T^{*n}(H^{(\delta)})\leqslant\delta^{d_1}$, for all $\delta > 0$ small enough, all proper closed subgroups $H$ and $n \approx d_2\frac{\log{\frac{1}{\delta}}}{\log{\frac{1}{\varepsilon}}}.$}

This improvement is obtained using the pigeonhole principle and the freeness of the elements of $S$.

\begin{itemize}
\item In Section 4 we extend the {\bf $\ell^2$-flattening} lemma from \cite{BdS14}.
\end{itemize}
Our generalization of the flattening lemma \cite[Lemma 2.5]{BdS14} to the locally compact setting does not require much additional effort. However,  it only applies for measures with controlled support. But we anticipated this issue in part 1 above, by controlling the speed of escape in terms of $\varepsilon$. Indeed, we want to apply the flattening lemma to the measure $\mu_T^{\ast n}$, with $n \approx d_2\frac{\log{\frac{1}{\delta}}}{\log{\frac{1}{\varepsilon}}}$. Now, the support of $\mu_T^{\ast n}$  is contained in  $B_{\delta^{-\beta}}(1)$, with $\beta>0$ arbitrarily small. Since the ``controlled support'' condition that we require is soft enough, we are in position to apply our flattening Lemma \ref{BdS}.

Thus, our main result (Corollary \ref{flattening}) shows that the measure $\mu_T$ produced in Section 3 will flatten rather quickly:
given $\alpha>0$, we have $\|\mu_{T}^{*n}*P_{\delta}\|_2\leqslant\delta^{-\alpha}$, for $\delta$ small enough and $n\sim\log\frac{1}{\delta}$.

\begin{itemize}
\item In Section 5 we prove a {\bf Mixing inequality} and combine it with the above to conclude.
\end{itemize}

More precisely, we show that if  $\mu_T$ is the measure produced in Section 3, then the convolution operator $F \in L^2(B) \mapsto (\mu_T  \ast F) \in L^2(G)$ has norm less than $1/2$, when restricted to the orthogonal complement of a finite dimensional subspace $V \subset L^2(B)$. The first observation  is that this flexibility of discarding a finite dimensional subspace $V$ when trying to bound the norm of $\Vert \mu_T \ast F\Vert_2$, allows us to restrict our study to functions $F$ that live at a ``small scale''. Namely, it will be enough to consider functions $F$ that do not change much when ``discretizing'' the group with high accuracy. This reduction is achieved via a Littlewood-Paley type decomposition (Theorem \ref{L-P} and Corollary \ref{level_delta}). Then we are left to show a mixing inequality (Theorem \ref{rho}). This is inspired by \cite[Lemma 10.35]{BG10}, and should be thought of as an analogue of the well-known mixing inequality for finite groups (see e.g. \cite[Proposition 1.3.7]{Ta15}), after discretizing the group. We will then be able to conclude restricted spectral gap by combining this inequality with the flattening obtained in Section 4.



\section{Escape from  subgroups}

The goal of this section is to prove the following:

\begin{theorem}[escape from subgroups]\label{escape}
Let $G$ be a connected simple Lie group with trivial center, and
\text{Ad}$:G\rightarrow GL(\frak g)$ be its adjoint representation. Let $\Gamma<G$ be a countable dense subgroup. Assume that there there is a basis $\frak B$ of $\frak g$ such that the matrix of Ad$(g)$ in the basis $\frak B$ has algebraic entries, for every $g\in\Gamma$. 

Then there are constants $d_1,d_2 > 0$ depending on $\Gamma$ only such that the following holds.

Given $\varepsilon_1>0$, we can find $0<\varepsilon<\varepsilon_1$ and a finite set $T\subset \Gamma\cap B_{\varepsilon}(1)$ which freely generates a subgroup of $\Gamma$ such that for any small enough $\delta>0$, the probability measure $\mu={\frac{1}{2|T|}\sum_{g\in T}(\delta_g+\delta_{g^{-1}})}$ satisfies $$\mu^{*2n}(H^{(\delta)})\leqslant\delta^{d_1},\;\;\;\text{where $n=\Big\lfloor{d_2\frac{\log{\frac{1}{\delta}}}{\log{\frac{1}{\varepsilon}}}}\Big\rfloor$},$$ for any proper closed connected subgroup $H<G$.

\end{theorem}

\subsection{Ping-pong} The first ingredient in the proof of Theorem \ref{escape} is a proposition which, roughly speaking, asserts the existence of representations $\rho_i:\Gamma\rightarrow$ GL$(V_i)$, $i\in I$, and $M\geqslant 2$ such that 
\begin{itemize}
\item the intersection of $\Gamma$ with any proper closed subgroup of $G$ stabilizes a line in some $V_i$, and 
\item we can find a set $S\subset\Gamma$ of simultaneous ``ping-pong players"  for all the $\rho_i$'s in any given neighborhood of the identity in $G$ such that $|S|=M$.
\end{itemize}

\begin{proposition}\label{ping-pong}
Let $G$ be a connected simple real Lie group with trivial center. Let $\Gamma<G$ be a finitely generated dense subgroup. Assume that there there is a basis $\frak B$ of the Lie algebra $\frak g$ of $G$ such that the matrix of Ad$(g)$ in the basis $\frak B$ has algebraic entries, for every $g\in\Gamma$. 

Then there exist finitely many vector spaces $V_i$, $i\in I$, defined over local fields $K_i$,  representations $\rho_i: \Gamma \to \GL(V_i)$, and an integer $M\geqslant 2$ such that the following properties hold true:
\begin{enumerate}
\item For any proper closed subgroup $H < G$ such that $\Gamma \cap H$ is non-discrete, there exist $i\in I$ and $[v] \in \P(V_i)$ such that $\rho_i(g)([v])=[v]$, for all $g\in \Gamma\cap H$.
\item For any $\eta>0$, there is a finite set $S\subset\Gamma$ satisfying $|S|=M$ and $S\subset B_{\eta}(1)$ such that for all $i\in I$ and every $g\in\tilde S:=S\cup S^{-1}$, we can find two sets $K_g^{(i)}\subset U_g^{(i)}\subset \P(V_i)$ such that the following conditions hold:
\begin{enumerate}[(a)]
\item For every  $g\in\tilde S$ we have $\rho_i(g)(U_g^{(i)})\subset K_g^{(i)}$.
\item Every line $[v]\in\P(V_i)$ is contained in at least two of the sets $\{U_g^{(i)}\}_{g\in\tilde S}$.
\item For every $g_1,g_2\in\tilde S$ we have $K_{g_1}^{(i)}\subset U_{g_2}^{(i)}$, unless $g_1g_2=1$.
\item For every $g_1,g_2\in\tilde S$ we have  $K_{g_1}^{(i)}\cap K_{g_2}^{(i)}=\emptyset$, unless $g_1=g_2$.
\end{enumerate}
\end{enumerate}
\end{proposition}

Before proving Proposition \ref{ping-pong}, let us record a simple observation that will be used later.

\begin{lemma}\emph{\cite{SGV11}}\label{SGV}
In the setting from Proposition \ref{ping-pong}, let $i\in I$ and $v\in V_i\setminus\{0\}$. 
Let $g=g_ng_{n-1}...g_1$ be a reduced word on length $n$ in $\tilde S$. Assume that $\rho_i(g)([v])=[v]$ and let $1\leqslant j<n$.

\begin{enumerate}
\item If $\rho_i(g_jg_{j-1}...g_1)([v])\in U_{g_{j+1}}$, then $g_{j+1},...,g_n$ are uniquely determined by $v$.

\item If $\rho_i(g_jg_{j-1}...g_1)([v])\notin U_{g_{j+1}}$, then $\rho_i(g_lg_{l-1}...g_1)([v])\notin U_{g_{l+1}}$, for all $1\leqslant l\leqslant j$.
\end{enumerate}
\end{lemma}

{\it Proof.} For simplicity, denote $\rho=\rho_i$ and $K_g=K_g^{(i)}$, $U_g=U_g^{(i)}$, for all $g\in\tilde S$. 
Assume that $\rho(g_jg_{j-1}...g_1)([v])\in U_{g_{j+1}}$.
Since $\rho(g_{j+1})(U_{g_{j+1}})\subset K_{g_{j+1}}$, we get $\rho(g_{j+1}g_j...g_1)([v])\in K_{g_{j+1}}$. Since $g_{j+1}\not=g_{j+2}^{-1}$, we have that $K_{g_{j+1}}\subset U_{g_{j+2}}$, and hence  $\rho(g_{j+1}g_j...g_1)([v])\in U_{g_{j+2}}$. Using induction it follows that $\rho(g_pg_{p-1}...g_1)([v])\in K_{g_p}$, for all $j+1\leqslant p\leqslant n$. Thus, $[v]=\rho(g_n...g_1)([v])\in K_{g_n}$. Since the sets $\{K_g\}_{g\in\tilde S}$ are mutually disjoint, $g_n$ is therefore determined by $v$. Further, we have that $\rho(g_n^{-1})([v])=\rho(g_{n-1}...g_1)\in K_{g_{n-1}}$. Since $\rho(g_n^{-1})([v])$ is determined by $v$, we deduce that $g_{n-1}$ is also determined by $v$. The first assertion now follows by induction.
Since the beginning of the proof implies the second assertion, the proof is complete.
 \hfill$\blacksquare$

The rest of this subsection is devoted to proving Proposition \ref{ping-pong}. The proof is very similar to the the proof of \cite[Proposition 21]{SGV11}. Consider a connected simple real Lie group $G$ with trivial center, together with a finitely generated dense subgroup $\Gamma$ as in the statement of Proposition \ref{ping-pong}.

By identifying $G$ with $\Ad(G)$, we can assume that $G$ is the connected component of a real algebraic group $\G = \overline{\Ad(G)}^Z \subset \GL(\fg)$. By our assumptions on $\Gamma$ we can find a number field $k$ with an embedding $k \subset \R$ and a basis of $\fg$ such that $\Gamma \subset \GL_d(k) \subset \GL_d(\R) \cong \GL(\fg)$. Since $\Gamma$ is Zariski dense in $\G$, we see that $\G$ is in fact defined over $k$.

Now, note that in order to check item (1) of Proposition \ref{ping-pong} for a subgroup $H < G$, it suffices to check it for the closure of $\Gamma\cap H$ (in the real topology). This shows that we only have to deal with proper closed subgroups $H < G$ which are non-discrete in $G$ and such that $\Gamma \cap H$ is dense in $H$. But if $H$ is such a subgroup, then its Zariski closure $\bbh \subset \G$ is a proper algebraic subgroup which is defined over $k$, because $\Gamma \cap H \subset \GL_d(k)$ is a Zariski dense subgroup of it. Hence the Lie algebra $\fh \subset \fg$ of $\bbh$ is a non-trivial proper subspace of $\fg$ defined over $k$ which is globally invariant under $\bbh$ but not under $\G$, because $G$ is simple. Altogether, we find that the line in 
$\bigwedge_{j=1}^{\dim \fh} \fg(k)$ corresponding to the subspace $\fh(k) \subset \fg(k)$ is invariant under $\bbh(k)$, but not under $\G(k)$.

Next, define the finite set of representations $\rho_i$, $i \in I$, to be the collection of all (non-trivial) irreducible subrepresentations of the representations of $\G$ on $\bigwedge_{j=1}^m \fg$, $m < d = \dim G$. These are algebraic representations, defined over a finite extension $k'$ of $k$. We will show that Proposition \ref{ping-pong} holds when we view these representations $\rho_i$ as defined over appropriate places $K_i$ of $k'$. Note that no matter how we choose the places $K_i$, we still get representations of $\Gamma \subset \G(k)$ which satisfy  item (1) of the proposition, by the above paragraph.
Let us now choose the places $K_i$ for which we will be able to prove that item (2) of the proposition also holds true.

\begin{lemma}
\label{choosefields}
Use the above notation. 
Then for every $i \in I$, there are a local field $K_i$ and a sequence $(h_n)_n$ in $\Gamma$ which converges to $1$ in the real topology such that $(\rho_i(h_n))_n$ goes to infinity in the $K_i$-topology. 
\end{lemma}
{\it Proof.}
Fix $i \in I$. First we claim that there is a sequence $(h_n)_n \subset \Gamma$ which converges to $1$ in the real topology such that the elements $\rho_i(h_n)$ are pairwise distinct. 

To prove the claim, we view the representations $\rho_i$ as representations over $\C$ by fixing an embedding $k' \subset \C$. This way, it makes sense to talk about $\rho_i(G)$. Note that the image $\rho_i(B_1(1))$ of the unit ball of $G$ is connected. Since the representation $\rho_i$ is non-trivial and $\Gamma$ is dense in $G$, there exists a sequence $(h_n)_n \subset \Gamma \cap B_1(1)$ such that $\rho_i(h_n)$ is non-trivial and converges to $1$ in the complex topology. In particular, after passing to a subsequence of $h_n$ if necessary, we get that the elements $\rho_i(h_n)$ are distinct. 
As $(h_n)_n$ is a bounded sequence and $\rho_i(h_n)$ converges to $1$, we deduce that $h_n$ converges to $1$ as well, proving our claim.

Next, denote by $R$ the ring generated by the coefficients of the elements $\Ad(g)$, $g \in \Gamma$. Since $\Gamma$ is finitely generated, $R\subset k'$ is a finitely generated subring. The discrete diagonal embedding of $k'$ in its ad\`ele group gives a discrete embedding of $R$ in a product of finitely many places $K_\nu$, $\nu \in \cS$, of $k'$.

From this we obtain a discrete embedding ${\rho_i(\G(R)) \hookrightarrow \Pi_{\nu \in \cS} \rho_i(\G(K_\nu))}$. In particular, $\rho_i(\Gamma)$ is discrete inside $\Pi_{\nu \in \cS} \rho_i(\G(K_\nu))$. Therefore there exists a field $K_i := K_\nu$ such that the infinite set $\{ \rho_i(h_n) \}$ is unbounded as a subset of $\rho_i(\G(K_\nu))$.
\hfill$\blacksquare$

Below, we denote by $\oneG$ the set of sequences $(h_n) \subset \Gamma$ which converge to $1$ in the real topology. For $i\in I$,  we view $\rho_i: \G \to \GL(V_i)$ as a representation over $K_i$, and equip $\GL(V_i)$ with the operator norm $\Vert \cdot \Vert_i$ corresponding to the absolute value on $K_i$.
 We denote by $A_i$ the set of cluster points in the $K_i$-topology of sequences of the form $(\rho_i(h_n)/\Vert \rho_i(h_n) \Vert_i)_n$, where $(h_n) \in\oneG$. 
Finally, we put $r_i = \min_{b \in A_i} \rk(b)$, where $\rk(b)$ is the rank of $b$. 

A key fact that we will use is that since $G$ is simple, we have that $\rho_i(g)$ has determinant $1$ for all $g \in G$, and in particular for all $g \in \Gamma$. Hence, if $(h_n)_n\subset\Gamma$ and $(\rho_i(h_n))_n$ is unbounded in the $K_i$-topology, then the normalized sequence $(\rho_i(h_n)/\Vert \rho_i(h_n) \Vert_{K_i})_n$ has a non-invertible cluster point. So by our choice of $K_i$, we have $r_i < d_i := \dim(V_i)$.

Let us mention the following stability result for the sets $A_i$.

\begin{lemma}
\label{stability}
If $b$ and $b'$ belong to $A_i$ and $bb' \neq 0$, then some scalar multiple of $bb'$ belongs to $A_i$. In particular $\rk(bb') \geq r_i$.
\end{lemma}
{\it Proof.}
If $b = \lim_n \rho_i(g_n)/\Vert \rho_i(g_n) \Vert_i$ and $b' = \lim_n \rho_i(h_n)/\Vert \rho_i(h_n) \Vert_i$, with $(g_n)_n, (h_n)_n \in \oneG$, then the product sequence $(g_nh_n)_n \subset \Gamma$ converges to 1 in the real topology. Moreover, 
\[\lim_n \frac{\rho_i(g_nh_n)}{\Vert \rho_i(g_n) \Vert_i \Vert \rho_i(h_n) \Vert_i} = bb', \qquad  \text{so that} \qquad \lim_n \frac{\Vert  \rho_i(g_nh_n)\Vert_i}{\Vert \rho_i(g_n) \Vert_i \Vert \rho_i(h_n) \Vert_i} =\Vert bb' \Vert_i.\]
Therefore
\[ \lim_n \frac{\rho_i(g_nh_n)}{\Vert \rho_i(g_nh_n) \Vert_i} = \frac{bb'}{\Vert bb' \Vert_i}. \qedhere \]
\hfill$\blacksquare$


Now, we turn to the construction of the set $S$ from Proposition \ref{escape}. The following lemma will produce the first element of $S$. In the context of Lemma \ref{firstplayer}, it will be $g_n$ with $n$ large enough, depending on $\eta$. The other elements of $S$ will arise as appropriate conjugates of this first element.

\begin{lemma}\label{firstplayer}
There exists a sequence $(g_n)_n \in \oneG$ such that for all $i \in I$
\begin{enumerate}[(1)]
\item $\lim_n \frac{\rho_i(g_n)}{\Vert \rho_i(g_n)\Vert_{i}} = b_i$ for some $b_i$ with $\rk(b_i) = r_i$;
\item $\Range(b_i) \cap \Ker(b_i) = \{0\}$.
\end{enumerate}
\end{lemma}
{\it Proof.} We proceed in three steps.

{\bf Step 1.} There exists a sequence $(h_n)_n \in \oneG$ which satisfies (1) above for all $i \in I$.

We proceed by induction. Enumerate the set $I = \{1,\cdots,\vert I \vert \}$. Assume that $(k_n)_n \in \oneG$ satisfies (1) for all indices $i < i_0$, for some $1 \leqslant i_0 < |I|$. Taking a subsequence if necessary, we can assume that the sequence $(\rho_{i_0}(k_n)/\Vert \rho_{i_0}(k_n) \Vert_{i_0})$ converges to some element $b_{i_0} \in \End(V_i)$. 

Since the rank of $b_{i_0}$ could be greater than $r_{i_0}$, we also consider a sequence $(k_n')_n \in \oneG$ such that $(\rho_{i_0}(k_n')/\Vert \rho_{i_0}(k_n') \Vert_{K_{i_0}})_n$ converges to some $b_{i_0}'$ with rank $r_{i_0}$. Taking a subsequence we can assume $(\rho_i(k_n')/\Vert \rho_{i}(k_n') \Vert_{i})_n$ converges to some element $b_i'$ for all $i < i_0$, with possibly $\rk(b_i') > r_i$.

We will prove the existence of an element $g \in \Gamma$ such that the sequence $(gk_ng^{-1}k_n')_n$ satisfies (1) for all $i \leqslant i_0$. Note that no matter how we choose $g$, this sequence is inside $\oneG$. In fact it suffices to find $g \in \Gamma$ such that $\rho_i(g^{-1})b_i\rho_i(g)b_i' \neq 0$, for all $i \leqslant i_0$. Indeed, then Lemma \ref{stability} implies that the sequence $(\rho_i(gk_ng^{-1}k_n')/\Vert \rho_i(gk_ng^{-1}k_n') \Vert_i)_n$ converges to some non-zero multiple of $\rho_i(g^{-1})b_i\rho_i(g)b_i'$, which has rank at most equal to $\min(\rk(b_i),\rk(b_i')) = r_i$.

For each $i \leqslant i_0$, the set $X_i = \{g \in \G(K_i) \, \vert \, \Range(\rho_i(g)b_i') \nsubseteq \Ker(b_i)\}$ is a Zariski open set in $\G$ which is non-empty because $\rho_i$ is irreducible. Therefore $\Gamma \bigcap (\cap_{i \leqslant i_0} X_i)$ is nonempty. This proves Step 1.

{\bf Step 2.} There exists a sequence $(g_n)_n \in \oneG$ such that $(1)$ is true for any $i \in I$ and the corresponding elements $b_i$ satisfy $b_i^2 \neq 0$.

Consider a sequence $(h_n)$ as in Step 1, and denote by $b_i'$ the corresponding elements. We will find an element $g \in \Gamma$ such that the sequence of elements $g_n := gh_ng^{-1}h_n$ does what we want.

As above, for any $i$, the set $X_i$ of elements $g \in \G(K_i)$ such that $\Range(\rho_i(g)b_i') \nsubseteq \Ker(b_i')$ is a non-empty Zariski-open set. So is the set $Y_i$ of $g \in \G(K_i)$ such that $\Range(\rho_i(g^{-1})b_i') \nsubseteq \Ker(b_i')$, for all $i \leqslant |I|$. Take $g \in \Gamma \bigcap (\cap_{i\leqslant |I|} (X_i \cap Y_i))$ so that the element $a_i := \rho_i(g)b_i'\rho_i(g^{-1})b_i'$ is non-zero. Then for all $i$, the sequence $(\rho_i(g_n)/\Vert \rho_i(g_n)\Vert_i)_n$ converges to some nonzero multiple $b_i$ of $a_i$. We claim that $a_i^2 \neq 0$. 

Indeed, Lemma \ref{stability} implies that the rank of $a_i$ is equal to $r_i = \rk(\rho_i(g)b_i')$. This means that the range of $a_i$ is equal to the range of $\rho_i(g)b_i'$. Since  $g \in X_i$, it follows that $b_i'a_i$ is non-zero. Using again Lemma \ref{stability}, we get that the rank of $b_i'a_i$ is equal to $r_i = \rk(b_i')$. This means that the range of $b_i'a_i$ is equal to the range of $b_i'$. But since $g \in Y_i$ we see that $b_i'\rho_i(g^{-1})b_i'a_i \neq 0$. This shows that $a_i^2 \neq 0$.

{\bf Step 3.} The sequence from Step 2 satisfies the conclusion of the lemma.

We just need to check that for all $i$, any element $b \in A_i$ with rank $r_i$ and such that $b^2 \neq 0$ satisfies $\Range(b) \cap \Ker(b) = \{0\}$. Indeed, if $b^2 \neq 0$ then Lemma \ref{stability} implies that some multiple of $b^2$ belongs to $A_i$. Hence $\rk(b^2) = r_i = \rk(b)$. This precisely means that $\Range(b) \cap \Ker(b) = \{0\}$.
\hfill$\blacksquare$

Before actually proving Proposition \ref{ping-pong}, let us give two easy lemmas.

\begin{lemma}\label{inverse}
Given a local field $\cK$, consider a sequence of invertible elements  $(g_n)_n \subset \GL_d(\cK)$, such that 
\[ \lim_n \frac{g_n}{\Vert g_n \Vert} = b \qquad \text{and} \qquad \lim_n \frac{g_n^{-1}}{\Vert g_n^{-1} \Vert} = b',\]
for some non-invertible elements $b,b' \in M_d(\cK)$.
Then $bb' = 0$, so that $\Range(b') \subset \Ker(b)$.
\end{lemma}
{\it Proof.}
Note that $bb'$ is a scalar matrix, being the limit of the sequence $(1/ \Vert g_n \Vert \Vert g_n^{-1} \Vert)_n$. Since it is non-invertible, it must be $0$.
\hfill$\blacksquare$

\begin{lemma}\label{conjugators}
Let $\rho:\G(\cK)\rightarrow\GL(W_{\rho})$ be an irreducible algebraic representation over a local field $\cK$. Let $V_1^+,V_1^-,V_2^+,V_2^-\subseteq W_\rho$ be non-zero, proper subspaces such that $V_1^+\cap V_1^-=V_2^+\cap V_2^-=\{0\}$ and $V_1^+\subseteq V_2^-$ and $V_2^+\subseteq V_1^-$.
For $M \geqslant 1$, denote by $X_M \subseteq \G(\cK)^M$ the set of $M$-tuples $(h_1,\ldots,h_M)$ satisfying the following two conditions.
\begin{enumerate}
\item For $1 \leqslant s \neq t \leqslant M$, we have $\rho(h_s)V_1^+ \nsubseteq \rho(h_t)(V_1^{-} \cup V_2^-)$ and $\rho(h_s)V_2^+ \nsubseteq \rho(h_t)(V_1^{-} \cup V_2^-)$;
\item For any subset $S \subset \{1,\dots,M\}$ and any choice of $V_s \in \{V_1^-,V_2^-\}$, $s \in S$, we have
\[\dim(\cap_{s \in S} \rho(h_s)V_s) \leqslant \max(0,\dim(W_\rho) - \vert S \vert).\]
\end{enumerate}
Then $X_M$ is a nonempty Zariski-open set.
\end{lemma}
{\it Proof.}
Denote by $A_M$ (resp. $B_M$) the set of $M$-tuples satisfying condition (1) (resp. (2)). Then $A_M$ is clearly a finite intersection of Zariski open sets, which are non-empty by irreducibility of $\rho$.

Let us prove by induction over $M$ that $B_M$ is a non-empty Zariski open set. For $M = 1$, the condition (2) is empty, so this is clearly true. Assuming the result for $M$, let us check it for  $M+1$. 
Consider the finite collection of all vector spaces of the form $E_\alpha = \cap_{s \in S} \rho(h_s)V_s \subset W_\rho$, where $S \subset \{1,\dots,M\}$ and $V_s \in \{V_1^-,V_2^-\}$, for all $s \in S$. Then $B_{M+1}$ is equal to
 \[\Big\{(h_1,\dots,h_{M+1}) \, \Big\vert \, (h_1,\dots,h_{M}) \in B_M \text{ and } E_\alpha \nsubseteq \rho(h_{M+1})(V_1^- \cup V_2^-) \text{ for all } \alpha \text{ with } E_\alpha \neq\emptyset \Big\}.\]
Using this, it can be easily seen that $B_{M+1}$ is a finite intersection of non-empty Zariski open sets. Therefore, $B_{M+1}$ is Zariski open, as well as non-empty by  the Zariski connectedness of $\G$.
\hfill$\blacksquare$

{\it Proof of Proposition \ref{ping-pong}.}
Consider the representations $\rho_i$ over local fields $K_i$, $i \in I$, defined above. For $i\in I$, we consider the representation $\rho_{i'}: g \mapsto \rho_i(g^{-1})^t$. Note that by the definition of $r_i$, we clearly have $r_{i'} = r_i$.

Applying Lemma \ref{firstplayer} to the set
 of representations $\{\rho_i\}_{i\in I}\cup\{\rho_{i'}\}_{i\in I}$,  we obtain a sequence $(g_n)_n \subset \Gamma$ which converges to $1$ in the real topology, and elements $b_i, b_{i'}  \in \End(V_i)$ such that for all $i \in I$,
\begin{itemize}
\item $\lim_n \frac{\rho_i(g_n)}{\Vert \rho_i(g_n)\Vert_{K_i}} = b_i$ and  $\lim_n \frac{\rho_i(g_n^{-1})}{\Vert \rho_i(g_n^{-1})\Vert_{K_i}} = b_{i'}$,
\item $\rk(b_i) = \rk(b_{i'}) = r_i$, and
\item $\Range(b_i) \cap \Ker(b_i) =  \Range(b_{i'}) \cap \Ker(b_{i'}) =  \{0\}$.
\end{itemize}
By Lemma \ref{inverse}, we can add the following property to the above list:
\begin{itemize}
\item $\Range(b_i) \subset \Ker(b_{i'})$ and $\Range(b_{i'}) \subset \Ker(b_i)$.
\end{itemize}
Now, for $i \in I$, the sets $V_{i,1}^+ = \Range(b_i)$, $V_{i,1}^- = \Ker(b_i)$, $V_{i,2}^+ = \Range(b_{i'})$ and $V_{i,2}^- = \Ker(b_{i'})$ satisfy the hypothesis of Lemma \ref{conjugators}. Put $M := \max_i(\dim(\rho_i)) + 1$. For $i\in I$, denote by $X_i$ the non-empty Zariski open subset of $\G(K_i)^M$ given by Lemma \ref{conjugators} applied to these sets. Pick an $M$-tuple $(h_1,\cdots,h_M) \in \Gamma \bigcap (\cap_i X_i)$. 

Before going further, let us mention that $\rho_i(h_s)V_{i,1}^+ = \Range(h_sb_{i}h_s^{-1})$, $\rho_i(h_s)V_{i,1}^- = \Ker(h_sb_{i}h_s^{-1})$, whereas $\rho_i(h_s)V_{i,2}^+ = \Range(h_sb_{i'}h_s^{-1})$, $\rho_i(h_s)V_{i,2}^- = \Ker(h_sb_{i'}h_s^{-1})$.

Then by the definition of $X_i$, for every $i \in I$ and $1 \leqslant s \neq t \leqslant N$, we have $\rho_i(h_s)V_{i,1}^+ \nsubseteq \rho_i(h_t)V_{i,2}^{-}$. This means that $(h_tb_{i'}h_t^{-1}).(h_sb_{i}h_s^{-1}) \neq 0$. But both $(h_sb_{i}h_s^{-1})$ and $(h_tb_{i'}h_t^{-1})$ belong to $A_i$ and have rank $r_i$. Thus, their product has rank equal to $r_i$ by Lemma \ref{stability}. From this we deduce that $\rho_i(h_s)V_{i,1}^+ \cap \rho_i(h_t)V_{i,2}^{-} = \{0\}$.
Similarly, $\rho_i(h_s)V_{i,1}^+ \cap \rho_i(h_t)V_{i,1}^{-} = \{0\}$ and $\rho(h_s)V_{i,2}^+ \cap \rho(h_t)(V_{i,1}^{-} \cup V_{i,2}^-)= \{0\}$.

Using the above properties,  for every $i \in I$ and $1 \leqslant s \leqslant N$, we can find compact neighborhoods $K_{i,s}, K_{i,s}'  \subset \P(V_i)$ of $\P(\rho_i(h_s)V_{i,1}^+)$ and $\P(\rho_i(h_s)V_{i,2}^+)$ respectively, and open sets $U_{i,s}, U_{i,s}' \subset \P(V_i)$ which are complements of neighborhoods of $\P(\rho_i(h_s)V_{i,1}^-)$ and $\P(\rho_i(h_s)V_{i,2}^-)$, respectively, such that:
\begin{itemize}
\item $K_{i,s} \subset U_{i,s}$ and $K_{i,s}' \subset U_{i,s}'$ for all $s$;
\item $K_{i,s}' \cap U_{i,s} = \emptyset = K_{i,s} \cap U_{i,s}'$; 
\item For all $s \neq t$, $K_{i,s} \subset U_{i,t} \cap U_{i,t}'$ and $K_{i,s}' \subset U_{i,t} \cap U_{i,t}'$;
\item For any $x \in \P(V_i)$, we can find at least two indices $s$ for which $x \in U_{i,s}$ or $x \in U_{i,s}'$.
\end{itemize}
The last fact is due to property (2) from Lemma \ref{conjugators}, which implies that for any set $S \subset \{1,\cdots,M\}$ with $\vert S \vert = M - 1$ and any choice of $V_s \in \{V_{i,1}^-,V_{i,2}^-\}$, $s \in S$, we have $\cap_{s \in S} \rho(h_s)V_s = \{0\}$.

Finally, given $\eta > 0$, we can find $n$ large enough so that for all $i \in I$ and all $s$, we have $h_sg_nh_s^{-1}(U_{i,s}) \subset K_{i,s}$, $h_sg_n^{-1}h_s^{-1}(U_{i,s}') \subset K_{i,s}'$ and $h_sg_nh_s^{-1}, h_sg_n^{-1}h_s^{-1} \in B_\eta(1)$. We define $S$ to be the set of elements $\{h_sg_nh_s^{-1} \, \vert 1 \leqslant s \leqslant M\}$. If $g =  h_sg_nh_s^{-1} \in S$, define $K_g^{(i)} = K_{i,s}$ and $U_g^{(i)} = U_{i,s}$, and if $g =  h_sg_n^{-1}h_s^{-1} \in S^{-1}$, define $K_g^{(i)} = K_{i,s}'$ and $U_g^{(i)} = U_{i,s}'$.
These sets are easily seen to satisfy the desired properties.
\hfill$\blacksquare$

\subsection{From subgroups to neighborhoods of subgroups}
The goal of this section is to prove the following proposition, which roughly says that algebraic points with small logarithmic height cannot be very close to a proper algebraic subgroup. Our method is fairly similar to~\cite[Proposition 16]{SGV11} (see also \cite[Proposition 3.11]{BdS14} or \cite[Proposition 4]{Va10}).

\begin{proposition}\label{neigh}
Let $G$ be a connected simple Lie group and $T \subset G$ a finite subset. Assume that there there is a basis $\Bfr$ of the Lie algebra $\gfr$ of $G$ such that the matrix of Ad$(g)$ in the basis $\Bfr$ has algebraic entries, for every $g\in T$. 

Then there exists a constant $C>0$ (depending on $T$)  such that for every integer $n\geqslant 1$ and any non-discrete proper closed subgroup $H<G$, we can find a proper closed subgroup $H'<G$ such that
\[W_{\leqslant n}(T) \cap H^{(e^{-Cn})}\subseteq H',\]
where $W_{\leqslant n}(T) = \{g_1g_2...g_n \, | \, g_1,g_2,...,g_n \in T \cup T^{-1}\}$.
\end{proposition}

{\bf Notation}. In this subsection, we use the notation $O_X(a)$ to denote a positive quantity bounded by $Ca$, for some constant $C>0$ depending only on $X$. We also use the notation $a\gg_{X}b$ to mean the existence  of some constant $C>0$ depending only on $X$ such that $a\geqslant Cb$. 

\begin{lemma}\label{l:OneLargeCoordinate}
Let $X\subseteq {\rm M}_n(\bbr)$ be a finite subset. Suppose the $\bbr$-span $A$ of $X$ is an $\bbr$-algebra, and $V:=\bbr^n$ is a simple $A$-module. Then there exists $c_0 > 0$ such that for every $\lbf \in V^{\ast}$ and $\vbf \in V$ 
\[\max_{x\in X} |\lbf (x\vbf)|\geqslant c_0 \|\lbf\|_2\|\vbf\|_2.\]
\end{lemma}

{\it Proof}.
Let $H_{X}(\lbf,\vbf):=\max_{x\in X}|\lbf(x\vbf)|$. We need to show that 
the infimum of $H_X(\lbf,\vbf)$ on the pair of unit vectors is positive. Suppose the contrary. So by the continuity of $H_X:V^{\ast}\times V\rightarrow \bbr$, there are unit vectors $\lbf_0$ and $\vbf_0$ such that $H_X(\lbf_0,\vbf_0)=0$. This implies that for any $a\in A$ we have $\lbf_0(a \vbf_0)=0$. Hence the $A$-module generated by $\vbf_0$ is a proper subspace which contradicts the simplicity of $V$. 
\hfill$\blacksquare$

\begin{lemma}\label{l:InProperSubvariety}
Let $G$ be a simple Lie group and $T \subseteq G$ be a finite symmetric set such that $\Gamma = \langle T \rangle$ is a dense subgroup of $G$.  Suppose that the matrix of Ad$(g)$ with respect to a basis $\Bfr$ of the Lie algebra $\frak g$ of $G$ has algebraic entries, for every $g\in T$. Then there exists $C_1 > 0$ such that the following holds:

If $n\geqslant 1$ is an integer, then for any proper non-discrete closed subgroup $H$ of $G$, there are non-zero vectors $\vbf\in \gfr\otimes_{\bbr}\bbc$ and $\lbf\in\gfr^{\ast}\otimes_{\bbr}\bbc$ such that 
$\lbf(\Ad(\gamma)(\vbf))=0$, for any $\gamma \in W_{\leqslant n}(T) \cap H^{(e^{-C_1n})}$.
\end{lemma}

{\it Proof.}
Since $\Gamma$ is a dense subgroup of $G$, the $\bbr$-span $A$ of $\Ad(\Gamma)$ in ${\rm End}_{\bbr}(\gfr)$ is equal to the $\bbr$-span of $\Ad(G)$. Denote by $d$ the dimension of $G$. It is easy to see that the $\bbr$-span of $W_{\leqslant d^2}(T)$ is equal to $A$. Hence by Lemma~\ref{l:OneLargeCoordinate}, there exists $c_0 > 0$ such that for any $\lbf\in \gfr^{\ast}$ and $\vbf\in \gfr$ we have
\be\label{e:LargeCoordinate}
\max_{\gamma\in W_{\leqslant d^2}(T)} |\lbf (\Ad(\gamma)(\vbf))|\geqslant c_0\|\lbf\|_2\|\vbf\|_2,
\ee
as the adjoint representation is irreducible. 

Let $H<G$ be a proper non-discrete closed subgroup and fix $n \geq 1$. Let $\vbf\in \gfr, \lbf\in \gfr^{\ast}$ such that 
\begin{enumerate}
	\item $\|\vbf\|_2=1$ and $\|\lbf\|_2=1$.
	\item $\vbf\in \hfr$ and $\hfr\subseteq \ker \lbf$ where $\hfr:=\Lie(H)$ is the Lie algebra of $H$.
\end{enumerate} 
By using (\ref{e:LargeCoordinate}) and rescaling $\vbf$, we find $\gamma_0\in W_{\leqslant d^2}(T)$, $\vbf_H\in \gfr$, and $\lbf_H\in \gfr^{\ast}$ such that
\begin{enumerate}
	\item $\|\lbf_H\|_2=1$, $\|\vbf_H\|_2 \leqslant 1/c_0$.
	\item $\lbf_H(\Ad(\gamma_0)(\vbf_H))=1$.
	\item for any $h\in H$, $\lbf_H(\Ad(h)(\vbf_H))=0$.
\end{enumerate} 
By the hypothesis, $\gfr$ has a basis $\Bfr:=\{\vbf_1,\ldots, \vbf_d\}$ such that $\vbf_i^{\ast}(\Ad(\gamma)(\vbf_j))\in \overline{\bbq}$, for any $\gamma\in \Gamma$ and $1\leqslant i,j\leqslant d$, where $\Bfr^{\ast}:=\{\vbf_1^{\ast},\ldots,\vbf_d^{\ast}\}$ is the dual basis. Since $\Gamma$ is finitely generated, there are a number field $k$ and a finite set of places $\cS$ of $k$ such that $\vbf_i^{\ast}(\Ad(\gamma)(\vbf_j))\in \ocal_k(\cS)$, for any $\gamma\in \Gamma$. 

For any $g\in G$, let $\eta_g(\underline{x},\underline{y}) \in \bbr[x_1,\ldots,x_d,y_1,\ldots,y_d]$ be the polynomial $\eta_g([\lbf]_{\Bfr^{\ast}},[\vbf]_{\Bfr}):=\lbf(\Ad(g)(\vbf))$, where $[\lbf]_{\Bfr^{\ast}}$ (resp. $[\vbf]_{\Bfr}$) is the vector of coordinates of $\lbf$ in the basis $\Bfr^{\ast}$ (resp. $\Bfr$). It is clear that $\eta_g$ is a degree $2$ polynomial in $2d$ variables. Fix a constant $C_1 > 0$ large enough, depending only on $T$ (we will be more specific later).
Now, suppose that the following system of polynomial equations do not have a common solution over $\bbc$
\begin{align*}
 \eta_{\gamma}(\underline{x},\underline{y}) &=0 \h\h\text{ for any }\gamma\in  W_{\leq n}(T)\cap H^{(e^{-C_1n})},\\
\eta_{\gamma_0}(\underline{x},\underline{y})-1&=0.
\end{align*}
We notice that the coefficients of $\eta_{\gamma}$ are in $\ocal_k(\cS)$. We view $\ocal_k(\cS)$ as a discrete subring of $\prod_{\pfr\in V_k(\infty)\cup \cS}k_{\pfr}$, where $V_k(\infty)$ is the set of Archimedean places of $k$. It is clear that the $\cS$-norm (the maximum norm in $\prod_{\pfr\in V_k(\infty)\cup \cS}k_{\pfr}$) of the coefficients of $\eta_{\gamma}$ for $\gamma\in W_{\leqslant n}(T)$ is at most $e^{O_{T}(n)}$. Then by the effective Nullstellensatz~\cite[Theorem IV]{MW} there are polynomials $q_{\gamma}(\underline{x},\underline{y}),q_{\gamma_0}(\underline{x},\underline{y})\in \ocal_k[x_1,\ldots,x_d,y_1,\ldots,y_d]$ and $a\in \ocal_k$ such that
\begin{enumerate}
	\item $\sum_{\gamma\in W_{\leqslant n}(T)\cap H^{(e^{-C_1n})}} q_{\gamma}(\underline{x},\underline{y}) \eta_{\gamma}(\underline{x},\underline{y})+ q_{\gamma_0}(\underline{x},\underline{y}) \eta_{\gamma_0}(\underline{x},\underline{y})=a.$
	\item $\deg q_{\gamma}, \deg q_{\gamma_0}\ll_{d,\deg k} 1$. 
	\item The $\cS$-norms of the coefficients of $q_{\gamma}$ and $q_{\gamma_0}$ are at most $e^{O_{T}(n)}$. 
	\item The $\cS$-norm of $a$ is at most $e^{O_{T}(n)}$, and it is non-zero. 
\end{enumerate}

Since $a\in \ocal_k$ is non-zero, we have $1\leqslant |N_{k/\bbq}(a)|=\prod_{\pfr\in V_k(\infty)}|a|_{\pfr}\leqslant (\min_{\pfr\in V_k(\infty)}|a|_{\pfr})\|a\|_S^{\deg k-1}$. Thus 
\be\label{e:LowerBoundOna}
\min_{\pfr\in V_k(\infty)}|a|_{\pfr}\geqslant e^{-O_{T}(n)}.
\ee
Suppose $\pfr_0\in V_k(\infty)$ is the place which gives us the embedding of $\Ad(\Gamma)$ into ${\rm End}_{\bbr}(\gfr)$. 

So by the properties of $\lbf_H$ and $\vbf_H$ mentioned above we have that
\begin{align*}
|\eta_{\gamma}(\lbf_H,\vbf_H)|_{\pfr_0}&\leqslant e^{-C_1n/2}, \\
|q_{\gamma}(\lbf_H,\vbf_H)|_{\pfr_0}&\leqslant e^{O_{T}(n)},
\end{align*}
for any $\gamma\in \textstyle W_{\leqslant n}(T) \cap H^{(e^{-C_1n})}$. Hence we have
\[\textstyle
|\sum_{\gamma\in W_{\leqslant n}(T)\cap H^{(e^{-C_1n})}} q_{\gamma}(\lbf_H,\vbf_H) \eta_{\gamma}(\lbf_H,\vbf_H)+ q_{\gamma_0}(\lbf_H,\vbf_H) \eta_{\gamma_0}(\lbf_H,\vbf_H)|_{\pfr_0}\leqslant e^{O_{T}(n)-C_1n/2}\leqslant e^{-C_1n/4}
\]
if we chose $C_1$ large enough. But if we chose $C_1$ perhaps even larger (but still depending only on $T$) this contradicts \eqref{e:LowerBoundOna}.
\hfill$\blacksquare$

{\it Proof of Proposition~\ref{neigh}.}
Let $\bbg$ be the Zariski-closure of $\Ad(G)$ in $\GL(\gfr)$. By Lemma~\ref{l:InProperSubvariety}, there exists a constant $C_1 > 0$ such that for any $n$ and any non-discrete proper closed subgroup $H$ of $G$ there is a variety $X$ (depending on $H$ and $n$) of $\bbg$ whose dimension is strictly less than $\dim \bbg$ such that $\Ad(W_{\leqslant n}(T) \cap H^{(e^{-C_1n})})\subseteq X$. Using the generalized Bezout theorem it was proved in \cite[Proposition 3.2]{EMO05} that there is $N(X)\geqslant 1$ such that
$W_{\leqslant N(X)}(A) \not\subseteq X$ whenever $A$ generates a Zariski-dense subgroup of $\bbg$. 
Moreover, by the proof of \cite[Proposition 3.2]{EMO05}, $N(X)$ is bounded above by some bound depending on the number of irreducible components of $X$ and the maximal degree of an irreducible component of $X$. Since $X$ is the intersection of $\bbg$ with a hyperplane, we conclude that $N:=\sup_{X}N(X)<\infty$. This number $N$ only depends on $T$.

Next, we show that there exists $C > 0$ (depending only on $T$) such that for all multiple $n$ of $N$,
\be\label{e:ProperSubgroup}
\textstyle W_{\leqslant N}\Big(W_{\leqslant n/N}(T)) \cap H^{(e^{-Cn})}\Big) \subseteq W_{\leqslant n}(T) \cap H^{(e^{-C_1n})}.
\ee
This - coupled with the above paragraph - implies that for all multiple $n$ of $N$ and all proper closed subgroup $H$ of $G$, the set $W_{\leqslant n/N}(T) \cap H^{(e^{-Cn})}$ is contained in a proper algebraic subgroup of $\bbg$.  

For any $\gamma_i\in W_{\leqslant n/N}(T) \cap H^{(e^{-Cn})}$, there are $h_i\in H$ such that $\|\Ad(\gamma_i)-\Ad(h_i)\|_2\leqslant e^{-Cn}$ and $\|\Ad(\gamma_i)\|_2\leqslant e^{O_{T}(n)}$. Hence, $\|\Ad(h_i)\|_2\leqslant e^{O_{T}(n)}$ and 
\begin{align*}
\|\Ad(\gamma_1\cdots \gamma_N)-\Ad(h_1\cdots h_N)\|_2&=\|\sum_{i=0}^{N-1}(\Ad(\gamma_1\cdots \gamma_{N-i}h_{N-i+1} \cdots h_N)-\Ad(\gamma_1\cdots \gamma_{N-i-1}h_{N-i} \cdots h_N)) \|_2\\
&\leqslant \sum_{i=0}^{N-1} (\prod_{j=1}^{N-i-1} \|\Ad(\gamma_j)\|_2) (\prod_{j=N-i+1}^N\|\Ad h_j\|_2) (\|\Ad \gamma_{N-i}-\Ad h_{N-i}\|_2) \\
& \leqslant e^{O_{T}(n)-Cn} \leqslant e^{-C_1n}
\end{align*}
if $C'\gg_{T} 1$, which implies (\ref{e:ProperSubgroup}).
\hfill$\blacksquare$

\subsection{Proof of Theorem \ref{escape}} 
By \cite[Corollary 2.5]{BrG02}, $\Gamma$ contains a finitely generated subgroup which is dense in $G$. Thus, we may assume that $\Gamma$ is finitely generated. 
Let $\rho_i:\Gamma\rightarrow\text{GL}(V_i)$, $i\in I$, be the representations and $M\geqslant 2$ be the integer given by Theorem \ref{ping-pong}.

By a result of Kazhdan and Margulis (see \cite[Theorem 8.16]{Ra72}), there is a neighborhood  $U$ of the identity in $G$  such that for any discrete subgroup $\Sigma<G$,  $\Sigma\cap U$ is contained in a connected nilpotent subgroup of $G$.
Let $U_0\subset U$ be an open set such that $U$ contains the closure of $U_0^{-1}U_0$.

Throughout the proof, we fix two constants $\kappa>1$ and $\eta>0$ (depending on $G$ only) such that
\begin{enumerate}[(a)]
\item\label{a}  $B_R(1)$ can be covered by at most $R^{\kappa}$  of the sets $\{gU_0\}_{g\in G}$, whenever $R>2$,
\item\label{b} $B_R(1)$ can be covered by at most ${\Big(\frac{R}{r}\Big)^{\kappa}}$ balls in $G$ of radius ${\frac{r}{2}}$, whenever $R>2r>0$,
\item\label{c} $\|x^{-1}\|_2\leqslant \|x\|_2^{\kappa}$, for every $x\in G$, and
\item\label{d} ${(1+\eta)^{(3\kappa +4)\kappa}<\Big(\frac{2M-1}{2M-2}\Big)^{\frac{1}{13}}}$.
\end{enumerate}

Let $S\subset\Gamma$ be a set satisfying Theorem \ref{ping-pong} such that $\tilde S=S\cup S^{-1}\subset B_{\eta}(1)$ and $|S|=M$. 
For $i\in I$, let  $K_g^{(i)}\subset U_g^{(i)}$ ($g\in\tilde S$) be the subsets of $V_i$ provided by Theorem \ref{ping-pong}.
The usual ping-pong lemma implies that $S$ freely generates a subgroup of $\Gamma$, which we denote by $\langle S\rangle$. 
Let $|g|_{S}$ be the length of an element  $g\in\langle S\rangle$ with respect to $\tilde S$. We denote by $W_n(S)$ the set of elements of length $n$, and by $W_{\leqslant n}(S)$ the set of elements of length at most $n$.


Let $\ell\geqslant 1$ be an integer and put ${\varepsilon=(1+\eta)^{-\ell}}$. 
In part 1 of the proof, we construct a  finite set $T\subset \Gamma\cap B_{\varepsilon}(1)$. Our construction is inspired by the proof of \cite[Lemma 3]{BY11}.
In the rest of the proof (parts 2-4), we provide constants $d_1,d_2 > 0$ and show that $T$ satisfies the conclusion of Theorem \ref{escape}, whenever $\ell$ is large enough.
This will clearly imply Theorem \ref{escape}.

\vskip 0.1in
{\bf Part 1: construction of the set $T$.}\label{part1} 
\vskip 0.05in

Let $a,b\in S$ with $a \not= b$ and define $$Y=\{w=s_1s_2...s_{\ell}|s_1=a,s_{\ell}=b,s_2,...,s_{\ell-1}\in\tilde S,s_{i+1}\not=s_i^{-1},\text{for all}\;1\leqslant i<\ell\}.$$

 Let $Z=\{w^3|w\in Y\}.$
Since $|\tilde S|=2M$, we get that  $|Z|=|Y|\geqslant (2M-1)^{\ell-3}$. 
Since $\tilde S\subset B_{\eta}(1)$,  it follows that $W_n(S)\subset B_{(1+\eta)^{n}}(1)$, for all $n\geqslant 1$.
Since $Z\subset W_{3\ell}(S)$, we get that $Z\subset B_{(1+\eta)^{3\ell}}(1)$. 
By using (\ref{b}), $Z$ can be covered by at most ${\Big[\frac{(1+\eta)^{3\ell}}{\frac{\varepsilon}{(1+\eta)^{3\kappa\ell}}}\Big]^{\kappa}=(1+\eta)^{(3\kappa+4)\kappa\ell}}$ balls of radius ${\frac{\varepsilon}{2(1+\eta)^{3\kappa\ell}}}$.

From this we deduce that there is $g_0\in Z$ such that \begin{equation} |B_{\frac{\varepsilon}{(1+\eta)^{3\kappa\ell}}}(g_0)\cap Z|\geqslant\frac{|Z|}{(1+\eta)^{(3\kappa+4)\kappa\ell}}\geqslant\frac{(2M-1)^{\ell-3}}{(1+\eta)^{(3\kappa+4)\kappa\ell}} \end{equation}

We define $T=g_0^{-1}(B_{\frac{\varepsilon}{(1+\eta)^{3\kappa\ell}}}(g_0)\cap Z)\setminus \{1\}$ and $\tilde T=T\cup T^{-1}$. Then $|T|\geqslant \frac{(2M-1)^{\ell-3}}{(1+\eta)^{(3\kappa+4)\kappa\ell}}-1$.  Since by inequality (\ref{d}) we have that  ${\frac{2M-1}{(1+\eta)^{(3\kappa+4)\kappa}}>(2M-2)^{\frac{1}{13}}(2M-1)^{\frac{12}{13}}}$, we get that \begin{equation}\label{T}|T|\geqslant \frac{[(2M-2)^{\frac{1}{13}}(2M-1)^{\frac{12}{13}}]^{\ell}}{(2M-1)^4}\geqslant [(2M-2)^{\frac{1}{13}}(2M-1)^{\frac{12}{13}}]^{\ell-5},\;\;\text{for all}\;\;\ell\geqslant 1.\end{equation}

If $g\in T$, then ${\|g_0g-g_0\|_2\leqslant\frac{\varepsilon}{(1+\eta)^{3\kappa\ell}}}$. Since $\|g_0\|_2\leqslant (1+\eta)^{3\ell}$, we get $\|g_0^{-1}\|_2\leqslant \|g_0\|_2^{\kappa}\leqslant (1+\eta)^{3\kappa \ell}$.
Altogether, it follows that $\|g-1\|_2\leqslant \|g_0^{-1}\|_2\|g_0g-g_0\|_2\leqslant\varepsilon$, for all $g\in T$. Hence $T\subset B_{\varepsilon}(1)$. 

We end part (1) of the proof by recording a useful property of $T$.

\vskip 0.05in
{\bf Claim 1.}
If  $g\in W_n(T)$, then $n\ell\leqslant |g|_S\leqslant 6n\ell$. Thus, $T$ freely generates a free subgroup of $\Gamma$.

{\it Proof.} It is  enough to show that  $n\ell\leqslant |g|_S\leqslant 3n\ell$, for all $g\in W_n(Z)$ and $n\geqslant 1$.
Let $g=g_n^{\varepsilon_n}...g_1^{\varepsilon_1}$, where $n\geqslant 1$ and $g_1,...,g_n\in Z$, $\varepsilon_1,...,\varepsilon_n\in\{\pm1\}$ are such that $g_{i+1}^{\varepsilon_{i+1}}g_{i}^{\varepsilon_{i}}\not=1$, for all $1\leqslant i\leqslant n-1$. 
Let $w_1,...,w_n\in  Y$ such that $g_1=w_1^3,...,g_n=w_n^3$. Then $$g=w_n^{2\varepsilon_n}(w_n^{\varepsilon_n}w_{n-1}^{\varepsilon_{n-1}})w_{n-1}^{\varepsilon_{n-1}}(w_{n-1}^{\varepsilon_{n-1}}w_{n-2}^{\varepsilon_{n-2}})...(w_{2}^{\varepsilon_{2}}w_1^{\varepsilon_1})w_1^{2\varepsilon_1}.$$

Since $w_{i+1}^{\varepsilon_{i+1}}w_{i}^{\varepsilon_{i}}\not=1$ and $w_i^{\varepsilon}w_i^{\varepsilon}$ is already reduced,  after making all the possible cancellations, the middle $w_i^{\varepsilon_i}$ from
$g_i^{\varepsilon_i}=w_i^{\varepsilon_i}w_i^{\varepsilon_i}w_i^{\varepsilon_i}$ will not be affected. This implies the conclusion. \hfill$\square$

\vskip 0.1in
{\bf Part 2: bounding the number of returns.} 
\vskip 0.05in

We continue by showing that the number of elements $g\in W_n(T)$ which fix a given line $[v]\in\P(V_i)$, for some $i\in I$, is bounded above by $|W_n(T)|^{1-c_0}$, for a constant $c_0>0$.

\vskip 0.05in
{\bf Claim 2.}
There exist $n_0 \geqslant 1$,  $\ell_0\geqslant 1$ and $c_0>0$ such that given $i\in I$ and $v\in V_i$, we have $$|\{g\in W_n(T)|\rho_i(g)([v])=[v]\} \vert \leqslant |W_n(T)|^{1-c_0},\;\;\;\text{for all $\ell\geqslant\ell_0$ and $n\geqslant n_0$.}$$

{\it Proof of Claim 2.} 
Let $i\in I$ and $v\in V_i$.
For simplicity, we denote $\rho=\rho_i$ and $K_g=K_g^{(i)}$, $U_g=U_g^{(i)}$, for $g\in\tilde S$. Fix  $n\geqslant 1$ and define $A=\{g\in W_n(T)|\rho(g)([v])=[v]\}$. Denote $N=\lfloor{\frac{n\ell}{2}}\rfloor$.

In order to estimate $|A|$, we partition $A$ into two subsets according to the reduced form of $g$.   Let $g\in A$ and $g=k_pk_{p-1}...k_1$ be its reduced form with respect to $\tilde S$, where $p=|g|_{S}$ and $k_1,...,k_p\in\tilde S$.
By Claim 1 we get that $2N\leqslant n\ell\leqslant p\leqslant 6n\ell\leqslant 12N+6$.  We define $$B=\{g\in A|\rho(k_Nk_{N-1}...k_1)([v])\in U_{k_{N+1}}\}\;\;\;\text{and}\;\;\;C=A\setminus B.$$

We proceed by estimating $|B|$ and $|C|$ separately.

\vskip 0.05in
{\bf Claim 3.} ${|B|\leqslant (2|T|)^{\frac{11 n}{12}+1}}$, for all $n >12$.

{\it Proof of Claim 3.}
Assume that  $g=k_pk_{p-1}...k_1\in B$.  Then the first part of Lemma \ref{SGV} implies that $k_{N+1},...,k_{p-1},k_{p}$ are uniquely determined by $v$. 

Now, since $g\in W_n(T)$, we can  write $g=g_n^{\varepsilon_n}g_{n-1}^{\varepsilon_{n-1}}...g_1^{\varepsilon_1}$, where $g_1,...,g_n\in T$, $\varepsilon_1,...,\varepsilon_n\in\{\pm 1\}$ and $g_j^{\varepsilon_j}\not=g_{j+1}^{\varepsilon_{j+1}}$, for all $1\leqslant j<n$. Let $w_0\in Y$ such that $g_0=w_0^3$ and $w_1,...,w_n\in Y\setminus\{w_0\}$ such that $g_1=w_0^{-3}w_1^{3},...,g_n=w_0^{-3}w_n^{3}$.
Then the reduced form can be written as $g=h_nw_n^{\varepsilon_n}h_{n-1}...h_1w_1^{\varepsilon_1}h_0$, where $h_i\in\Gamma$ satisfies $|h_i|_{S}\leqslant 5\ell$, and the factor $w_i^{\varepsilon_i}$ corresponds to the middle $w_i^{\varepsilon_i}$ from $w_i^{3\varepsilon_i}$.

 We claim that $g_q^{\varepsilon_q}$ is uniquely determined, for any $q$ such that $n\geqslant q\geqslant\frac{11n}{12}+1$.
 More precisely, we will show by induction that
 $\varepsilon_q$, $w_q$ and $h_q$  are uniquely determined, for all  $q$ with  $n\geqslant q\geqslant{\frac{11n}{12}}+1$. 

First, if $q=n$, we have that  either $h_n=w_0^{-3}w_n$, if $\varepsilon_n=1$, or $h_n=w_n^{-1}$, if $\varepsilon_n=-1$.
Since $|h_n|_{S}\leqslant 4\ell$ and $p-N\geqslant n\ell-N\geqslant \frac{n\ell}{2}>4\ell$, it follows that $w_n$ and $\varepsilon_n$ are determined. 
Specifically, there are two cases: (1) $k_p...k_{p-\ell+1}=w_0^{-1}$ or (2) $k_p...k_{p-\ell+1}=w_n^{-1}$. In case (1)  $\varepsilon_n=1$, $w_n=k_{p-3\ell}...k_{p-4\ell+1}$, and $h_n=w_0^{-3}w_n$, while in case  (2)  $\varepsilon_n=-1$, $w_n=k_{p-\ell+1}^{-1}...k_p^{-1}$ and $h_n=w_n^{-1}$.

Assume that  $\varepsilon_n,...,\varepsilon_{q+1}$, $w_n,...,w_{q+1}$, $h_n,...,h_{q+1}$ are determined, for some $q$ with  $n\geqslant q\geqslant {\frac{11 n}{12}}+1$.
Since $|h_nw_n^{\varepsilon_n}...h_{q+1}w_{q+1}^{\varepsilon_{q+1}}|_{S}\leqslant 6(n-q)\ell$ and $p-N\geqslant \frac{n\ell}{2}\geqslant 6(n-q+1)\ell=6(n-q)\ell+6\ell$, we deduce that the first $6\ell$ letters from the left in the reduced word $h_qw_q^{\varepsilon_q}...h_1w_1^{\varepsilon_1}h_0$ with respect to $S$ are determined. Note that $h_q\in\{w_{q+1}w_0^{-3}w_q, w_{q+1}w_q^{-1}\}$, if $\varepsilon_{q+1}=1$, and $h_q\in\{w_{q+1}^{-1}w_q, w_{q+1}^{-1}w_0^{3}w_q^{-1}\}$, if  $\varepsilon_{q+1}=-1$. Since $\varepsilon_{q+1}$, $w_{q+1}$ are determined and $|h_qw_q^{\varepsilon_q}|_{S}\leqslant 6\ell$, it follows easily that $\varepsilon_q,w_q$ and $h_q$ are determined. This finishes the proof of our assertion.

Therefore, if $q=\lfloor{{\frac{11n}{12}}}\rfloor+2$, then $g_n^{\varepsilon_n},...,g_q^{\varepsilon_q}\in\tilde T$ are  uniquely determined  for every $g=g_n^{\varepsilon_n}...g_1^{\varepsilon_1}\in B$. Since $g_1^{\varepsilon_1},...,g_{q-1}^{\varepsilon_{q-1}}$ can each take at most $2|T|$ values, we get that $|B|\leqslant (2|T|)^{q-1}\leqslant (2|T|)^{\frac{11 n}{12}+1}$.  \hfill$\square$

\vskip 0.05in
{\bf Claim 4.} $|C|\leqslant [4(2M-2)^{\ell}]^{{\frac{n}{6}}}(2|T|)^{n+2-{\frac{n}{6}}}$, for all $n\geqslant 1$.

{\it Proof of Claim 4.} Assume that $g=k_pk_{p-1}...k_1\in C$. Then the second part of Lemma \ref{SGV} implies that $\rho(k_jk_{j-1}...k_1)([v])\notin U_{k_{j+1}}$, for all $1\leqslant j\leqslant N$.
Below we will use this fact as follows. Suppose that $k_1,...,k_j$ are already determined, for some $1\leqslant j\leqslant N$. Since $\rho(k_jk_{j-1}...k_1)([v])$ belongs to at least $2$ of the sets $\{U_g\}_{g\in\tilde S}$ and $|\tilde S|=2M$,  we derive that $k_{j+1}\in\tilde S$ can take at most $2M-2$ values. 

Now, since $g\in W_n(T)$, we can  write $g=g_n^{\varepsilon_n}g_{n-1}^{\varepsilon_{n-1}}...g_1^{\varepsilon_1}$, where $g_1,...,g_n\in T$, $\varepsilon_1,...,\varepsilon_n\in\{\pm 1\}$ and $g_j^{\varepsilon_j}\not=g_{j+1}^{\varepsilon_{j+1}}$, for all $1\leqslant j<n$. Let $w_1,...,w_n\in Y\setminus\{w_0\}$ such that $g_1=w_0^{-3}w_1^{3},...,g_n=w_0^{-3}w_n^{3}$.

Let $q$ with $1\leqslant q\leqslant{\frac{n}{12}}-1$ and assume that $g_1^{\varepsilon_1},...,g_q^{\varepsilon_q}$ are already determined. In other words, assume that $w_1,...,w_q$ and $\varepsilon_1,...,\varepsilon_q$ are determined. Our goal is to estimate the number of possible values of $g_{q+1}^{\varepsilon_{q+1}}\in\tilde T$. 
Depending on the values of $\varepsilon_q,\varepsilon_{q+1}\in\{\pm1\}$ we are in one of four cases. 
We assume that $\varepsilon_q=\varepsilon_{q+1}=1$, since the estimates in the other three cases are entirely similar. In this case, we have
$g=g_n^{\varepsilon_n}...g_{q+2}^{\varepsilon_{q+2}}(w_0^{-3}w_{q+1}^2)(w_{q+1}w_0^{-1})(w_0^{-2}w_q^3)g_{q-1}^{\varepsilon_{q-1}}...g_1^{\varepsilon_1}$. Let $j$ such that we have  $(w_0^{-2}w_q^3)g_{q-1}^{\varepsilon_{q-1}}...g_1^{\varepsilon_1}=k_jk_{j-1}...k_1$. Then  $j$ and $k_1,...,k_j$ are determined. Note that $j\leqslant 6q\ell$.

Write $w_0=r_1...r_{\ell}$, $w_{q+1}=s_1...s_{\ell}$, where $r_1,...,r_{\ell},s_1,...,r_{\ell}\in\tilde S$. Notice that $|w_{q+1}w_0^{-1}|_{S}$ is even and $2\leqslant |w_{q+1}w_0^{-1}|_{S}\leqslant 2\ell-2$. 
Let $1\leqslant \ell'\leqslant \ell-1$ such that $|w_{q+1}w_0^{-1}|_{S}=2\ell'$. 

Assume that $1\leqslant\ell'\leqslant\ell-1$ is determined.
Then  $s_{\ell'+1}=r_{\ell'+1},...,s_{\ell}=r_{\ell}$, hence $s_{\ell'+1},...,s_{\ell}$ are determined. Since $w_{q+1}w_0^{-1}=s_1...s_{\ell'}r_{\ell'}^{-1}...r_1^{-1}$, we get that $k_{j+1}=r_1^{-1},...,k_{j+\ell'}=r_{\ell'}^{-1}$ and $k_{j+\ell'+1}=s_{\ell'},...,k_{j+2\ell'}=s_1$. Hence $k_1,...,k_{j+\ell'}$ are determined.
As $j+2\ell'\leqslant 6\ell q+2(\ell-1)<N$ we get $\rho(k_{j+\ell'}...k_1)([v])\notin U_{k_{j+\ell'+1}}$. The beginning of the proof implies that $k_{j+\ell'+1}$ and hence $s_{\ell'}$ can take at most $2M-2$ values. 
Moreover, if $k_{j+\ell'+1},...,k_{j+\ell'+p}$ are determined, for some $1\leqslant p\leqslant \ell'-1$, then since $\rho(k_{j+\ell'+p}...k_1)([v])\notin U_{k_{j+\ell'+p+1}}$, we deduce that $k_{j+\ell'+p+1}$ and therefore $s_{\ell'-p}$ can take at most $2M-2$ values. It follows that there are at most $(2M-2)^{\ell'}$ possibilities for $s_1,...,s_{\ell'}$.

We derive that in the case $\varepsilon_q=\varepsilon_{q+1}=1$,  the total number of possible values of $w_{q+1}$ is at most ${\sum_{\ell'=1}^{\ell-1}(2M-2)^{\ell'}\leqslant (2M-2)^{\ell}}$. By adapting the above argument, it follows
  that the number of possible values of $w_{q+1}$ is at most $(2M-2)^{\ell}$ in the other three cases as well. Altogether, we get that if $g_1^{\varepsilon_1},...,g_q^{\varepsilon_q}$ are determined and $1\leqslant q\leqslant{\frac{n}{12}}-1$, then $g_{q+1}^{\varepsilon_{q+1}}$ can take at most $4(2M-2)^{\ell}$ values. 
  
  Let $q=\lfloor{\frac{n}{12}}\rfloor$. Thus, if $g_1^{\varepsilon_1}$ is determined, then $g_2^{\varepsilon_2}...g_q^{\varepsilon_q}$ can take at most $[4(2M-2)^{\ell}]^{q-1}$ values. 
 As $g_1^{\varepsilon_1},g_{q+1}^{\varepsilon_{q+1}},...,g_{n}^{\varepsilon_{n}}$ can each take at most $2|T|$ values, we get that $|C|\leqslant [4(2M-2)^{\ell}]^{{\frac{n}{12}}}(2|T|)^{n+2-{\frac{n}{12}}}.$ This finishes the proof of Claim 4.
\hfill$\square$

\vskip 0.05in
{\it End of proof of Claim 2}.
By combining Claims 3 and 4, and using that $|T|\leqslant (2M-1)^{\ell-1}$, we get \begin{equation}\label{A}|A|\leqslant (2|T|)^{\frac{11 n}{12}+1}+
 [4(2M-2)^{\ell}]^{{\frac{n}{12}}}(2|T|)^{n+2-{\frac{n}{12}}}\leqslant [(2M-2)^{\frac{1}{12}}(2M-1)^{\frac{11}{12}}]^{(n+24)(\ell+2)}\end{equation}
 for all $n > 12$ and every $\ell\geqslant 1$. Equations \eqref{T} and \eqref{A} together imply that there exist $n_0 \geqslant 1$, $c_0>0$ and $l_0\geqslant 1$
such that $|A|\leqslant |T|^{(1-c_0)n}$, for all $n\geqslant n_0$ and every $\ell\geqslant\ell_0$. Since $|W_n(T)|=2|T|(2|T|-1)^{n-1}>|T|^{n}$, the conclusion of Claim 2 follows.
\hfill$\square$

\vskip 0.1in
{\bf Part 3: bounding the probability of return.}\label{part3}
\vskip 0.05in

 Define ${\mu=\frac{1}{2|T|}\sum_{g\in T}(\delta_g+\delta_{g^{-1}})}$.
By using Part 2 and following closely the proof of \cite[Proposition 9]{Va10} (see also \cite[Proposition 7]{SGV11}) we next estimate $\mu^{*n}(\{g\in\Gamma|\rho_i(g)([v])=[v]\})$.

\vskip 0.05in
{\bf Claim 5.} There exist $n_1 \geqslant 1$ and $c > 0$ such that for every $i\in I$, $v\in V_i$ and $\ell\geqslant\ell_0$, we have that
\[\mu^{*n}(\{g\in\Gamma|\rho_i(g)([v])=[v]\})\leqslant |T|^{-cn}, \;\;\;\text{for all $n\geqslant n_1$}.\]

{\it Proof of Claim 5.} Denote $\rho=\rho_i$ and $A=\{g\in\Gamma|\rho(g)([v])=[v]\}$. Let $n\geqslant 10n_0$, where $n_0$ is as in Part 2. For every $k\geqslant 1$, fix $g_k\in W_k(T)$. Then $\mu^{*n}(\{g\})=\mu^{*n}(\{g_k\})$, for all $g\in W_k(T)$. Since $\mu^{*n}$ is supported on words of length at most $n$ in $T$, we get
\[\mu^{*n}(A)=\sum_{k=0}^n\mu^{*n}(A\cap W_k(T))=\sum_{k=0}^n|A\cap W_k(T)|\mu^{*n}(\{g_k\}).\]

Let us now majorize each of the terms involved. First, by Kesten's theorem \cite{Ke59} we have that 
\[\mu^{*n}(\{g\})\leqslant {\Big(\frac{\sqrt{2|T|-1}}{|T|}\Big)}^{n},\;\;\text{for all $g\in\Gamma$}.\]
Moreover, we deduce from Part 2 that for $n \geqslant 10n_0$, we have that 
\[|A\cap W_k(T)| \leqslant |W_k(T)|^{1-c_0} \leqslant (2|T|-1)^{-\frac{c_0n}{10}}|W_k(T)| \text{ for all } k \geqslant n/10.\] 
When $k < n/10$, we use the brutal bound $|A\cap W_k(T)| \leqslant |W_k(T)| \leqslant (2|T|)^k$. Altogether, we get
\begin{align*}
\mu^{*n}(A) &\leqslant\sum_{1\leqslant k<\frac{n}{10}}(2|T|)^k{\Big(\frac{\sqrt{2|T|-1}}{|T|}\Big)}^{n}+(2|T|-1)^{-\frac{c_0n}{10}}\sum_{\frac{n}{10}\leqslant k\leqslant n}|W_k(T)|\mu^{*n}(\{g_k\})\\
&\leqslant\sum_{1\leqslant k<\frac{n}{10}} (2|T|)^k{\Big(\frac{\sqrt{2|T|-1}}{|T|}\Big)}^{n}+(2|T|-1)^{-\frac{c_0n}{10}}\\ 
&\leqslant (2|T|)^{\frac{n}{10}}{\Big(\frac{\sqrt{2|T|-1}}{|T|}\Big)}^{n}+(2|T|-1)^{-\frac{c_0n}{10}}.
\end{align*}

The conclusion of Claim 5 is now immediate.
\hfill$\square$

\vskip 0.1in
{\bf Part 4: end of the proof.} 
\vskip 0.05in

We are now ready to conclude the proof.
Let $n_1,c,\ell_0$ be given as above and let $C$ be the constant given by Proposition \ref{neigh}. Since $M\geqslant 2$, we get that $(2M-1)^{\frac{12}{13}}\geqslant 3^{\frac{12}{13}}>e$. 
By using \eqref{T}, and after taking a larger $\ell_0$, we may assume that $|T|\geqslant e^{\ell}$ and that $(4n+1)2^{(\kappa+1)n}\leqslant e^{\frac{n\ell}{4}}$, for any $n\geqslant 1$ and $\ell\geqslant\ell_0$.

\vskip 0.05in
{\bf Claim 6.} Let $\delta>0$ be small enough and $n$ be an integer such that $\frac{\log(1+\eta)}{7C}\frac{\log{\frac{1}{\delta}}}{\log{\frac{1}{\varepsilon}}}\leqslant n\leqslant\frac{\log(1+\eta)}{6C}\frac{\log{\frac{1}{\delta}}}{\log{\frac{1}{\varepsilon}}}$.

Then $\mu^{*n}(H^{(\delta)})\leqslant\delta^{\frac{\min\{c,\frac{1}{4}\}}{7C}}$, for every proper closed connected subgroup $H<G$.
\vskip 0.05in

{\it Proof of Claim 6}.
Fix $n$ as in the claim and let $H<G$ be a proper closed connected subgroup. 
Thanks to Proposition \ref{neigh}, we can find a proper closed subgroup $H'<G$ such that \[W_{\leqslant 6n\ell}(S)\cap H^{(e^{-C6n\ell})}\subset H'.\]

Let $g\in$ supp$(\mu^{*n})\cap H^{(\delta)}$. Then $g\in W_{\leqslant n}(T)$ and since $T\subset W_{6\ell}(S)$, we deduce that $g\in W_{\leqslant 6n\ell}(S)$. Since $\varepsilon=(1+\eta)^{-\ell}$, hence $\ell=\frac{\log{\frac{1}{\varepsilon}}}{\log(1+\eta)}$, the assumption on $n$ implies that $\delta\leqslant e^{-6n\ell C}$. 
By using the previous  paragraph, we derive that $g\in H'$. 
Since $\mu$ is supported on $T$, we also have $g\in\langle T\rangle$.
Denoting $\Gamma_0=\langle T\rangle\cap H'$, we therefore get that

\begin{equation}\label{muneigh}\mu^{*n}(H^{(\delta)})\leqslant\mu^{*n}(\Gamma_0).  \end{equation}

We continue by treating two separate cases:

{\bf Case 1}. $\Gamma_0$ is non-discrete in $G$.

In this case, since $\Gamma_0\subset\Gamma\cap H'$, we get that $\Gamma\cap H'$ is non-discrete. Proposition \ref{ping-pong} implies the existence of $i\in I$ and $[v]\in\mathbb P(V_i)$ such that $\rho_i(g)([v])=[v]$, for all $g\in\Gamma\cap H'$. 
Since $|T|\geqslant e^{\ell}$, by combining \ref{muneigh} with Part 3, we get that for $\delta>0$ small enough so that $n\geqslant n_1$,
\begin{equation}\label{est1}\mu^{*n}(H^{(\delta)})\leqslant\mu^{*n}(\Gamma\cap H')\leqslant \mu^{*n}(\{g\in\Gamma|\rho_i(g)([v])=[v]\})\leqslant |T|^{-cn}\leqslant e^{-cn\ell}.\end{equation}

Since $n \geqslant \frac{\log{(1+\eta)}\log{\frac{1}{\delta}}}{7C\log{\frac{1}{\varepsilon}}}$, we get $n\ell\geqslant\frac{\log{\frac{1}{\delta}}}{7C}$. This implies that $e^{-cn\ell}\leqslant\delta^{\frac{c}{7C}}$, proving the claim. 

{\bf Case 2.} $\Gamma_0$ is discrete in $G$.

In this case, by the definition of $U$, we have that $\Gamma_1:=\langle\Gamma_0\cap U\rangle$ is a nilpotent group. 
Since $\Gamma_1<\langle T\rangle$ and $\langle T\rangle$ is a free group, $\Gamma_1$ must be a cyclic group. As a consequence, we have that $$|\Gamma_0\cap U\cap\text{supp}(\mu^{*2n})|\leqslant |\Gamma_1\cap\text{supp}(\mu^{*2n})|=|\Gamma_1\cap W_{\leqslant 2n}(T)|\leqslant 4n+1.$$

Next, if $N=\lfloor (1+\varepsilon)^{\kappa n}\rfloor$, then by (\ref{a}) we can find $g_1,...,g_N\in G$ such that $B_{(1+\varepsilon)^n}(1)\subset\cup_{i=1}^Ng_iU_0$. Since supp$(\mu^{*n})\subset B_{(1+\varepsilon)^n}(1)$, we thus get that $\Gamma_0\cap\text{supp}(\mu^{*n})\subset\cup_{i=1}^N(\Gamma_0\cap g_iU_0\cap\text{supp}(\mu^{*n})).$ Recall that $U_0^{-1}U_0\subset U$. So, if $1\leqslant i\leqslant N$ and $x,y\in \Gamma_0\cap g_iU_0\cap\text{supp}(\mu^{*n})$, then $x^{-1}y\in \Gamma_0\cap U\cap\text{supp}(\mu^{*2n})$. This implies that $|\Gamma_0\cap g_iU_0\cap\text{supp}(\mu^{*n})|\leqslant |\Gamma_0\cap U\cap\text{supp}(\mu^{*2n})|\leqslant 4n+1$, for every $1\leqslant i\leqslant N$. 

Altogether, we get that $|\Gamma_0\cap\text{supp}(\mu^{*n})|\leqslant (4n+1)N\leqslant (4n+1)(1+\varepsilon)^{\kappa n}.$
In combination with Kesten's theorem, we derive that $$\mu^{*n}(\Gamma_0)\leqslant (4n+1)(1+\varepsilon)^{\kappa n}{\Big(\frac{\sqrt{2|T|-1}}{|T|}\Big)}^{n}.$$

Since $|T|\geqslant e^{\ell}$, we get $\frac{\sqrt{2|T|-1}}{|T|}\leqslant\frac{2}{\sqrt{|T|}}\leqslant 2e^{-\frac{\ell}{2}}$. Since $(4n+1)(1+\varepsilon)^{\kappa n}2^n\leqslant (4n+1)2^{(\kappa+1)n}\leqslant e^{\frac{n\ell}{4}}$, by using \ref{muneigh} we conclude that \begin{equation}\label{est2}\mu^{*n}(H^{(\delta)})\leqslant\mu^{*n}(\Gamma_0)\leqslant e^{-\frac{n\ell}{4}}.\end{equation}

Since $n\geqslant\frac{\log(1+\eta)}{7C}\frac{\log{\frac{1}{\delta}}}{\log{\frac{1}{\varepsilon}}}=\frac{\log{\frac{1}{\delta}}}{7C\ell}$, by combining \ref{est1} and \ref{est2} we get that $$\mu^{*n}(H^{(\delta)})\leqslant e^{-\min\{c,\frac{1}{4}\}n\ell}\leqslant e^{-\frac{\min\{c,\frac{1}{4}\}}{7C}\log{\frac{1}{\delta}}}=\delta^{\frac{\min\{c,\frac{1}{4}\}}{7C}},$$

which proves Claim 6. \hfill$\square$

 Finally, put $d_1 = \frac{\min\{c,\frac{1}{4}\}}{7C}$ and $d_2 = \frac{\log(1+\eta)}{12C}.$
Then Claim 6 implies that $d_1,d_2 > 0$ satisfy the conclusion of Theorem \ref{escape}.
\hfill$\blacksquare$

\vskip 0.1in


\section{$\ell^2$-flattening} 
A key step in Bourgain and Gamburd's remarkable strategy \cite{BG05} for proving spectral gap is the so-called  {\bf $\ell^2$-flattening lemma}.  In \cite{BG06} and \cite{BG10}, Bourgain and Gamburd established a flattening lemma for probability measures on $SU(2)$  and $SU(d), d\geqslant 2$, respectively.   Bourgain and Yehudayoff then proved a flattening lemma for probability measures on $SL_2(\mathbb R)$ whose support is large but ``controlled" \cite{BY11}.  All of these results rely on {\bf product theorems} for the respective Lie groups. 
In an important recent development, de Saxc\'{e} obtained a product theorem for arbitrary connected simple Lie groups \cite{dS14}. 
This allowed Benoist and de Saxc\'{e} \cite{BdS14} to extend the flattening lemmas of \cite{BG06,BG10} to any compact connected simple Lie group. 

In this section, we first note that the product theorem of \cite{dS14} allows to derive a flattening lemma in the spirit of \cite[Lemma 4.1]{BY11} for arbitrary connected simple Lie groups.

\begin{lemma}[$\ell^2$-flattening, \cite{BdS14}]\label{BdS} Let $G$ be a connected simple Lie group with trivial center. 
Given $\alpha,\kappa>0$, there exist $\beta,\gamma>0$ such that the following holds for any $\delta>0$ small enough.

Suppose that $\mu$ is a symmetric Borel probability measure on $G$ such that
\begin{enumerate}
\item supp$(\mu)\subset B_{\delta^{-\beta}}(1)$,
\item $\|\mu*P_{\delta}\|_2\geqslant \delta^{-\alpha}$, and
\item $(\mu*\mu)(H^{(\rho)})\leqslant\delta^{-\gamma}\rho^{\kappa}$, for all $\rho\geqslant\delta$ and any proper closed connected subgroup $H<G$.
\end{enumerate}
Then $\|\mu*\mu*P_{\delta}\|_2\leqslant\delta^{\gamma}\|\mu*P_{\delta}\|_2.$

\end{lemma}
Lemma \ref{BdS} follows by adapting the proof of \cite[Lemma 2.5]{BdS14} in order to deal with non-compact Lie groups $G$ and measures $\mu$ with large controlled support (in the sense of (1)). Nevertheless, for completeness, we include the details of proof in the Appendix. 

 For now, we assume this lemma, and continue towards proving our main results.
More precisely, by applying Lemma \ref{BdS} repeatedly we obtain:

\begin{corollary}\label{flattening}
Let $G$ be a connected simple Lie group with trivial center, and $d_1,d_2 > 0$ be given. 
Then for every $\alpha>0$, there exist $\varepsilon_0>0$ and $c_0>0$ such the following holds. 

Let $0<\varepsilon<\varepsilon_0$ and $\mu$ be a Borel probability measure on $G$ such that 
supp$(\mu)\subset B_{\varepsilon}(1)$. Assume that
 for any $\delta>0$ small enough we have $\mu^{*2n}(H^{(\delta)})\leqslant\delta^{d_1}$, for any proper connected closed subgroup $H<G$, where  $n=\Big\lfloor{d_2\frac{\log{\frac{1}{\delta}}}{\log{\frac{1}{\varepsilon}}}}\Big\rfloor$.

Then for any $\delta>0$ small enough we have $\|\mu^{*n}*P_{\delta}\|_2\leqslant\delta^{-\alpha}$, for any integer $n\geqslant c_0\frac{\log{\frac{1}{\delta}}}{\log{\frac{1}{\varepsilon}}}.$
\end{corollary}

{\it Proof.} Let $\alpha>0$. By Lemma \ref{BdS} there are $\beta,\gamma>0$ such that for any $\delta>0$ small enough the following holds: if $\nu$ is a symmetric Borel probability measure on $G$ which satisfies

  \begin{itemize}
 \item [(a)] supp$(\nu)\subset B_{\delta^{-\beta}}(1)$, and 
\item [(b)] $(\nu*\nu)(H^{(\rho)})\leqslant\delta^{-\gamma}\rho^{\frac{d_1}{4}}$, for all $\rho\geqslant\delta$ and any proper closed connected subgroup $H<G$, \end{itemize}

then either $\|\nu*P_{\delta}\|_2\leqslant\delta^{-\alpha}$, or
 $\|\nu*\nu*P_{\delta}\|_2\leqslant\delta^{\gamma}\|\nu*P_{\delta}\|_2.$

We first claim that there is a constant $C>1$ depending only on $G$ such that the following holds. Let $\rho\in (0,1)$, $R>2$, $a,b,x\in B_R(1)$ and $h,k\in G$ such that $\|x^{-1}a-h\|_2\leqslant\rho$ and $\|x^{-1}b-k\|_2\leqslant\rho$. Then $\|b^{-1}a-k^{-1}h\|_2\leqslant R^C\rho$.
Indeed, the claim follows since there is a constant $c>1$ depending only on $G$ such that $\|y^{-1}\|_2\leqslant (\|y\|_2+1)^{c}$, for any $y\in G$, and we have that
\begin{align*}
\|b^{-1}a-k^{-1}h\|_2&=\|(x^{-1}b)^{-1}(x^{-1}a)-k^{-1}h\|_2\\&\leqslant\|x^{-1}\|_2\|a\|_2\|(x^{-1}b)^{-1}-k^{-1}\|_2+\|k^{-1}\|_2\|x^{-1}a-h\|_2\\&\leqslant \|x^{-1}\|_2\|a\|_2\|(x^{-1}b)^{-1}\|_2\|k^{-1}\|_2\|x^{-1}b-k\|_2+\|k^{-1}\|_2\|x^{-1}a-h\|_2.
\end{align*}

Let $k\geqslant 1$ be the smallest integer such that $\delta^{k\gamma}\|P_{\delta}\|_2\leqslant\delta^{-\alpha}$, for any $\delta>0$ small enough.  Let $\varepsilon>0$ small enough such that $\frac{2^kd_2}{4\log{\frac{1}{\varepsilon}}}<\frac{\min\{\beta,\frac{\gamma}{d_1C}\}}{\varepsilon}$. Let $\mu$ be a Borel probability measure on $G$ which is supported on $B_{\varepsilon}(1)$ and satisfies the hypothesis. The proof relies on the following:

{\bf Claim.} If $\delta>0$ is small enough and $n$ is an integer such that $\Big\lfloor{d_2\frac{\log{\frac{1}{\delta}}}{4\log{\frac{1}{\varepsilon}}}}\Big\rfloor\leqslant n\leqslant\min\{\beta,\frac{\gamma}{d_1C}\}\frac{\log{\frac{1}{\delta}}}{\varepsilon}$, then the measure $\nu=\mu^{*n}$  satisfies conditions (a) and (b).

{\it Proof of the claim.}
Since supp$(\nu)\subset B_{(1+\varepsilon)^n}(1)$ and $(1+\varepsilon)^n\leqslant [(1+\varepsilon)^{\frac{1}{\varepsilon}}]^{\beta\log{\frac{1}{\delta}}}<e^{\beta\log{\frac{1}{\delta}}}=\delta^{-\beta}$, we get that $\nu$ satisfies (a). 
To verify (b), let $\rho\geqslant\delta$ and $H<G$ be a proper closed connected subgroup. We may assume that $\rho\leqslant\delta^{\frac{4\gamma}{d_1}}$, because otherwise $\delta^{-\gamma}\rho^{\frac{d_1}{4}}>1$ and (b) is trivially satisfied.

Let $m=\Big\lfloor{d_2\frac{\log{\frac{1}{\rho^{\frac{1}{2}}}}}{\log{\frac{1}{\varepsilon}}}}\Big\rfloor=\Big\lfloor{d_2\frac{\log{\frac{1}{\rho}}}{2\log{\frac{1}{\varepsilon}}}}\Big\rfloor$. Then $m\leqslant 2n$ and the hypothesis implies that 
\[\mu^{*2m}(H^{\rho^{(\frac{1}{2})}})\leqslant\rho^{\frac{d_1}{2}}\].
For $x\in G$, denote $A_x=xH^{(\rho)}\cap$ supp$(\mu^{*m})$. Since $\mu^{*2n}=\mu^{*(2n-m)}*\mu^{*m}$, we have

\begin{equation}\label{nu_1}
\nu*\nu(H^{(\rho)})=\mu^{*2n}(H^{(\rho)})\leqslant\sup_{x\in \text{supp}(\mu^{*(2n-m)})}\mu^{*m}(xH^{(\rho)})=\sup_{x\in \text{supp}(\mu^{*(2n-m)})}\mu^{*m}(A_x).
\end{equation}

Further, since $\mu^{*m}$ is symmetric, Lemma \ref{A^{-1}A} implies that

\begin{equation}\label{nu_2}
\mu^{*m}(A_x)\leqslant\mu^{*2m}(A_x^{-1}A_x)^{\frac{1}{2}},\;\;\;\;\text{for any $x\in G$}.
\end{equation}

 Let $x\in$ supp$(\mu^{*(2n-m)})$ and $a,b\in A_x$.  Since supp$(\mu^{*k})\subset B_{(1+\varepsilon)^k}(1)$, for any $k\geqslant 1$, we have that $a,b,x\in B_{(1+\varepsilon)^{2n}}(1)$. By the definition of $A_x$, we can find $h,k\in H$ such that $\|x^{-1}a-h\|_2\leqslant\rho$ and $\|x^{-1}b-k\|_2\leqslant\rho$. The earlier claim implies that $\|b^{-1}a-k^{-1}h\|_2\leqslant (1+\varepsilon)^{2Cn}\rho$. 
 Since $n<\frac{\gamma}{d_1C}\frac{\log{\frac{1}{\delta}}}{\varepsilon}$, we get that $(1+\varepsilon)^{2Cn}<e^{2Cn\varepsilon}<\delta^{-\frac{2\gamma}{d_1}}\leqslant\rho^{-\frac{1}{2}}.$ Thus $\|b^{-1}a-k^{-1}h\|_2\leqslant\rho^{\frac{1}{2}}$. 
 
 Since $k^{-1}h\in H$ and $a,b\in A_x$ are arbitrary, we deduce that $A_x^{-1}A_x\subset H^{(\rho^{\frac{1}{2}})}$. 
 By combining \eqref{nu_1} and \eqref{nu_2} we therefore derive that $$\nu*\nu(H^{(\rho)})\leqslant\mu^{*2m}(H^{(\rho^{\frac{1}{2}})})^{\frac{1}{2}}\leqslant\rho^{\frac{d_1}{4}},$$
which finishes the proof of the claim.
\hfill$\square$

Let $\delta>0$ and put $n_0=\Big\lfloor{d_2\frac{\log{\frac{1}{\delta}}}{4\log{\frac{1}{\varepsilon}}}}\Big\rfloor$ and $n_1=2^kn_0$.
We claim that $\|\mu^{*n_1}*P_{\delta}\|_2\leqslant\delta^{-\alpha}$, for any small enough $\delta>0$.
Once this claim is proven, the conclusion follows for $c_0=2^{k-2}d_2$  since $n_1\leqslant 2^{k-2}d_2\frac{\log{\frac{1}{\delta}}}{\log{\frac{1}{\varepsilon}}}$ and $\|\mu^{*n}*P_{\delta}\|_2=\|\mu^{*(n-n_1)}*(\mu^{*n_1}*P_{\delta})\|_2\leqslant\|\mu^{*n_1}*P_{\delta}\|_2$, for any $n\geqslant n_1$.

 Assume by contradiction that the claim is false and let $0\leqslant i\leqslant k$. Then $2^in_0\leqslant n_1$ and therefore $\|\mu^{*2^in_0}*P_{\delta}\|_2\geqslant\|\mu^{*n_1}*P_{\delta}\|_2>\delta^{-\alpha}$. 
 On the other hand, 
 $\Big\lfloor{d_2\frac{\log{\frac{1}{\delta}}}{4\log{\frac{1}{\varepsilon}}}}\Big\rfloor\leqslant 2^in_0\leqslant\min\{\beta,\frac{\gamma}{d_1C}\}\frac{\log{\frac{1}{\delta}}}{\varepsilon}$. The claim implies that $\mu^{2^in_0}$ satisfies conditions (a) and (b).  As $\|\mu^{*2^in_0}*P_{\delta}\|_2>\delta^{-\alpha}$ we must have $$\|\mu^{*2^{i+1}n_0}*P_{\delta}\|_2\leqslant\delta^{\gamma}\|\mu^{*2^in_0}*P_{\delta}\|_2,\;\;\;\text{for every $0\leqslant i\leqslant k$.}$$

By combining these inequalities we deduce that $\|\mu^{*n_1}*P_{\delta}\|_2\leqslant\delta^{k\gamma}\|\mu^{*n_0}*P_{\delta}\|_2\leqslant\delta^{k\gamma}\|P_{\delta}\|_2\leqslant\delta^{-\alpha},$
which is a contradiction.
\hfill$\blacksquare$


\section{A mixing inequality}
The goal of this section is to prove an analogue for simple Lie groups of the well-known mixing inequality for quasirandom finite groups (see \cite[Proposition 1.3.7]{Ta15}). In the next section, we will combine this mixing inequality with Corollary \ref{flattening} and a Littlewood-Paley decomposition on simple Lie groups to deduce Theorem \ref{restricted}.

\begin{theorem}[mixing inequality]\label{rho}
Let $G$ be a connected simple Lie group with trivial center. Denote by $d$ the dimension of $G$, and let $B \subset G$ be a measurable set with compact closure.

Then there exist constants $a,b,\kappa>0$ such that for every $F\in L^2(B)$ with $\|F\|_2=1$, we have
$$\|f*F\|_2^{16d}\leqslant a\|P_{\delta}*F\|_2 + b\delta^\kappa,$$

for all $f\in L^2(G)$ with $\|f\|_2=1$ and all $0<\delta<1.$
\end{theorem}

This result and its proof are inspired by \cite[Lemma 10.35]{BG10}, which dealt with the case $G=SU(d)$, for $d\geqslant 2$. In particular, we borrow from \cite{BG10} the  idea of reducing to functions $F$ that satisfy an additional ``symmetry", i.e. are eigenvectors for a maximal torus of $G$.  This reduction is crucial, as it will allow us to exploit certain cancellations appearing in the integrals.

Turning to the proof of Theorem \ref{rho}, we start with a classical lemma which can be easily deduced from \cite[Section 2]{RS88}. We denote by $C_{\text{c}}^1(G)$  the space of compactly supported $C^1$-functions on $G$.

\begin{lemma}
\label{LieRC}
Let $G$ and $H$ be two Lie groups of dimensions $n$ and $m$. Assume that $n \geqslant m$. Consider an analytic function $\phi : G \to H$ such that the derivative $d\phi_x:\mathfrak{g} \to \mathfrak{h}$ has rank $m$,  at almost every point $x \in G$. Let $\psi \in C_{\text{c}}^1(G)$ and denote by $\mu=\phi_*(\psi \cdot dm_G)$ the push-forward measure of the measure $\psi \cdot dm_G$ on $G$ through $\phi$.

Then $\mu$ is absolutely continuous with respect to $m_H$, and the Radon-Nykodym derivative $\rho:H \to \R$ is {\it$L^1$-H\"older}: there exist $\alpha > 0$ and $C > 0$ such that
\[\int_H \vert \rho(g^{-1}h) - \rho(h) \vert\text{d}h \leqslant C\Vert g - 1\Vert_2^\alpha, \;\;\;\;\text{for every $g \in H$}.\]
\end{lemma}

By applying this lemma, we obtain the following:

\begin{lemma}
\label{abs}
Let $G$ be a connected simple Lie group and $H<G$ be a connected compact Lie subgroup of dimension $1$.  Define $\pi:G\times H^2\rightarrow G$ by letting $\pi(g,t_1,t_2)=t_1gt_1^{-1}t_2g^{-1}t_2^{-1}$, for all $g\in G, t_1,t_2\in H$.  Let $\psi \in C_{\text{c}}^1(G)$ and define $\nu=\pi_*((\psi \cdot dm_G)\times m_{H} \times m_H)$.

Then $\nu^{*n}$ is absolutely continuous with respect to $m_G$, and the corresponding Radon-Nykodym derivative is $L^1$-H\"older, for every integer $n\geqslant$ dim$(G)$.
\end{lemma}

{\it Proof.} 
Let $n\geqslant 1$.
Then $\nu^{*n}=\pi^{(n)}_*({(\psi \cdot dm_G)}^n\times m_H^{2n})$, where $\pi^{(n)}:G^n\times H^{2n}\rightarrow G$ is defined as 
\[\pi^{(n)}(g_1,...,g_n,t_1,...,t_{2n})=\prod_{i=1}^n\pi(g_i,t_{2i-1},t_{2i})=\prod_{i=1}^n(t_{2i-1}g_it_{2i-1}^{-1}t_{2i}g_i^{-1}t_{2i}^{-1}).\]
By Lemma \ref{LieRC}, we only have to check that the derivative of the analytic function $\pi^{(n)}$ has rank $d$, at almost every point, as soon as $n \geqslant d := \dim(G)$.

Fix $n \geqslant d$. Let $\frak g$ and $\frak h$ be  the Lie algebras of $G$ and $H$, respectively. Let Ad$:G\rightarrow GL(\frak g)$ be the adjoint representation of $G$.  Since dim$(H)=1$ and $H$ is connected, there is $b\in\frak g$ such that $\frak h=\{ub \, \vert \, u \in \R\}$ and $H=\{\exp(ub) \, \vert \, u \in \R\}$.

Let $X_n$ be the set of  $(g_1,...,g_n,t_1,...,t_{2n})\in G^n\times H^{2n}$ such that the following set spans $\frak g$: 
\[\{\text{Ad}(\prod_{j=1}^{i-1}\pi(g_j,t_{2j-1},t_{2j}))(b)-\text{Ad}((\prod_{j=1}^{i-1}\pi(g_j,t_{2j-1},t_{2j}))t_{2i-1}g_it_{2i-1}^{-1})(b)\;|\; 1\leqslant i\leqslant n\}.\]

{\bf Claim 1.}
 $\rk(d(\pi^{(n)})_x)=d$, for every $x\in X_n$.

{\it Proof of Claim 1.} Take $x=(g_1,...,g_n,t_1,...,t_{2n})\in X_n$. Proving the claim amounts to showing that the map $\tilde \pi_n : y \mapsto \pi^{(n)}(y)\pi^{(n)}(x)^{-1}$ is such that $d(\tilde \pi_n)_x$ has rank $d$. For all $1 \leqslant i \leqslant n$, define a map $\varphi_i : \R \to G$ by the formula
\[\varphi_i(u)= \tilde \pi_n(g_1,...,g_n,t_1,...,t_{2i - 2},\exp(ub)t_{2i-1},t_{2i} \cdots, t_{2n})\]
The derivative $\varphi_i'(0) \in \frak g$ belongs to the range of the derivative $d(\tilde \pi_n)_x$, while an easy computation gives that 
 \[\varphi_i'(0) = \text{Ad}(\prod_{j=1}^{i-1}\pi(g_j,t_{2j-1},t_{2j}))(b)-\text{Ad}((\prod_{j=1}^{i-1}\pi(g_j,t_{2j-1},t_{2j}))t_{2i-1}g_it_{2i - 1}^{-1})(b).\]
 Since $x\in X_n$, the set $\{\varphi_i'(0) \, \vert \, 1 \leqslant i \leqslant n \}$ spans $\frak g$, and $d(\tilde \pi_n)_x$ is therefore onto.
 \hfill$\square$
 
 {\bf Claim 2.} $X_n$ is a nonempty Zariski open subset of $G^n\times H^{2n}$, for every $n\geqslant d$.
 
 {\it Proof of Claim 2.} Since $X_n$ is clearly Zariski open, for every $n\geqslant 1$, it remains to argue that $X_n$ is nonempty, whenever $n\geqslant d$. Since $G$ is simple, $\frak g$ is the only non-trivial Ad$(G)$-invariant subspace of $\frak g$. Thus, the span of $\{\text{Ad}(g)(b)-\text{Ad}(h)(b)|g,h\in G\}$ is equal to $\frak g$. Equivalently, we derive that the span of $\{\text{Ad}(g)(b)-b|g\in G\}$ is also equal to $\frak g$. 
 We can therefore find $g_1,...,g_d\in G$ such that $\{\text{Ad}(g_i)(b)-b|1\leqslant i\leqslant d\}$ spans $\frak g$. Define $g_{d+1}=....=g_n=t_1=...=t_{2n}=1$. Then it is clear that $(g_1,...,g_n,t_1,...,t_{2n})\in X_n$, which shows that $X_n$ is nonempty, as claimed. \hfill$\square$
 
Finally, if $n\geqslant d$, then Claim 2 implies that $X_n$ a co-null subset of $G^n\times H^{2n}$. 
\hfill$\blacksquare$

\vskip 0.1in
We are now ready to prove Theorem \ref{rho}.

{\it Proof of Theorem \ref{rho}.} Let $F\in L^2(B)$ and $f\in L^2(G)$ with $\|f\|_2=1$. Since $F*\check{F}$ is supported on $BB^{-1}$, we have
 $\|f*F\|_2^2=\langle f*F,f*F\rangle=\langle\check{f}*f,F*\check{F}\rangle\leqslant \|\check{f}*f\|_{2,BB^{-1}}\|F*\check{F}\|_2$. Since  $\|\check{f}*f\|_{\infty}\leqslant 1$, we  get that $\|\check{f}*f\|_{2,BB^{-1}}\leqslant |BB^{-1}|^{1/2}$. Moreover, for every $g\in G$, we have that $F*\check{F}(g)=\int_G\overline{F(g^{-1}x)}F(x)\;\text{d}x=\overline{\langle\lambda_g(F),F\rangle}.$
By putting these facts together, we get that $$\|f*F\|_2^{16d}\leqslant |BB^{-1}|^{4d}(\int_G|\langle\lambda_g(F),F\rangle|^2\;\text{d}g)^{4d}.$$
Thus, the conclusion reduces to proving the following:

$(*)$ there exist constants $a,b,\kappa>0$ such that for every $F\in L^2(B)$ with $\|F\|_2=1$, we have that $$\big(\int_G|\langle\lambda_g(F),F\rangle|^2\;\text{d}g\big)^{4d}\leqslant a\|P_{\delta}*F\|_2+b\delta^{\kappa},\;\;\;\;\text{for all $0<\delta<1$.}$$

To this end, we fix a compact connected Lie subgroup $H$ of $G$ with dimension $1$.
Below, we denote by $x,y,z,g$ elements of $G$ and by $t,t_1,t_2$ elements of $H$.
Writing $\text{d}x$ (respectively, $\text{d}t$) will refer to integration against the Haar measure of $G$ (respectively, $H$).

Let $\widetilde B\subset G$ be an open set with compact closure which contains $B^{-1}B$. Let $\psi \in C_{\text{c}}^1(G)$  be a non-negative function which is equal to $1$ on $\widetilde B$.
Define $\pi:G\times H^2\rightarrow G$ by  $\pi(x,t_1,t_2)=t_1xt_1^{-1}t_2x^{-1}t_2^{-1}$, for all $x\in G$ and $t_1,t_2\in H$. Let $\nu=\pi_*((\psi \cdot dm_G) \times m_H^2)$. Lemma \ref{abs} implies that $\nu^{*d}$ is absolutely continuous with respect to $m_G$ and the corresponding Radon-Nykodym derivative $\rho$ is $L^1$-H\"older. In other words, there exist $\kappa > 0$ and $C > 0$ such that
\[\int_G \vert \rho(g^{-1}h) - \rho(h) \vert dh \leqslant C\Vert g - 1\Vert_2^{2\kappa}, \;\; \forall g \in G.\]
For $x\in G$, we define an operator $R_x:L^2(G)\rightarrow L^2(G)$ by the formula 
\[(R_xf)(z)=\int_{H}f(ztxt^{-1})\;\text{d}t,  \; f\in L^2(G), \; z\in G.\] 

{\bf Claim 1}. For every $f\in L^2(G)$, we have $\displaystyle{\int_{\widetilde{B}}\|R_x(f)\|_2^2\;\text{d}x} \leqslant \Vert f*\nu \Vert_2\Vert f \Vert_2$.
\begin{proof}[Proof of Claim 1]
 Let $f\in L^2(G)$. 
Since $(R_x^*R_xf)(z) =\displaystyle{\int_{H^2}}f(zt_1x^{-1}t_1^{-1}t_2xt_2^{-1})\;\text{d}t_1\;\text{d}t_2,$ the claim follows from the following calculation 
\begin{align*}
\int_{\widetilde{B}}\|R_x(f)\|_2^2\;\text{d}x \leqslant \int_{G}\|R_x(f)\|_2^2\psi(x)\;\text{d}x & = \int_{G} \langle R_x^*R_x(f),f\rangle\psi(x)\; \text{d}x\\
&= \int_{G}\Big(\int_{G \times H^2}f(zt_1x^{-1}t_1^{-1}t_2xt_2^{-1})\; \psi(x)\;\text{d}x\;\text{d}t_1\;\text{d}t_2\Big)\overline{f(z)}\; \text{d}z\\
& =\int_{G}(f*\nu)(z)\overline{f(z)} \; \text{d}z \leqslant \Vert f*\nu \Vert_2\Vert f \Vert_2. \qedhere
\end{align*}
\end{proof}
Next, using that $\rho$ is $L^1$-H\"older, we deduce the following claim:

{\bf Claim 2}. There is $c>0$ such that $\|P_{\delta}*f*\rho-f*\rho\|_2\leqslant c\delta^{\kappa}\|f\|_2$, for all $f\in L^2(B)$ and  $0<\delta<1$.

{\it Proof of Claim 2}.
Take $f \in L^2(B)$ and $\delta > 0$. Note that for $x \in G$, we have
\[(P_{\delta}*f*\rho-f*\rho)(x) = \frac{1}{\vert B_\delta\vert} \int_{B_\delta(1) \times B} f(z) (\rho(z^{-1}y^{-1}x) - \rho(z^{-1}x)) \; \d y \; \d z.\]
Using the Cauchy-Schwarz inequality and $L^1$-H\"older condition for $\rho$, we get that $\Vert P_{\delta}*f*\rho-f*\rho\Vert_2^2 $ is at most equal to
\begin{align*}
\frac{1}{\vert B_\delta\vert^2}\int_G \Big(\int_{B_\delta(1) \times B}  \vert & f(z) \vert^2 \vert \rho(z^{-1}y^{-1}x) - \rho(z^{-1}x)\vert \; \d y \; \d z\Big)\Big(\int_{B_\delta(1) \times B} \vert \rho(z^{-1}y^{-1}x) - \rho(z^{-1}x)\vert \; \d y \; \d z\Big) \; \d x\\
& \leqslant \frac{2\Vert \rho \Vert_1}{\vert B_\delta\vert} \int_{G \times B_\delta(1) \times B} \vert f(z) \vert^2 \vert \rho(z^{-1}y^{-1}x) - \rho(z^{-1}x)\vert \; \d x \; \d y \; \d z\\
& \leqslant \frac{2\Vert \rho \Vert_1}{\vert B_\delta\vert}  \int_B \vert f(z) \vert^2 \Big( \int_{G \times B_\delta(1)} \vert \rho(z^{-1}y^{-1}zx) - \rho(x)\vert \; \d x \; \d y\Big) \; \d z\\
&  \leqslant 2C \Vert \rho \Vert_1  \int_B \vert f(z) \vert^2 \sup_{y \in B_\delta(1)} \Vert z^{-1}yz - 1 \Vert_2^{2\kappa} \; \d z \; \leqslant \; c^2\delta^{2\kappa} \Vert f \Vert_2^2,
\end{align*}
for some constant $c > 0$ independent of $f$ and $\delta$.
\hfill$\square$

\vskip 0.05in
Let $F\in L^2(B)$ with $\|F\|_2=1$.
The proof of $(*)$ splits into two cases.

{\bf Case 1}. 
We first prove assertion $(*)$ in the following case:
there is a character $\eta:H\rightarrow\mathbb T$ such that for all $t \in H$, $F(xt)=\eta(t)F(x)$, for almost every $x\in G$.

Then for almost every $(x,y,t)\in G^2\times H$ we have that $F(xt)\overline{F(yt)}=F(x)\overline{F(y)}$. 
By using this fact we get that 
\begin{align*}
\int_{G}|\langle\lambda_g(F),F\rangle|^2\;\text{d}g & = \int_{G^3}F(g^{-1}x)\overline{F(x)}\,\overline{F(g^{-1}y)}F(y)\;\text{d}g\;\text{d}x\;\text{d}y\\
& = \int_{G^3}\int_{H}F(g^{-1}xt)\overline{F(x)}\,\overline{F(g^{-1}yt)}F(y)\;\text{d}g\;\text{d}x\;\text{d}y\;\text{d}t
\end{align*}
Using left invariance of the Haar measure and unimodularity on the $g$ and $y$ variables, we get
\begin{align*}
\int_{G}|\langle\lambda_g(F),F\rangle|^2\;\text{d}g & = \int_{G^3}\int_{H}F(gt^{-1}y^{-1}xt)\overline{F(x)}\,\overline{F(g)}F(y) \;\text{d}g\;\text{d}x\;\text{d}y\;\text{d}t\\
& = \int_{G^3}\int_{H}F(gt^{-1}y^{-1}t)\overline{F(x)}\,\overline{F(g)}F(xy) \;\text{d}g\;\text{d}x\;\text{d}y\;\text{d}t\\
& = \int_{G} (\check{F}*F)(y)\; \langle R_{y^{-1}}F,F\rangle\;\text{d}y = \int_{G} (\check{F}*F)(y)\; \langle F,R_yF\rangle\;\text{d}y
\end{align*}

Since $\|\check{F}*F\|_{\infty}\leqslant 1$ and $\check{F}*F\in L^2(\widetilde{B})$, we conclude from Claim 1 and Lemma \ref{powers} that 
\begin{align*}
\int_{G}|\langle\lambda_g(F),F\rangle|^2\;\text{d}g\leqslant \int_{\widetilde{B}}\|R_x(F)\|_2\;\text{d}x &\leqslant |\widetilde{B}|^{1/2}\Big(\int_{\widetilde B}\|R_x(F)\|_2^2\;\text{d}x\Big)^{1/2}\\& \leqslant |\widetilde{B}|^{\frac{1}{2}}\|F*\nu\|_2^{\frac{1}{2}}\leqslant  |\widetilde{B}|^{\frac{1}{2}}\|\nu\|^{\frac{1}{4}}\|F*\nu^{*d}\|_2^{\frac{1}{4d}}\\&=|\widetilde{B}|^{\frac{1}{2}}\|\nu\|^{\frac{1}{4}}\|F*\rho\|_2^{\frac{1}{4d}}.
\end{align*}
On the other hand, Claim 2 yields
\[\|F*\rho\|_2\leqslant \big(\|P_{\delta}*F*\rho\|_2+\|F*\rho-P_{\delta}*F*\rho\|_2\big)\leqslant \|\rho\|_1\|P_{\delta}*F\|_2+c\delta^\kappa.\]
Thus, if we let $a=|\widetilde{B}|^{2d}\|\nu\|^{d}\|\rho\|_1$ and $b=c|\widetilde{B}|^{2d}\|\nu\|^d$, the desired inequality $(*)$ follows in Case 1. Moreover, notice the crucial fact that $a$ and $b$ are independent of the character $\eta$.

{\bf Case 2.} We now prove $(*)$ for an arbitrary function $F\in L^2(B)$ with $\|F\|_2=1$. 

Consider the unitary representation $H \actson^{\sigma} L^2(G)$ corresponding to the right multiplication action $H \actson G$.
Since $H$ is compact and abelian,  we can decompose
\[L^2(G)=\underset{\eta\in\text{Char}(H)}\bigoplus\mathcal H_{\eta},\]
where $\mathcal H_{\eta}$ denotes the eigenspace of $\sigma$ corresponding to a character $\eta:H\rightarrow\mathbb T$.
 
Thus, we can decompose $F = \sum_\eta F_{\eta}$, where $F_\eta(xt)=\eta(t)F(x)$, for almost every $x\in G$, $t \in H$. Note that the functions $F_\eta$ do not necessarily belong to $L^2(B)$. However, $F_\eta$ belongs to the closure of the linear span of $\sigma(H)F$, and therefore to $L^2(BH)$, for every $\eta$.

By applying Case 1 with $BH$ instead of $B$,  and using homogeneity, we get that there exist constants $a,b,\kappa>0$ (independent of $F$) such that for all $\eta \in \Char(H)$ we have
\[(\int_{G}|\langle\lambda_g(F_\eta),F_\eta\rangle|^2\;\text{d}g)^{4d} \leqslant a\|P_{\delta}*F_\eta\|_2\Vert F_\eta \Vert_2^{16d-1} + b\delta^\kappa \Vert F_\eta \Vert_2^{16d}.\]
Since all the norms on $\R^2$ are equivalent, we find $a' , b' > 0$ (only depending on $a,b,d$) such that for all $\eta\in \Char(H)$ we have that
\[(\int_{G}|\langle\lambda_g(F_\eta),F_\eta\rangle|^2\;\text{d}g)^{1/2} \leqslant a'\|P_{\delta}*F_\eta\|_2^{1/8d} \Vert F_\eta\Vert_2^{2 - 1/8d} + b'\delta^{\kappa/8d} \Vert F_\eta \Vert_2^{2}.\]
But since $\lambda_g(F_\eta) \in \cH_\eta$ and $P_\delta \ast F_\eta \in \cH_\eta$ for all $\eta$, by using the triangle inequality for $\|.\|_2$ and H\"older's inequality we get that
\begin{align*}
\Big(\int_{G}|\langle\lambda_g(F),F\rangle |^2\;\text{d}g\Big)^{1/2} & = \Big(\int_{G} \Big\vert \sum_\eta \langle\lambda_g(F_\eta),F_\eta\rangle \Big\vert^2\;\text{d}g\Big)^{1/2}\\ &\leqslant\sum_{\eta} \Big(\int_{G}|\langle\lambda_g(F_\eta),F_\eta\rangle|^2\;\text{d}g\Big)^{1/2}\\
& \leqslant \Big(\sum_\eta a'\|P_{\delta}*F_\eta\|_2^{1/8d} \Vert F_\eta\Vert_2^{2 - 1/8d}\Big) + b'\delta^{\kappa/8d}\\
& \leqslant a' \Big( \sum_\eta \Vert P_\delta \ast F_\eta \Vert_2^2 \Big)^{1/16d} \Big( \sum_\eta \Vert F_\eta \Vert_2^2 \Big)^{1 - 1/16d} + b'\delta^{\kappa/8d}\\
& = a' \Vert P_\delta \ast F \Vert_2^{1/8d} + b'\delta^{\kappa/8d}.
\end{align*}
Using again the equivalence of norms in $\R^2$ and modifying the values of $a$ and $b$ if necessary, the conclusion follows. \hfill$\blacksquare$


\section{Proofs of Theorem \ref{restricted} and Corollary \ref{by}}

\subsection{A Littlewood-Paley decomposition on Lie groups} 
Let $G$ be a connected simple Lie group with trivial center. 
In order to prove Theorem \ref{restricted}, we next introduce a Littlewood-Paley decomposition on $G$.  This is analogous to the Littlewood-Paley decomposition on $G=SU(d)$ defined by  Bourgain and Gamburd  in \cite[Section 10]{BG10}. 
As before, we endow $G$ with the $\|.\|_2$ metric and denote by $\mathcal C(G)$ the family of measurable subsets of $G$ with compact closure.

We define bounded linear operators $\Delta_i:L^2(G)\rightarrow L^2(G)$, $i\geqslant 0$, as follows
\begin{align*}
\Delta_0(F) & = P_{1/2} \ast F\\
\Delta_i(F) & = P_{2^{-(i+1)}} \ast F - P_{2^{-i}} \ast F, \text{ for all } i \geqslant 1.
\end{align*}

\begin{remark}
The decomposition $F=\sum_{i\geqslant 0}\Delta_i(F)$ is analogous to the classical Littlewood-Paley decomposition on $\mathbb R^n$, in the following sense. 
For any $i\geqslant 0$, the function $\Delta_i(F)$ ``lives" at scale $2^{-i}$: it is essentially constant at scales $\ll 2^{-i}$ and essentially has mean zero on balls of radius $\gg 2^{-i}$. 
\end{remark}

We now prove that the operators $\Delta_i, i\geqslant 0,$ yield an almost orthogonal decomposition of $L^2(G)$. This will allow us to reduce to functions living at an arbitrary small scale in the proof of restricted spectral gap Theorem \ref{restricted}.

\begin{theorem}\label{L-P}
There exists a constant $C > 0$ such that for all $F\in L^2(G)$ and any  $\mu\in\mathcal M(G)$ with $\supp(\mu) \subset B_1(1)$, we have that
\begin{enumerate}
\item $\sum_{i \geqslant 0} \|\Delta_i(F)\|_2^2\leqslant C\|F\|_2^2$.
\item $\|\mu*F\|_2^2\leqslant C\sum_{i\geqslant 0}\|\mu*\Delta_i(F)\|_2^2$.
\item $\sum_{i \geqslant 0} 2^{i/2}\Vert P_{2^{-2i}} \ast \Delta_i(F) - \Delta_i(F) \Vert_2^2 \leqslant C\sum_{i \geqslant 0} \Vert \Delta_i(F)\Vert_2^2$.
\item $\sum_{i \geqslant 0} 2^{i/2}\Vert P_{2^{-i/2}} \ast \Delta_i(F) \Vert_2^2 \leqslant C\sum_{i \geqslant 0} \Vert \Delta_i(F)\Vert_2^2$.
\end{enumerate}
\end{theorem}

The first ingredient of the proof of Theorem \ref{L-P} is the following lemma. This lemma and its proof are a variation of \cite[Lemma 11]{KS71} due to Knapp and Stein.

\begin{lemma}[Cotlar-Stein]
\label{cotlar-stein}
Consider a Hilbert space $\mathcal H$ and bounded operators $T_i : \mathcal H \to\mathcal H$, $i \geqslant 0$. Assume that there exists $\varphi : \Z \to \R_+$ with $\Phi := \sum_{n \in \Z} \varphi(n) < \infty$ such that for all $i,j \geqslant 0$, we have $\Vert T_j^*T_i \Vert^{1/2} \leqslant \varphi(j-i)$ and $\Vert T_iT_j^* \Vert^{1/2} \leqslant \varphi(i-j)$. For $k\geqslant 0$, denote  $\Phi_k:= \sum_{\vert n \vert \geqslant k} \varphi(n)$.

Then for all $\xi \in\mathcal H$ and all $k \geqslant 0$ we have  
\[ \sum_{i,j : \, \vert i - j \vert \geqslant k} \vert \langle T_i\xi,T_j\xi \rangle \vert \, \leqslant \,  \Phi_k \Phi \Vert \xi \Vert^2.\]
\end{lemma}

{\it Proof.}
Fix $\xi \in\mathcal H$ and $k \geqslant 0$. For every $i,j \geqslant 0$, we choose a scalar $\alpha_{i,j}$ in such a way that $\vert \langle T_i\xi,T_j\xi \rangle \vert= \alpha_{i,j}  \langle T_i\xi,T_j\xi \rangle$, and that $\alpha_{i,j} = 0$ whenever $\vert \langle T_i\xi,T_j\xi \rangle \vert=0$.
Then for all $N \geqslant 0$, the operator $R_N := \sum_{0 \leqslant i,j \leqslant N: \; \vert i - j \vert \geqslant k} \alpha_{i,j}T_j^*T_i$ is self-adjoint.

In order to prove the lemma, it is sufficient to check that the operator norm of $R_N$ is at most $\Phi_k \Phi$,  for all $N\geqslant 0$.
Take $N \geqslant 0$. Since $R_N$ is self-adjoint, $\Vert R_N \Vert^p = \Vert R_N^p \Vert$,  for all integers $p\geqslant 1$. This leads to the estimate:
\[\Vert R_N \Vert^p  \leqslant \sum_{0 \leqslant i_1,j_1,\dots,i_p,j_p \leqslant N: \, \vert i_l - j_l \vert \geqslant k, \, \forall 1\leqslant l\leqslant p} \Vert T_{j_1}^*T_{i_1}T_{j_2}^*T_{i_2}\cdots T_{j_p}^*T_{i_p}\Vert.\]

Since 
 the general term of this sum is bounded by the following two quantities
\begin{align*}
\Vert T_{j_1}^*T_{i_1}T_{j_2}^*T_{i_2}\cdots T_{j_p}^*T_{i_p}\Vert & \leqslant \Vert T_{j_1}^*T_{i_1}\Vert\Vert T_{j_2}^*T_{i_2}\Vert \cdots \Vert T_{j_p}^*T_{i_p}\Vert\;\;\;\text{and}\\
\Vert T_{j_1}^*T_{i_1}T_{j_2}^*T_{i_2}\cdots T_{j_p}^*T_{i_p}\Vert & \leqslant \Vert T_{j_1}^* \Vert \Vert T_{i_1}T_{j_2}^* \Vert \cdots  \Vert T_{i_{p-1}}T_{j_p}^* \Vert \Vert T_{i_p}\Vert,
\end{align*}
 we get that
\begin{align*}
\Vert R_N \Vert^p & \leqslant \sum_{0 \leqslant i_1,j_1,\dots,i_p,j_p \leqslant N: \, \vert i_l - j_l \vert \geqslant k, \, \forall 1\leqslant l\leqslant p} (\Vert T_{j_1}^* \Vert \Vert T_{j_1}^*T_{i_1}\Vert \Vert T_{i_1}T_{j_2}^*\Vert \cdots \Vert T_{j_p}^*T_{i_p}\Vert \Vert T_{i_p} \Vert)^{1/2}\\
& \leqslant N(\max_{0 \leqslant i \leqslant N} \Vert T_i \Vert)(\max_{1 \leqslant j \leqslant N}\sum_{1 \leqslant i \leqslant N: \, \vert i-j \vert \geqslant k} \Vert T_j^*T_i\Vert^{1/2})^p(\max_{0 \leqslant i \leqslant N}\sum_{1 \leqslant j \leqslant N} \Vert T_iT_j^*\Vert^{1/2})^{p-1}\\
& \leqslant N(\max_{0 \leqslant i \leqslant N} \Vert T_i \Vert)\Phi_k^p \Phi^{p-1}.
\end{align*}
Since $p\geqslant 1$ is arbitrary, we indeed get that $\Vert R_N \Vert \leqslant \Phi_k\Phi$.
\hfill$\blacksquare$

\begin{remark}
The case $k=0$ of Lemma \ref{cotlar-stein} recovers the classical Cotlar-Stein lemma (see \cite[Chapter VII]{St93}) which asserts that, under the same assumptions as above, the sum $\sum_{i \geqslant 0} T_i$ converges in the strong operator topology. Lemma \ref{cotlar-stein} also implies that the sum $\sum_{i \geqslant 0} \Vert T_i\xi\Vert^2$ is finite, for all $\xi \in\mathcal H$. Later on, we will use the following inequalities, which follow easily from Lemma \ref{cotlar-stein}
\begin{equation}\label{CS1}
\sum_{i \geqslant 0} \Vert T_i\xi\Vert^2 \leqslant \Phi^2\Vert \xi \Vert^2,
\end{equation}
\begin{equation}\label{CS2}
\Vert \sum_{i \geqslant 0} T_i\xi\Vert_2^2 \leqslant k \sum_{i \geqslant 0} \Vert T_i\xi\Vert^2 + \Phi_k\Phi\Vert \xi \Vert^2, \;\; \text{for all}\;\; k \geqslant 0.
\end{equation}
\end{remark}

In order to prove Theorem \ref{L-P} we will also need the following lemma which allows us to quantify the ``orthogonality" between the operators $\Delta_i$, $i\geqslant 0$.
For a Borel  probability measure $\mu\in\mathcal M(G)$  we denote by $T_{\mu}:L^2(G)\rightarrow L^2(G)$ the contractive operator given by $T_{\mu}(F)=\mu*F$.

\begin{lemma} 
\label{almostorthogonal}
There exists a constant $C_0>0$ such that for any Borel probability measure $\mu\in\mathcal M(G)$ with $\supp(\mu) \subset B_1(1)$ we have that  
\begin{align*}\|(P_{\delta_1}-P_{\delta_2}) \ast \mu \ast P_{\delta_3}\|_1\leqslant C_0\frac{\delta_2}{\delta_3},\;\;\text{for all}\;\; 0<\delta_1\leqslant\delta_2\leqslant\delta_3<1,\;\;\text{and} \\
\Vert \Delta_j^*T_\mu^*T_\mu \Delta_i \Vert \leqslant \frac{C_0}{2^{\vert i - j \vert}} \; \text{ and } \; \Vert T_\mu \Delta_i  \Delta_j^* T_\mu^* \Vert\leqslant \frac{C_0}{2^{\vert i - j \vert}},\;\;\text{for all}\;\; i,j\geqslant 0.
\end{align*}
\end{lemma}

{\it Proof.}
Denote by $B_1$, $B_2$ and $B_3$ the balls centered at $1$ with respective radii $\delta_1$, $\delta_2$ and $\delta_3$.
Note that $\Vert (P_{\delta_1}-P_{\delta_2}) \ast \mu \ast P_{\delta_3}\Vert_1 \leqslant{\int_{G} \Vert (P_{\delta_1}-P_{\delta_2}) \ast \delta_x \ast P_{\delta_3}\Vert_1\;\text{d}\mu(x)}$. So it suffices to prove that lemma for Dirac measures $\mu = \delta_x$, with $x \in B_1(1)$.

Fix $x \in B_1(1)$. Then for all $y \in G$, we have $\vert(P_{\delta_1}-P_{\delta_2}) \ast \delta_x \ast P_{\delta_3}(y)\vert \leqslant \|P_{\delta_1}-P_{\delta_2}\|_1/ |B_3| \leqslant 2/ |B_3|$. Let us now bound the measure of the support of $(P_{\delta_1}-P_{\delta_2}) \ast \delta_x \ast P_{\delta_3}$. One easily checks that this support is contained in $B_2xB_3 \cap B_2x(G \setminus B_3)$.

Firstly, if $y \in B_2xB_3$, we write $y = axb$, where $a \in B_2$, $b \in B_3$. Then we have that $\Vert y \Vert_2 \leqslant 8$ and
\begin{align*}
\Vert x^{-1}y - 1 \Vert_2 \leqslant \Vert x^{-1}y - b \Vert_2 + \delta_3
& \leqslant \Vert x^{-1} \Vert_2 \Vert y - xb \Vert_2 + \delta_3\\
& \leqslant \Vert x^{-1} \Vert_2 \Vert xb \Vert_2 \delta_2 + \delta_3 \leqslant C_1 \delta_2 + \delta_3,
\end{align*}
where $C_1>0$ is independent of $x \in B_1(1)$. 

Secondly, if $y\in B_2x(G\setminus B_3)$, we write $y = a'xb'$, where $a'\in B_2$ and $b' \notin B_3$. Then we see that $\Vert xb' \Vert_2 \leqslant \Vert a^{-1} \Vert_2 \Vert y \Vert_2 \leqslant 8\Vert a^{-1} \Vert_2$ and
\begin{align*}
\Vert x^{-1}y - 1 \Vert_2  \geqslant \Vert b' - 1 \Vert_2 - \Vert x^{-1}y - b' \Vert_2 & \geqslant \delta_3 - \Vert x^{-1} \Vert_2\Vert y - xb' \Vert_2\\
& \geqslant \delta_3 - \Vert x^{-1} \Vert_2\Vert xb' \Vert_2\delta_2 \geqslant  \delta_3 - C_2\delta_2,
\end{align*}
where $C_2 > 0$ is independent of $x\in B_1(1)$. 

Therefore, the support of $(P_{\delta_1}-P_{\delta_2}) \ast \delta_x \ast P_{\delta_3}$ is contained in $x (B_{\delta_3+C_1\delta_2}\setminus B_{\delta_3-C_2\delta_2})$.
Altogether, we get that \[\|(P_{\delta_1}-P_{\delta_2}) \ast \delta_x \ast P_{\delta_3}\|_1\leqslant{\frac{2\;|B_{\delta_3+C_1\delta_2}\setminus B_{\delta_3 - C_2\delta_2}|}{|B_{\delta_3}|}},\]
which implies the first inequality. 
By using the fact that $\Vert x^{-1} - 1 \Vert_2 \leqslant \Vert x^{-1}\Vert_2 \Vert x - 1\Vert_2$ and arguing similarly to the above, it follows that the quantities $\|(\check{P}_{\delta_1} - \check{P}_{\delta_2}) \ast \mu \ast P_{\delta_3}\|_1$ and $\|(P_{\delta_1}-P_{\delta_2}) \ast \mu \ast \check{P}_{\delta_3}\|_1$ are bounded above by $C_0\delta_2 / \delta_3$, for some possibly larger constant $C_0>0$. Since $\|f*g\|_2\leqslant \|f\|_1\|g\|_2$, for all $f \in L^1(G)$, $g \in L^2(G)$, these estimates imply the rest of the asserted inequalities.\hfill$\blacksquare$

\bigskip
{\it Proof of Theorem \ref{L-P}}.
Let $C_0>0$ be the constant provided by Lemma \ref{almostorthogonal} and define $\varphi: \Z \to \R_+$ by letting $\varphi(n) = \frac{C_0^{1/2}}{2^{\vert n \vert/2}}$. Then Lemma \ref{almostorthogonal} gives that for any finitely supported probability measure on $G$ with supp$(\mu)\subset B_1(1)$, the operators $T_i := T_\mu \Delta_i$ on $L^2(G)$ satisfy $\Vert T_j^*T_i \Vert^{1/2} \leqslant \varphi(j-i)$ and $\Vert T_iT_j^* \Vert^{1/2} \leqslant \varphi(i-j)$, for all $i,j \geqslant 0$.

Let $\Phi$ and $\Phi_k$ be as defined  in Lemma \ref{cotlar-stein} and take $k$ large enough so that $\Phi_k\Phi < 1$.
Let $F \in L^2(G)$. Since $\lim\limits_{\delta\rightarrow 0}\|P_{\delta}*F-F\|_2=0$, we get that $F = \sum_{i \geqslant 0} \Delta_i(F)$. By combining this fact with equations \eqref{CS1} and \eqref{CS2}, we derive that 
\begin{equation}\label{eqsum}
\frac{1}{\Phi^2} \sum_{i \geqslant 0} \Vert \Delta_i(F) \Vert_2^2 \; \leqslant \; \Vert F \Vert_2^2 \; \leqslant \; \frac{k}{1 - \Phi_k\Phi} \sum_{i \geqslant 0} \Vert \Delta_i(F) \Vert_2^2
\end{equation}
Similarly, for all $\mu \in \Prob(G)$ with supp$(\mu)\subset B_1(1)$, we have
\[\Vert \mu \ast F \Vert_2^2 \leqslant \frac{k}{1 - \Phi_k\Phi} \sum_{i \geqslant 0} \Vert \mu \ast \Delta_i(F) \Vert_2^2,\]

Further, Lemma \ref{almostorthogonal} implies that for all $i \geqslant 0$, we have 
\[\Vert P_{2^{-2i}}\ast \Delta_i(F)-\Delta_i(F) \Vert_2 \leqslant \frac{4C_0}{2^i}\Vert F \Vert_2 \qquad \text{and} \qquad \Vert P_{2^{-i/2}} \ast \Delta_i(F) \Vert_2 \leqslant \frac{C_0}{2^{i/2}}\Vert F \Vert_2.\]
Therefore,
\[\sum_{i \geqslant 0} 2^{i/2}\Vert P_{2^{-2i}} \ast \Delta_i(F) - \Delta_i(F) \Vert_2^2 \leqslant 2(4C_0)^2 \Vert F \Vert_2^2\]
and 
\[\sum_{i \geqslant 0} 2^{i/2}\Vert P_{2^{-i/2}} \ast \Delta_i(F) \Vert_2^2 \leqslant \frac{\sqrt 2}{\sqrt 2 - 1}C_0^2\Vert F \Vert_2^2.\]
It is now clear that the conclusion of Theorem \ref{L-P} holds for $C>0$ large enough (but still independent of $\mu$ and $F$).
\hfill$\blacksquare$

\subsection{Reduction to functions living at a small scale}
We continue with a consequence of Theorem \ref{L-P} that will allow us to reduce the problem of proving restricted spectral gap to functions that live at an arbitrarily small scale $\delta>0$.

\begin{corollary}\label{level_delta}
Let $C>0$ be the constant provided by Theorem \ref{L-P}.
Let $0 < r < 1$.  Let $B\in\mathcal C(G)$ and $\mu\in\mathcal M(G)$ be a Borel probability measure with supp$(\mu)\subset B_1(1)$. Assume that for any finite dimensional subspace $V\subset L^2(B)$, there is $F\in L^2(B)\ominus V$ such that $\|\mu*F\|_2>r \|F\|_2$. 
Let $\widetilde B\subset G$ be an open set with compact closure which contains the closure of $B$.

Then for every $\delta_0>0$, there exist $F\in L^2(\widetilde B)$ and $0<\delta<\delta_0$ such that
\begin{enumerate}
\item $\|\mu*F\|_2> r \|F\|_2/(2C)$.
\item $\|P_{\delta}*F-F\|_2<\delta^{1/16}\|F\|_2$.
\item $\|P_{\delta^{1/4}}*F\|_2<\delta^{1/16}\|F\|_2$.
\end{enumerate}
\end{corollary}

{\it Proof.}
Let $\delta_0>0$. Choose $N\geqslant 1$ such that $2^{-N}<\delta_0/2$ and $2^{-N/4} < r^2/(16C^3)$.
Since $B$ has compact closure, the operator $L^2(B)\ni F\mapsto P_{\delta}*F\in L^2(G)$ is compact, for any $\delta>0$. 
Hence, the operator $L^2(B)\ni F\mapsto \Delta_i(F)\in L^2(G)$ is compact, for all $i\geqslant 0$.
The hypothesis implies that we can find $F_0 \in L^2(B)$ such that $\|\mu*F_0\|_2>r\|F_0\|_2$ and $\sum_{i=0}^{N-1}\|\Delta_i*F_0\|^2_2 < r^2 \|F_0\|_2^2/(2C)$.

By using Theorem \ref{L-P} (2), we derive that $\sum_{i\geqslant 0}\|\mu*\Delta_i(F_0)\|_2^2\geqslant \|\mu*F_0\|_2^2/C > r^2 \| F_0 \|_2^2/C$. Since $\sum_{i=0}^{N-1}\|\mu*\Delta_i(F_0)\|_2 < r^2 \|F_0\|_2^2/(2C)$, we get that $\sum_{i\geqslant N}\|\mu*\Delta_i(F_0)\|_2^2 > r^2 \|F_0\|_2^2/(2C)$. In combination with Theorem \ref{L-P} (1) we deduce that $\sum_{i\geqslant N}\|\mu*\Delta_i(F_0)\|_2^2 > r^2 (\sum_{i\geqslant N}\|\Delta_i(F_0)\|_2^2)/(2C^2)$, or equivalently

\begin{equation}\label{eq1}\sum_{i\geqslant N}(\|\Delta_i(F_0)\|_2^2-\|\mu*\Delta_i(F_0)\|_2^2)\leqslant \big(1-\frac{r^2}{2C^2}\big)\sum_{i\geqslant N}\|\Delta_i(F_0)\|_2^2. \end{equation}

Since $\sum_{i\geqslant 0}\|\Delta_i(F_0)\|_2^2\geqslant \|F_0\|_2^2/C$ by Theorem \ref{L-P} (2) and $\sum_{i=0}^{N-1}\|\Delta_i*F_0\|^2_2 < \|F_0\|_2^2/(2C)$, we get that $\sum_{i\geqslant 0}\|\Delta_i(F_0)\|_2^2< 2(\sum_{i\geqslant N}\|\Delta_i(F_0)\|_2^2)$. By combining this inequality with Theorem \ref{L-P} (3) and using that $2^{-N/4}<r^2/(16C^3)$ we deduce that 
\begin{align}\label{eq2}
\sum_{i\geqslant N}2^{i/4}\|P_{2^{-2i}}*\Delta_i(F_0)-\Delta_i(F_0)\|_2^2&\leqslant 2^{-N/4}\sum_{i\geqslant 0}2^{i/2}\|P_{2^{-2i}}*\Delta_i(F_0)-\Delta_i(F_0)\|_2^2\\&\leqslant 2^{-N/4}C\Big(\sum_{i\geqslant 0}\|\Delta_i(F_0)\|_2^2\Big)\nonumber\\&< 2^{-N/4+1}C\Big(\sum_{i\geqslant N}\|\Delta_i(F_0)\|_2^2\Big)\leqslant \frac{r^2}{8C^2}\sum_{i\geqslant N}\|\Delta_i(F_0)\|_2^2.\nonumber
\end{align}
 
 Similarly, by using Theorem \ref{L-P} (4), we get that \begin{equation}\label{eq3} \sum_{i\geqslant N}2^{i/4}\|P_{2-i/2}*\Delta_i(F_0)\|_2^2< \frac{r^2}{8C^2}\sum_{i\geqslant N}\|\Delta_i(F_0)\|_2^2.\end{equation}
 
By combining equations \eqref{eq1}, \eqref{eq2}, and \eqref{eq3} we can find $i\geqslant N$ such that \[(\|\Delta_i(F_0)\|_2^2-\|\mu*\Delta_i(F_0)\|_2^2)+2^{i/4}\|P_{2^{-2i}}*\Delta_i(F_0)-\Delta_i(F_0)\|_2^2+2^{i/4}\|P_{2-i/2}*\Delta_i(F_0)\|_2^2<(1-\frac{r^2}{4C^2})\|\Delta_i(F_0)\|_2^2.\]

Let $F:=\Delta_i(F_0)$ and $\delta:=2^{-2i}$. Then $\delta\leqslant 2^{-2N}<\delta_0$. Moreover, the above inequality implies that $\|\mu*F\|_2^2> r^2 \|F\|_2^2/(4C^2)$, $\|P_{\delta}*F-F\|_2^2<\delta^{1/8}\|F\|_2^2$, and $\|P_{\delta^{1/4}}*F \|_2<\delta^{1/8}\|F \|_2^2$. Finally, notice that since $F_0\in L^2(B)$, the support of $F$ is contained in $B_{2^{-i+1}}(1)B \subset B_{\delta_0}(1)B$ and hence in $\widetilde B$, if $\delta_0>0$ is small enough. \hfill $\blacksquare$

\subsection{Proof of Theorem \ref{restricted}} Next, we prove the following ``quantitative restricted spectral gap" theorem for all measures with small support that escape subgroups at a controlled speed. 
It is clear that this result in combination with Theorem \ref{escape} immediately implies Theorem \ref{restricted}.

\begin{theorem}\label{restricted2}
Let $G$ be a connected simple Lie group with trivial center and $B\subset G$ a measurable set with compact closure. Let $d_1,d_2>0$ be given.

Then there exist $c>0$ and $\varepsilon_2>0$ such that the following holds true.
Let $0<\varepsilon<\varepsilon_2$ and $\mu\in\mathcal M(G)$ be a Borel probability on $G$ with $supp(\mu)\subset B_{\varepsilon}(1)$. Assume that for all $\delta > 0$ small enough, we have that for any proper connected closed subgroup $H<G$,
\[\mu^{*2n}(H^{(\delta)})\leqslant\delta^{d_1}, \text{ where } n=\Big\lfloor{d_2\frac{\log{\frac{1}{\delta}}}{\log{\frac{1}{\varepsilon}}}}\Big\rfloor.\]

Then there exists a finite dimensional subspace $V\subset L^2(B)$ such that $\|\mu*F\|_2<\varepsilon^{c}\|F\|_2$, for every $F\in L^2(B)\ominus V$.
\end{theorem}

Theorem \ref{restricted2} also implies the quantitative version of Theorem \ref{restricted} referred to in Remarks \ref{quant} and \ref{spectra}.

{\it Proof.}
Let $B \subset G$ and $d_1,d_2 >0$ be as in the statement of the theorem. 
Let $\widetilde B\subset G$ be an open set with compact closure which contains the closure of $B$. 
Denote $d=\text{dim}(G)$.

We start by quantifying how small $\varepsilon>0$ should be. 
First,
 Theorem \ref{rho} provides constants $a,b,\kappa > 0$ such that for any $F\in L^2(\widetilde B)$ with $\|F\|_2=1$ we have \begin{equation}\label{coeff}\|f*F\|_2^{16d}\leqslant a\|P_{\delta}*F\|_2+b\delta^\kappa,\;\;\text{ for all $f\in L^2(G)$ with $\|f\|_2=1$ and all $0<\delta<1$}. \end{equation} 

Put $q = \min\{\frac{1}{16},\kappa/4\}$ and let $C>0$ be the constant provided in Theorem \ref{L-P}. Choose 
\begin{itemize}
\item $0<\alpha < \frac{q}{16d}$ and denote by $c_0$ and $\varepsilon_0$ the corresponding constants given by Corollary \ref{flattening}.
\item $c>0$ such that $2c_0c<\min\{\frac{1}{16},\frac{q}{16d}-\alpha\}$.
\item $0<\varepsilon < \varepsilon_0$ small enough so that $2c_0(c+\frac{\log{2C}}{\log{\frac{1}{\varepsilon}}})<\min\{\frac{1}{16},\frac{q}{16d}-\alpha\}$.
\end{itemize}

Next, take a probability measure $\mu$ on $G$ supported on $B_{\varepsilon}(1)$ such that for all $\delta > 0$ small enough, we have $\mu^{*2n}(H^{(\delta)})\leqslant\delta^{d_1},$ where $n=\lfloor{d_2\frac{\log{\frac{1}{\delta}}}{\log{\frac{1}{\varepsilon}}}}\rfloor$, for any proper connected closed subgroup $H<G$.
By Corollary \ref{flattening}, there exists $\delta_0 > 0$ such that for all $0<\delta<\delta_0$, we have that
\begin{equation}\label{flattt} \|\mu^{*n}*P_{\delta}\|_2\leqslant\delta^{-\alpha},\;\;\text{ for all } n \geqslant \Big\lfloor{c_0\frac{\log{\frac{1}{\delta}}}{\log\frac{1}{\varepsilon}}}\Big\rfloor.\end{equation}
Taking $\delta_0$ smaller if necessary, we can assume that $\delta^{2c_0(c+\frac{\log{2C}}{\log{\frac{1}{\varepsilon}}})}>(a+b)^{\frac{1}{16d}}\delta^{\frac{q}{16d}-\alpha}+\delta^{\frac{1}{16}}$, for $\delta<\delta_0$.

Now, assume by contradiction that the measure $\mu$ does not satisfy the conclusion of the theorem. Then by Corollary \ref{level_delta}, there exists $F\in L^2(\widetilde{B})$ with $\|F\|_2=1$ and $0 < \delta < \delta_0$ such that
\begin{enumerate}
\item $\|\mu*F\|_2 > \varepsilon^c/(2C)$.
\item $\|P_{\delta}*F-F\|_2<\delta^{1/16}$.
\item $\|P_{\delta^{\frac{1}{4}}}*F\|_2<\delta^{1/16}$.
\end{enumerate}

Let $n=\Big\lfloor{c_0\frac{\log{\frac{1}{\delta}}}{\log\frac{1}{\varepsilon}}}\Big\rfloor$.
Since $\mu$ is symmetric, by using Lemma \ref{powers} we derive that 
\begin{align*}
\Big(\frac{\varepsilon^c}{2C}\Big)^{2n} \leqslant \|\mu*F\|^{2n} \leqslant \|\mu^{*n}*F\|_2 & \leqslant \|\mu^{*n} \ast P_{\delta} \ast F\|_2 + \|P_{\delta} \ast F-F\|_2
\end{align*}

On the other hand, by combining \eqref{coeff} and \eqref{flattt} we get that \begin{align*}\|\mu^{*n}*P_{\delta}*F\|_2&\leqslant \big(a\|P_{\delta^{1/4}}*F\|_2+\delta^{\kappa/4}\big)^{\frac{1}{16d}}\|\mu^{*n}*P_{\delta}\|_2\\&\leqslant(a\delta^{\frac{1}{16}}+b\delta^{\kappa/4})^{\frac{1}{16d}}\;\delta^{-\alpha}\leqslant(a+b)^{\frac{1}{16d}}\delta^{\frac{q}{16d}-\alpha}. \end{align*}

By putting the last two inequalities together we get that $\big(\frac{\varepsilon^c}{2C}\big)^{2n} \leqslant(a+b)^{\frac{1}{16d}}\delta^{\frac{q}{16d}-\alpha}+\delta^{\frac{1}{16}}.$
Since $\big(\frac{\varepsilon^c}{2C}\big)^{2n}\geqslant  \big(\frac{\varepsilon^c}{2C}\big)^{2c_0\frac{\log{\frac{1}{\delta}}}{\log{\frac{1}{\varepsilon}}}}=\delta^{2c_0(c+\frac{\log{2C}}{\log{\frac{1}{\varepsilon}}})},$ this contradicts the choice of $\delta>0$.
 \hfill$\blacksquare$

\subsection{Proof of Corollary \ref{by}}

Let $\Gamma$, $G$, $H$ and $B \subset G/H$ be as in the statement of Corollary \ref{by}. Recall that the measure $m_{G/H}$ on $G/H$ arises from a rho-function for the pair $(G,H)$ (see  \cite[Theorem B.1.4.]{BdHV08}).
 Thus, there exists a continuous function $\rho: G \to \R_+^*$ such that
\begin{equation}\label{intrho}
\int_G f(x)\rho(x) \d x = \int_{G/H}\int_H f(xh)\d h\, \d m_{G/H}(xH), \; \text{ for all } f \in C_c(G).
\end{equation}
Of course, equality \ref{intrho} holds more generally for any function $f \in L^1(G)$ with compact support. The measure $m_{G/H}$ is not necessarily $G$-invariant, but the function $\rho$ allows to determine the translates of $g \cdot m_{G/H}$ (see \cite[Lemma B.1.3]{BdHV08}):
\begin{equation}\label{RNrho}
\frac{d(g \cdot m_{G/H})}{dm_{G/H}}(xH) = \frac{\rho(gx)}{\rho(x)}, \; \text{ for all } x,g \in G.
\end{equation}

Put $B_1 = B_1(1) \cdot B \subset G/H$ and $B_2 = B_1(1)^{-1} \cdot B_1 \subset G/H$. Let $p:G\to G/H$ be the canonical projection.
Let $\tilde B_1, \tilde B_2 \subset G$ be open sets with compact closures such that $B_i \subset p(\tilde B_i)$ for $i = 1,2$. Replacing $\tilde B_i$ by $\tilde B_i \cdot K$ for some compact set $K \subset H$ with non empty interior, we may also assume that $\int_H 1_{\tilde B_i}(xh)\d h$ is bounded away from $0$ uniformly in $x \in \tilde B_i$.
Then using \eqref{intrho} there exists $\beta > 0$ such that $\Vert F \Vert_2/\beta \leqslant \Vert F \circ p \Vert_{2,\tilde B_i} \leqslant \beta \Vert F \Vert_2$ for all $F \in L^2(B_i)$, for both $i = 1$ and $i = 2$.

Fix $\varepsilon \in (0,1)$ small enough so that $\vert \sqrt{\frac{\rho(x)}{\rho(gx)}} - 1\vert \leqslant\frac{1}{4}$, for all $g \in B_\varepsilon(1)$ and $x \in \tilde B_1$. Take $r > 0$ such that $2r\beta^4 < 1/16$. 

By Theorem \ref{restricted2} there exist a finite dimensional space $V \subset L^2(\tilde B_2)$ and a finite set $T \subset\Gamma$ such that the measure $\mu = \frac{1}{2\vert T \vert} \sum_{g \in T} (\delta_g + \delta_{g^{-1}})$ satisfies supp$(\mu)\subset B_{\varepsilon}(1)$ and $ \Vert \mu \ast F\Vert_2 < r \Vert F \Vert_2$, for all $F \in L^2(\tilde B_2) \ominus V$.

Take a sequence of functions $F_n \in L^2(B)$ which converges weakly to $0$ and such that $\Vert F_n \Vert_2 = 1$ for all $n$. To prove the corollary, it is enough to show that eventually $\Vert \pi(\mu)(F_n) \Vert_2 < \frac{1}{2}$.

First, remark that by our choice of $\varepsilon$, we have $\vert 1 - \sqrt{\frac{\rho(gx)}{\rho(x)}}\vert \leqslant \frac{1}{4}\sqrt{\frac{\rho(gx)}{\rho(x)}}$ for all $g \in \supp(\mu)$ and $x \in\tilde B_1$. Thus, Equation \eqref{RNrho} gives for all $F \in L^2(B)$:
\begin{align*}
\Vert \pi(\mu)(F) -  \mu \ast F\Vert_2 & = \frac{1}{2\vert T \vert } \Vert \sum_{g \in T \cup T^{-1}} (\sqrt{\frac{\rho(g \, \cdot)}{\rho}} - 1)F(g^{-1} \, \cdot)\Vert_2\\
& \leqslant \frac{1}{4} \frac{1}{2\vert T \vert } \sum_{g \in T \cup T^{-1}} \Vert \sqrt{\frac{\rho(g \, \cdot )}{\rho}}F(g^{-1} \, \cdot )\Vert_2 = \frac{1}{4} \Vert F \Vert_2.
\end{align*}
Therefore, for all $n$ we have
\begin{equation}\label{boundRN}
\Vert \pi(\mu)(F_n)\Vert_2 \leqslant \Vert \mu \ast F_n\Vert_2 + \frac{1}{4}.
\end{equation}
So we are left to bound $\Vert \mu \ast F_n\Vert_2$ by $\frac{1}{4}$ for all $n$ large enough. 
Since $\mu \ast F_n$ is supported on $B_1$, by the definition of $\beta$, we have that
\begin{align*}
\Vert \mu \ast F_n \Vert_2^2 \leqslant \beta^2 \Vert (\mu \ast F_n) \circ p \Vert_{2,\tilde B_1}^2 & = \beta^2 \Vert \mu \ast (F_n \circ p) \Vert_{2,\tilde B_1}^2\\
& = \beta^2 \langle \mu \ast (F_n \circ p) 1_{\tilde B_1}, \mu \ast (1_{\tilde B_2}.(F_n \circ p)) \rangle\\
& \leqslant \beta^2 \Vert \mu \ast (F_n \circ p)\Vert_{2,\tilde B_1} \Vert \mu \ast (1_{\tilde B_2}.(F_n \circ p)) \Vert_2
\end{align*}
The second line above comes from the fact that $1_{\tilde B_1} \leqslant 1_{g \cdot \tilde B_2}$ for all $g \in \supp(\mu) \subset B_1(1)$. Using the same fact we moreover see that $\Vert \mu \ast (F_n \circ p)\Vert_{2,\tilde B_1} \leqslant \Vert F_n \circ p \Vert_{2,\tilde B_2} \leqslant \beta\Vert F_n \Vert_2 = \beta$. In summary,
\[\Vert \mu \ast F_n \Vert_2^2 \leqslant \beta^3  \Vert \mu \ast (1_{\tilde B_2}.(F_n \circ p)) \Vert_2.\]
Using \eqref{intrho} one easily checks that the sequence $(1_{\tilde B_2}.(F_n \circ p))_n \subset L^2(\tilde B_2)$ goes weakly to $0$. Hence, we deduce from the restricted spectral gap assumption on $\mu$ that for $n$ large enough,
\[\Vert \mu \ast (1_{\tilde B_2}.(F_n \circ p)) \Vert_2 < 2r \Vert 1_{\tilde B_2}.(F_n \circ p) \Vert_2 = 2r\Vert F_n \circ p \Vert_{2,\tilde B_2} \leqslant 2r\beta.\]
Altogether, we get that for $n$ large enough:
\[\Vert \mu \ast F_n \Vert_2^2 < 2r\beta^4 \leqslant \frac{1}{16}.\]
Combining this with \eqref{boundRN}, we indeed get that $\Vert \pi(\mu)(F_n) \Vert_2 < \frac{1}{2}$ for $n$ large enough.\hfill $\blacksquare$

 \section{The Banach-Ruziewicz problem}\label{6}
 This section is devoted to the proof of Theorem \ref{BR}. 
  Moreover, we will show the following:

\begin{theorem}\label{BRtext}

Let $G$ be a l.c.s.c. group and  $\Gamma<G$ a countable dense subgroup.  
 
 Then the following four conditions are equivalent:

\begin{enumerate}
\item If $\nu:\mathcal C(G)\rightarrow [0,\infty)$ is a $\Gamma$-invariant, finitely additive measure, then there exists $\alpha\geqslant 0$ such that $\nu(A)=\alpha |A|$, for all $A\in\mathcal C(G)$.
\item If $\Phi:L^{\infty}_{\text{c}}(G,m_G)\rightarrow\mathbb C$ is a $\Gamma$-invariant, positive linear functional, then there exists $\alpha\geqslant 0$ such that  $\Phi(f)=\alpha{\int_{G}f\;\text{d}m_G}$, for all $f\in L^{\infty}_{\text{c}}(G,m_G)$.
\item The translation action $\Gamma\curvearrowright (G,m_G)$ has local spectral gap with respect to a measurable set $B\subset G$ with compact closure and non-empty interior.
\item The translation action $\Gamma\curvearrowright (G,m_G)$ is strongly ergodic.
\end{enumerate}

\end{theorem}

\begin{remark}\label{cominv}
Suppose that $G$ is {\it compact}. It is clear that  (1) $\Longrightarrow$ (2). Further, it is well-known that  (2) $\Longleftrightarrow$ (3) and  (3) $\Longrightarrow$ (4) (see Theorem \ref{spec}).  The implication (4) $\Longleftrightarrow$ (3) was established recently in \cite[Theorem 4]{AE10}.
 Let us also explain why (2) $\Longrightarrow$ (1). Note that if (2) holds, then (3) does as well, hence $\Gamma$ is non-amenable.
 Since the action $\Gamma\curvearrowright G$ is free, it follows that $G$ is $\Gamma$-paradoxical (see Definition \ref{echidec}). The proof of \cite[Theorem 2.1.17]{Lu94} implies that any subsets $B,C\subset G$ with non-empty interior are equidecomposable.  Further, the proof of \cite[Proposition 2.1.12]{Lu94} gives that any finitely additive $\Gamma$-invariant measure $\nu:\mathcal C(G)\rightarrow [0,\infty)$ is absolutely continuous with respect to $m_G$. It follows readily that (1) holds.
 Theorem \ref{BRtext} is therefore contained in the literature when $G$ is compact. Our contribution is to show that it holds for {\it arbitrary} locally compact groups. 
\end{remark}

Turning to locally compact groups $G$, the non-trivial implications, which we will address below, are  (2) $\Longrightarrow$ (3), (2) $\Longrightarrow$ (1), and (4) $\Longrightarrow$ (3).

 \subsection{Local spectral gap and uniqueness of invariant means} In order to prove implication (2) $\Longrightarrow$ (3) from Theorem \ref{BRtext},
 we give an equivalent formulation of local spectral gap in terms of uniqueness of invariant linear functionals (see Theorem \ref{mean}). 
 
This generalizes a well-known result for probability measure preserving actions. 
Let $\Gamma\curvearrowright (X,\mu)$ be a probability measure preserving action of a countable group $\Gamma$.
Then integration against $\mu$ defines a $\Gamma$-invariant {\bf mean} (i.e. a unital positive linear functional) on $L^{\infty}(X,\mu)$.
In the early 1980's, it was realized that whether this is the unique $\Gamma$-invariant mean on $L^{\infty}(X,\mu)$ is equivalent to the spectral gap of the action. More precisely,  the following was shown:

\begin{theorem}\emph{\cite{Ro81, Sc81}}\label{spec} Let $\Gamma\curvearrowright (X,\mu)$ be an ergodic measure preserving action of a countable group $\Gamma$ on a probability space $(X,\mu)$. 
Consider the following conditions:

\begin{enumerate}
\item If $\Phi:L^{\infty}(X,\mu)\rightarrow\mathbb C$ is a $\Gamma$-invariant mean, then ${\Phi(f)=\int_Xf\;\text{d}\mu}$, for all $f\in L^{\infty}(X,\mu)$.
\item There does not exist a sequence $\{A_n\}$ of measurable subsets of $X$ such that $\mu(A_n)>0$, for all $n$, $\lim\limits_{n\rightarrow\infty}\mu(A_n)=0$, and $\lim\limits_{n\rightarrow\infty}\mu(gA_n\Delta A_n)/\mu(A_n)=0$, for all $g\in\Gamma$. 
\item If a sequence $\varphi_n\in L^1(X,\mu)$  of positive functions satisfies ${\int_{X}\varphi_n\;\text{d}\mu=1}$, for all $n$, and
 $\lim\limits_{n}\|g\cdot\varphi_n-\varphi_n\|_1=0$, for all $g\in\Gamma$, then $\lim\limits_{n}\|\varphi_n-1\|_{1}=0$.
\item The action $\Gamma\curvearrowright (X,\mu)$ has spectral gap.
\item The action $\Gamma\curvearrowright (X,\mu)$ is strongly ergodic.
\end{enumerate}

Then conditions (1)-(4) are equivalent and they all imply condition (5).
\end{theorem}

The equivalence of   (1) and (2) is due to  Rosenblatt  \cite[Theorem 1.4]{Ro81}. The equivalence of (1), (3) and (4), and the fact that (1) implies (5) are due to Schmidt \cite[Propositions 2.2 and 2.3]{Sc81} (see also \cite[Section 5]{FS99}, where a gap from \cite{Sc81} is fixed).

The main goal of this section is to generalize Theorem \ref{spec} to arbitrary measure preserving actions.  In order to state our result, we need to introduce some notation:

\begin{notation}\label{F_B} Let  $\Gamma\curvearrowright (X,\mu)$ be  a measure preserving action of a countable group $\Gamma$ on a standard measure space $(X,\mu)$. Let $B\subset X$ be a measurable set.
 We denote by
 \begin{itemize}
  \item $\mathcal C_B(X)$ the family of measurable subsets $C\subset X$ for which we can find $g_1,...,g_n\in\Gamma$ such that 
  $C\subset\cup_{i=1}^ng_iB$, almost everywhere, and by 
  
  \item $L^{\infty}_{B}(X,\mu)$  the set of functions $f\in L^{\infty}(X,\mu)$ whose essential support belongs to $\mathcal C_B(X)$.
  \end{itemize}
\end{notation} 


\begin{remark}\label{CB}
If $G$ is a l.c.s.c. group and $B\in\mathcal C(G)$ is a set with non-empty interior, then $\mathcal C_B(G)=\mathcal C(G)$ and $L^{\infty}_B(G)=L^{\infty}_{\text{c}}(G)$. 
\end{remark}

\begin{theorem}\label{mean}
Let $\Gamma\curvearrowright (X,\mu)$ be an ergodic measure preserving action of a countable group $\Gamma$ on a standard measure space $(X,\mu)$. 
Let $B\subset X$ be a measurable set with $0<\mu(B)<\infty$.
Consider the following conditions:

\begin{enumerate}
\item If $\Phi:L^{\infty}_{B}(X,\mu)\rightarrow\mathbb C$ is a $\Gamma$-invariant positive linear functional, then there exists $\alpha\geqslant 0$ such that $\Phi(f)=\alpha{\int_{X}f\;\text{d}\mu}$, for all $f\in L^{\infty}_{B}(X,\mu)$.
\vskip 0.03in
\item If a sequence $\{A_n\}$ of measurable subsets of $X$ satisfies that $\mu(A_n\cap B)>0$, for all $n$, and $\lim\limits_{n}\mu((gA_n\Delta A_n)\cap B)/\mu(A_n\cap B)=0$, for all $g\in\Gamma$, then  $\lim\limits_{n}\mu(A_n\cap B)=\mu(B)$.
\vskip 0.03in
\item If a sequence $\varphi_n\in L^1(X,\mu)$  of positive functions satisfies ${\int_{B}\varphi_n\;\text{d}\mu=\mu(B)}$, for all $n$, and
 $\lim\limits_{n}\|g\cdot\varphi_n-\varphi_n\|_{1,B}=0$, for all $g\in\Gamma$, then $\lim\limits_{n}\|\varphi_n-1\|_{1,B}=0$.
 \item The action $\Gamma\curvearrowright (X,\mu)$ has local spectral gap with respect to $B$.
 
  \item The action $\Gamma\curvearrowright (X,\mu)$ is strongly ergodic.
\end{enumerate}

Then conditions (1)-(4) are equivalent and they all imply condition (5).

\end{theorem}

Theorem \ref{mean} is motivated in part by Margulis' proof of  \cite[Theorem 3]{Ma82}.
Note that in the case $\mu(X)=1$ and $B=X$, it recovers Theorem \ref{spec}.

{\bf Proof}.
Let  $\lim\limits_n$  be a bounded linear functional on $\ell^{\infty}(\mathbb N)$ which extends the usual limit.

{\bf (1) $\Longrightarrow$ (2)}.  Assume  that (1) holds and (2) is false. Let $S$ be the set of $L\geqslant 0$
for which there exists a sequence $\{A_n\}$ of measurable subsets of $X$ such that $\mu(A_n\cap B)>0$, for all $n$,  $\lim\limits_{n}\mu((gA_n\Delta A_n)\cap B)/\mu(A_n\cap B)=0$, for all $g\in\Gamma$, and $\lim\limits_{n}\mu(A_n\cap B)=L$. It is easy to see that $S\subset [0,\infty)$ is a non-empty closed set. We denote by $\ell$
the minimum of $S$ and by $\{A_n\}$ the corresponding sequence of measurable subsets of $X$.  Since (2) is false,  $\ell<\mu(B)$.

 We define a positive linear functional $\Phi:L^{\infty}_{B}(X,\mu)\rightarrow\mathbb C$  by letting $$\Phi(f)=\lim\limits_{n}\frac{1}{\mu(A_n\cap B)}\int_{A_n}f\;\text{d}\mu,\;\;\;\text{for all $f\in L^{\infty}_{B}(X,\mu)$.}$$
 
We claim that $\Phi$ is well-defined and $\Gamma$-invariant. To this end, let $f\in L^{\infty}_{B}(X,\mu)$ and $g\in\Gamma$. Denote by $A$ the support of $f$.
Since $A\in\mathcal C_B(X)$, we can find $g_1,...,g_k\in\Gamma$ with $A\subset\cup_{i=1}^kg_iB$. Then $\mu(A_n\cap A)\leqslant\sum_{i=1}^k\mu(g_i^{-1}A_n\cap B)$ and we get that $\limsup\limits_{n\rightarrow\infty}\mu(A_n\cap A)/\mu(A_n\cap B)\leqslant k$. Since ${|\frac{1}{\mu(A_n\cap B)}\int_{A_n}f\;\text{d}\mu|\leqslant \|f\|_{\infty}\frac{\mu(A_n\cap A)}{\mu(A_n\cap B)}}$, it follows that $\Phi(f)$ is well-defined.
Further, we have \begin{align*}\label{fi2}|\int_{A_n}g\cdot f\;\text{d}\mu-\int_{A_n}f\;\text{d}\mu|\leqslant\int_{g^{-1}A_n\Delta A_n}|f|\;\text{d}\mu&=\int_{(g^{-1}A_n\Delta A_n)\cap A}|f|\;\text{d}\mu\\&\leqslant\|f\|_{\infty}\;\mu((g^{-1}A_n\Delta A_n)\cap A).\end{align*}

Since $A\in\mathcal C_B(X)$, we have  $\lim\limits_{n}\mu((gA_n\Delta A_n)\cap A)/\mu(A_n\cap B)=0$. In combination with the above, this implies that $\Phi(g\cdot f)=\Phi(f)$. Therefore, $\Phi$ is $\Gamma$-invariant.
Since $\Phi(1_B)=1$ and condition (1) is assumed true, we get that $\Phi(f)={\frac{1}{\mu(B)}\int_Xf\;\text{d}\mu}$, for all $f\in L^{\infty}_{B}(X,\mu)$. 

We are now ready to derive a contradiction. Firstly, assume that $\ell=\lim\limits_{n}\mu(A_n\cap B)=0$. In this case, after passing to a subsequence, we may assume  that $0<\mu(A_n\cap B)<3^{-n}\mu(B)$, for all $n\geqslant 1$. 
Let $C=B\setminus(\cup_{n\geqslant 1}A_n)$. It follows that $\mu(C)>0$ and hence $\Phi(1_C)=\mu(C)/\mu(B)>0$. On the other hand, it is clear from the definition of $\Phi$ that $\Phi(1_C)=0$. This gives a contradiction.

Secondly, assume that $0<\ell<\mu(B)$. If $C\in\mathcal C_B(X)$, then $\lim\limits_n\mu(A_n\cap C)=\ell\Phi(1_C)=(\ell\mu(C))/\mu(B)$. Denoting $A_{m,n}=A_m\cap A_n$, we further get that 

\begin{equation}\label{muu}\lim\limits_{m}\lim\limits_n \mu(A_{m,n}\cap C)=\lim\limits_{m}\frac{\ell\;\mu(A_m\cap C)}{\mu(B)}=\frac{\ell^2\mu(C)}{\mu(B)^2}. \end{equation}

In particular,
 ${\lim\limits_{m}\lim\limits_{n}\mu(A_{m,n}\cap B)=\ell^2/\mu(B)}$.
Let $g\in\Gamma$. Since $\ell>0$, the assumptions on $\{A_n\}$  imply that $\lim\limits_n\mu((gA_n\Delta A_n)\cap C)=0$. Since $g_{A_{m,n}}\Delta A_{m,n}\subset (gA_m\Delta A_m)\cup (gA_n\Delta A_n)$, we get that $\lim\limits_m\lim\limits_n\mu((gA_{m,n}\cap A_{m,n})\cap C)=0$. 
 It follows that there is a sequence of the form  $\tilde A_k=A_{m(k),n(k)}$ such that $\lim\limits_{k}\mu(\tilde A_k\cap B)=\ell^2/\mu(B)$ and  $\lim\limits_{k}\mu((g\tilde A_k\Delta\tilde A_k)\cap C)=0$, for all $C\in\mathcal C_B(X)$ and $g\in\Gamma$.
This  implies that $\ell^2/\mu(B)\in S$. Since $\ell^2/\mu(B)<\ell$, this contradicts the minimality of $\ell$. \hfill$\square$

{\bf (2) $\Longrightarrow$ (3)}. The proof relies on a variation of Namioka's trick. Suppose (2) is true. By contradiction, assume that there is a sequence of positive functions $\varphi_n\in L^1(X,\mu)$  satisfying ${\int_{B}\varphi_n\;\text{d}\mu=\mu(B)}$, for all $n$, and
 $\lim\limits_{n}\|g\cdot\varphi_n-\varphi_n\|_{1,B}=0$, for all $g\in\Gamma$, such that $\|\varphi_n-1\|_{1,B}\not\rightarrow 0$. After passing to a subsequence, assume that there is $\delta>0$ such that $\|\varphi_n-1\|_{1,B}\geqslant\delta$, for all $n\geqslant 1$.
 
 Let $c\in (0, \frac{\delta}{2\mu(B)})$ and denote $\delta_0:=\frac{\delta-2c\mu(B)}{2}>0$.
 Fix $n\geqslant 1$. We define $\psi_n:X\rightarrow\mathbb R$ by letting $\psi_n(x)=\varphi_n(x)-1$. Since ${\int_{B}\psi_n\;\text{d}\mu=0}$ and $c 1_{\{\psi_n\geqslant c\}}+\psi_n 1_{\{0\leqslant\psi_n<c\}}\leqslant c1_{\{\psi_n\geqslant 0\}}\leqslant c\varphi_n$, we have  
 
\begin{align*}\delta\leqslant \int_{B}|\psi_n|\;\text{d}\mu & =  2\int_{B}\psi_n1_{\{\psi_n\geqslant 0\}}\;\text{d}\mu=2\int_{B}\psi_n1_{\{\psi_n\geqslant c\}}\;\text{d}\mu+2\int_{B}\psi_n1_{\{0\leqslant\psi_n<c\}}\;\text{d}\mu\\ & =2\int_{B}(\psi_n-c)1_{\{\psi_n\geqslant c\}}\;\text{d}\mu+2\int_{B}(c 1_{\{\psi_n\geqslant c\}}+\psi_n 1_{\{0\leqslant\psi_n<c\}})\;\text{d}\mu\\ & \leqslant  2\int_{B}(\psi_n-c)1_{\{\psi_n\geqslant c\}}\;\text{d}\mu+2c\int_{B}\varphi_n\;\text{d}\mu\\ & = 2\int_{B}(\varphi_n-(1+c))1_{\{\varphi_n\geqslant (1+c)\}}\;\text{d}\mu+2c\mu(B).\end{align*}
 
It follows that ${\int_{B}(\varphi_n-(1+c))1_{\{\varphi_n\geqslant (1+c)\}}\;\text{d}\mu\geqslant\delta_0.}$
For  $t\geqslant 0$, we put $A_{t,n}=\{x\in X|\varphi_n(x)\geqslant t\}$. By combining this inequality with Fubini's theorem we get that \begin{equation}\label{ecu1}\int_{1+c}^{\infty}\mu(A_{t,n}\cap B)\;\text{d}t=\int_{B}\int_{1+c}^{\infty}1_{\{\varphi_n\geqslant t\}}\;\text{d}t\;\text{d}\mu=\int_{B}(\varphi_n-(1+c))1_{\{\varphi_n\geqslant (1+c)\}}\;\text{d}\mu\geqslant\delta_0. \end{equation}

Next, let $g\in\Gamma$. By using a similar calculation to the above we get that 

\begin{align*}\label{ecu2}\int_0^{\infty}\mu((gA_{t,n}\Delta A_{t,n})\cap B)\;\text{d}t & =\int_{B}\int_0^{\infty}|1_{gA_{t,n}}(x)-1_{A_{t,n}}(x)|\;\text{d}t\;\text{d}\mu(x) \\ & = \int_{B}\int_0^{\infty}|1_{\{\varphi_n(g^{-1}x)\geqslant t\}}-1_{\{\varphi_n(x)\geqslant t\}}|\;\text{d}t\;\text{d}\mu(x) \\ & = \int_{B}|\varphi_n(g^{-1}x)-\varphi_n(x)|\;\text{d}\mu(x)=\|g\cdot\varphi_n-\varphi_n\|_{1,B}.\end{align*}

Now, fix a finite set $F\subset\Gamma$ and $\varepsilon>0$.
We claim that  there is a  measurable set $A\subset X$  satisfying $\mu(A\cap B)\in (0,\frac{\mu(B)}{1+c}]$, and $\mu((gA\Delta A)\cap B)/\mu(A\cap B)<\varepsilon$, for all $g\in F$. To this end,  note that since $\|g\cdot\varphi_n-\varphi_n\|_{1,B}\rightarrow 0$, for all $g\in\Gamma$, we can find $n\geqslant 1$ such that $\sum_{g\in F}\|g\cdot\varphi_n-\varphi_n\|_{1,B}<\varepsilon\delta_0$. By combining \ref{ecu1} with the last displayed identity it follows that $$\varepsilon\int_{1+c}^{\infty}\mu(A_{t,n}\cap B)\;\text{d}t\geqslant\varepsilon\delta_0>\int_0^{\infty}(\sum_{g\in F}\mu((gA_{t,n}\Delta A_{t,n})\cap B))\;\text{d}t.$$

From this we deduce that there is $t\geqslant 1+c$ such that $\sum_{g\in F}\mu((gA_{t,n}\Delta A_{t,n})\cap B)<\varepsilon\mu(A_{t,n}\cap B)$. Also,  ${t\mu(A_{t,n})\leqslant\int_{B}\varphi_n\;\text{d}\mu=\mu(B)}$, hence $\mu(A_{t,n})\leqslant\frac{\mu(B)}{t}\leqslant\frac{\mu(B)}{1+c}$. 
Thus, $A=A_{t,n}$ verifies the claim.

By using the claim we can construct a sequence $\{A_n\}$ of measurable subsets of $X$ which satisfy that  $\mu(A_n\cap B)\in (0,\frac{\mu(B)}{1+c}]$,  for all $n$, and $\lim\limits_{n}\mu((gA_n\Delta A_n)\cap B)/\mu(A_n\cap B)=0$, for all $g\in\Gamma$.  Since  (2) is assumed true, we would get that $\lim\limits_{n}\mu(A_n\cap B)=\mu(B)$, which is a contradiction. \hfill$\square$

{\bf  (3) $\Longrightarrow$ (1)}. The proof follows the proof of \cite[Theorem 3]{Ma82} and relies on Day's convexity trick.
Suppose that (3) holds and let $\Phi:L^{\infty}_{B}(X,\mu)\rightarrow\mathbb C$ be a $\Gamma$-invariant, positive linear functional. If  $\Phi(1_B)=0$, then we clearly have that $\Phi\equiv 0$. Thus, we may assume that $\Phi(1_B)>0$. 
After replacing $\Phi$ with $\frac{\mu(B)}{\Phi(1_B)}\Phi$, we may further assume that $\Phi(1_B)=\mu(B)$.

Let $f_0\in L^{\infty}_{B}(X,\mu)$. We will show that $\Phi(f_0)={\int_X f_0\;\text{d}\mu}$.
To this end, we denote by $C$ the support of $f_0$ and fix $\varepsilon>0$. Let $F\subset\Gamma$ be a finite set. Define $Y=C\cup B\cup(\cup_{g\in F}g^{-1}B)$. Since $\mu(Y)<\infty$, we can find a net of positive functions $\varphi_i\in L^1(Y)$ such that \begin{equation}\label{equ1}\lim\limits_{i}\int_{Y}\varphi_if\;\text{d}\mu=\Phi(f),\;\;\text{for all $f\in L^{\infty}(Y)$}.  \end{equation}

Then for all $g\in F$ and every $f\in L^{\infty}(B)$  we have that $f,g^{-1}\cdot f\in L^{\infty}(Y)$ and hence \begin{equation}\label{equ2}\lim\limits_{i}\int_B (g\cdot\varphi_i-\varphi_i)f\;\text{d}\mu=\lim\limits_{i}\int_Y\varphi_i(g^{-1}\cdot f-f)\;\text{d}\mu=\Phi(g^{-1}\cdot f-f)=0. \end{equation}

Denote by $\mathcal A$ the set of positive functions $\varphi\in L^1(Y)$ such that ${|\int_{Y}\varphi f_0\;\text{d}\mu-\Phi(f_0)|<\varepsilon}$ and ${\int_{B}\varphi\;\text{d}\mu=\Phi(1_B)}$.
Then \ref{equ1} and \ref{equ2}  imply that $0=(0)_{g\in F}$ belongs to the weak closure of $$\mathcal A_{F}:=\{((g\cdot\varphi-\varphi)_{|B})_{g\in F}|\varphi\in\mathcal A\}\subset L^1(B)^{|F|}.$$

Since $\mathcal A_{F}$ is a convex set, its weak and norm closures concide.
Thus, $0$ belongs to the norm closure of $\mathcal A_{F}$, for every finite set $F\subset\Gamma$.   It follows that there is a sequence $\varphi_n\in\mathcal A$ such that $\lim\limits_{n\rightarrow\infty}\|g\cdot\varphi_n-\varphi_n\|_{1,B}=0$, for all $g\in\Gamma$.
Since ${\int_{B}\varphi_n\;\text{d}\mu=\Phi(1_B)=\mu(B)}$, for all $n$,  condition (3) gives that $\lim\limits_{n\rightarrow\infty}\|\varphi_n-1\|_{1,B}=0$. Further, we get that $\lim\limits_{n\rightarrow\infty}\|\varphi_n-1\|_{1,gB}=0$, for all $g\in\Gamma$. 

Since $Y\in\mathcal C_B(X)$, we deduce that $\lim\limits_{n\rightarrow\infty}\|\varphi_n-1\|_{1,Y}=0$. Since $|{\int_{Y}\varphi_n f_0\;\text{d}\mu-\Phi(f_0)|<\varepsilon}$, for all $n$, we get that ${|\int_{Y}f_0\;\text{d}\mu-\Phi(f_0)|\leqslant\varepsilon}$. As $\varepsilon>0$ is arbitrary, we conclude that $\Phi(f_0)={\int_{X}f_0\;\text{d}\mu}$. \hfill$\square$

{\bf  (3) $\Longrightarrow$ (4)}. Assume that (3) holds and that (4) is false. Then we can find a sequence $\xi_n\in L^2(X)$ such that ${\int_B\xi_n\;\text{d}\mu=0}$ and $\|\xi_n\|_{2,B}=\sqrt{\mu(B)}$, for all $n$, and $\lim\limits_{n\rightarrow\infty}\|g\cdot\xi_n-\xi_n\|_{2,B}=0$, for all $g\in\Gamma$. Moreover,  we may assume that $\xi_n$ is real-valued, for all $n$.

Define $\varphi_n:=|\xi_n|^2\in L^1(X)$. 
Then  ${\int_B\varphi_n\;\text{d}\mu=\mu(B)}$ and $\lim\limits_{n\rightarrow\infty}\|g\cdot\varphi_n-\varphi_n\|_{1,B}=0$, for all $g\in\Gamma$.
 By using condition (3) we deduce that $\lim\limits_{n\rightarrow\infty}\|\varphi_n-1\|_{1,B}=\lim\limits_{n\rightarrow\infty}\||\xi_n|^2-1\|_{1,B}=0$. 
  Since $\xi_n$ is real-valued, it follows that there exists a sequence $\{A_n\}$ of measurable subsets of $X$ such that $\lim\limits_{n\rightarrow\infty}\|\xi_n-(1_{A_n}-1_{X\setminus A_n})\|_{2,B}=0$. By using the defining properties of  $\{\xi_n\}$ we get that $\lim\limits_{n\rightarrow\infty}\mu((gA_n\Delta A_n)\cap B)=0$ and $\lim\limits_{n\rightarrow\infty}\mu(A_n\cap B)=\mu(B)/2$. This contradicts condition (2), which, as shown above, is equivalent to (3). \hfill$\square$

 {\bf (4) $\Longrightarrow$ (3)}. The proof of this implication is easy and so we leave it to the reader. \hfill$\square$

 {\bf (2) $\Longrightarrow$ (5)}. 
 To prove this final implication, we need a lemma whose proof we leave to the reader.

\begin{lemma}\label{invar}
Let $(X,\mu)$ be a standard measure space. Let $\nu$ be a Borel probability measure on $X$ which is equivalent to $\mu$.
Let $\{A_n\}$ be a sequence of measurable subsets of $X$. 
Then $\lim\limits_{n\rightarrow\infty}\nu(A_n)=0$ if and only if $\lim\limits_{n\rightarrow\infty}\mu(A_n\cap B)=0$, for every measurable set $B\subset X$ with $\mu(B)<\infty$.

\end{lemma}

Now, let $\nu$ be a Borel probability measure on $X$ which is equivalent to $\mu$. Let $\{A_n\}$ be a sequence of measurable subsets of $X$ such that $\lim\limits_{n\rightarrow\infty}\nu(gA_n\Delta A_n)=0$, for all $g\in\Gamma$. Lemma \ref{invar} implies that $\lim\limits_{n\rightarrow\infty}\mu((gA_n\Delta A_n)\cap B)=0$, for all $g\in\Gamma$.
 Since condition (2) holds, we get $\lim\limits_{n\rightarrow\infty}\mu(A_n\cap B)=\mu(B)$. 
 
 Using the almost $\Gamma$-invariance of $\{A_n\}$, it follows that $\lim\limits_{n\rightarrow\infty}\mu(A_n\cap C)=\mu(C)$ and hence that $\lim\limits_{n\rightarrow\infty}\mu(A_n^c\cap C)=0$, for all $C\in\mathcal F_B(X)$. 
 Since the action $\Gamma\curvearrowright (X,\mu)$ is ergodic, we have that $\mu(X\setminus(\cup_{g\in\Gamma}gV))=0$.
By combining the last two facts we deduce that $\lim\limits_{n\rightarrow\infty}\mu(A_n^c\cap B)=0$, for every measurable set $B\subset X$ with $\mu(B)<\infty$.
 Applying Lemma \ref{invar} again yields that $\lim\limits_{n\rightarrow\infty}\nu(A_n^c)=0$. This implies that the asymptotically invariant sequence $\{A_n\}$ is trivial.
 \hfill$\blacksquare$

\subsection{Absolute continuity of invariant finitely additive measures} Towards proving implication (2) $\Longrightarrow$ (1) from Theorem \ref{BRtext} we first establish that any $\Gamma$-invariant finitely additive measure on $\mathcal C(G)$ is absolutely continuous with respect to $m_G$. 

\begin{proposition}\label{tarski1} Let $G$ be a l.c.s.c. group and $\Gamma<G$ be a countable dense subgroup such that $\Gamma\cap G_0$ is non-amenable, for any open subgroup $G_0<G$.  Let $\nu:\mathcal C(G)\rightarrow [0,\infty)$ be a $\Gamma$-invariant, finitely additive measure.

If $B\in\mathcal C(G)$ and $|B|=0$, then $\nu(B)=0$.
\end{proposition}

Proposition \ref{tarski1} is proved by adapting an argument due to Tarski, who used the Banach-Tarski paradox to show that any $SO(n+1)$-invariant, finitely additive measure defined on the Lebesgue subsets of $S^n$,  $n\geqslant 2$, is absolutely continuous with respect to the Lebesgue measure (see \cite[Proposition 2.2.12]{Lu94}). 
As such, we need to recall the notion of equidecomposability:

\begin{definition}\label{echidec} Let $\Gamma\curvearrowright X$ be an action of a group $\Gamma$ on a set $X$.
Denote by $\mathcal P(X)$ the power set of $X$ and let $B,C\in\mathcal P(X)$.
We say that $B$ is {\it $\Gamma$-equidecomposable} to $C$ if we can decompose $B=\cup_{i=1}^kB_i$ and $C=\cup_{i=1}^kC_i$ with $B_i\cap B_j=C_i\cap C_j = \emptyset$, for all $1\leqslant i<j\leqslant k$, and there exist $g_1,...,g_k\in\Gamma$ such that $g_iB_i=C_i$, for all $1\leqslant i\leqslant k$. In this case, we write $B\sim C$. 

Also, we write $B\precsim C$ if $B\sim C_0$, for some subset $C_0\subset C$. More generally, if $m$ and $n$ are positive integers, we say that $mB$ is equidecomposable to $nC$ (and write $mB\prec nC$) if we can decompose $B$ in $m$ ways and use translations with elements of $\Gamma$ to rebuild $n$ copies of $C$. Note that the following cancellation law holds: if $mB\sim mC$, then $B\sim C$ (see \cite[page 13]{Lu94}). With this terminology, $B$ is said to be {\it paradoxical} if $B\sim 2B$.
\end{definition}

The proof of Proposition \ref{tarski1} relies on the following two lemmas, the first of which is due to Breuillard and Gelander.

\begin{lemma}\emph{\cite{BG04}}\label{bg}
Let $G$ be a l.c.s.c group and $\Gamma<G$ be a countable dense subgroup. Assume that $\Gamma\cap G_0$ is non-amenable, for every open subgroup $G_0<G$.

Then at least one of the following two conditions holds:
\begin{enumerate}
\item $\Gamma$ contains a non-abelian free subgroup which is not discrete in $G$.
\item There exists a compact subgroup $K<G$ such that $\Gamma\cap K$ is non-amenable. 
\end{enumerate}
\end{lemma}

{\it Proof.} Let $G^0$ be the connected component of the identity in $G$. Then $G/G^0$ is a totally disconnected lcsc group and therefore  admits an open compact subgroup $L$. Since $G_0:=L\cdot G^0$ is an open subgroup of $G$, it is enough to prove the conclusion of the lemma for the inclusion $\Gamma\cap G_0<G_0$. 
Thus,  after replacing $G$ by $G_0$, we may assume that $G$ is a compact extension of a connected group.

By the structure theory of locally compact groups \cite[Theorem 4.6]{MZ55}, we can  find a compact normal subgroup $K<G$ such that $H:=G/K$ is a real Lie group.
By replacing $G$ with an open subgroup, we may assume that $H$ is connected. Denote by $p:G\rightarrow H$ the quotient homomorphism. Let $\frak h$ be the Lie algebra of $H$ and consider the adjoint homomorphism $q:=$Ad$:H\rightarrow GL(\frak h)$. Since $H$ is connected, the kernel of $q$ is equal to the center $Z(H)$ of $H$.

By applying the topological Tits alternative \cite[Theorem 1.3]{BG04} to the inclusion of  $q(p(\Gamma))$ into $GL(\frak{h})\cong GL_{\text{dim(H)}}(\mathbb R)$ (endowed with the standard Euclidean topology), we are in one of two cases: (i) $q(p(\Gamma))$ contains a free dense subgroup, or (ii) $q(p(\Gamma))$ contains an open solvable subgroup $\Sigma$. 

In case (i), there is a free subgroup $\Delta<\Gamma$ such that  $q(p(\Delta))$ is dense in $q(p(\Gamma))$. Since $p(\Gamma)$ is dense in $H$, we get that $p(\Delta)Z(H)$ is dense in $H$. Hence, $p([\Delta,\Delta])=[p(\Delta),p(\Delta)]=[p(\Delta)Z(H),p(\Delta)Z(H)]$ is dense in $[H,H]$ (where as usual, $[\Delta,\Delta]$ denotes the commutator subgroup of $\Delta$). 
Note that $[H,H]$ is not discrete in $H$. Otherwise, since $H$ is connected, $H$ would be abelian. Thus, $p([\Delta,\Delta])$ is not discrete in $H$. Since $K$ is compact, we conclude that $[\Delta,\Delta]$ is a free subgroup of $\Gamma$ which is not discrete in $G$, which proves (1).

In case (ii), we have that $q^{-1}(\Sigma)$ is an open subgroup of $p(\Gamma)$. Since $H$ is connected and $p(\Gamma)<H$ is dense, it follows that $q^{-1}(\Sigma)=p(\Gamma)$.
Since the kernel of $q$ is equal to $Z(H)$, we get that $p(\Gamma)/(q^{-1}(\Sigma)\cap Z(H))\cong\Sigma$. Given that $\Sigma$ is solvable and $Z(H)$ is abelian, we conclude that $p(\Gamma)$ is amenable. Since $\Gamma$ is non-amenable, we derive that $\Gamma\cap K$ is non-amenable, which proves (2). \hfill$\blacksquare$

\begin{lemma}\label{equidec1}
Let $G$ be a lcsc group and $\Gamma<G$ be a countable dense subgroup such that $\Gamma\cap G_0$ is non-amenable, for every open subgroup $G_0<G$.

Then any sets $B,C$ with compact closures and non-empty interiors are $\Gamma$-equidecomposable.
\end{lemma}

{\it Proof.}   
Since $G$ is a l.c.s.c. group, we can find a left-invariant compatible metric $d$ on $G$ which is {\it proper}, in the sense that the closed ball $B_{r}(x):=\{y\in G|d(x,y)\leqslant r\}$ is compact, for every $r>0$ and $x\in G$ (see e.g. \cite{CdH14}). We denote $B_r=B_r(e)$, for all $r>0$.

 The idea of the proof is to show that there exists a bounded set $D\subset G$ which has non-empty interior and is paradoxical. By Lemma \ref{bg} we are in one of two cases.

{\bf Case 1.}  $\Gamma$ contains a non-abelian free subgroup $\Gamma_0$ which is not discrete in $G$.

We claim that $D:=B_1$ is paradoxical. 
Assume by contradiction that $B_1$ is not paradoxical. Then by a theorem of Tarski (see \cite[Corollary 9.2]{Wa93}) we can find a $\Gamma$-invariant, finitely additive measure $\varphi:\mathcal P(G)\rightarrow [0,\infty]$ such that $\varphi(B_1)=1$.
Note that if $r>0$, then we can find $g_1,...,g_n\in\Gamma$ such that $B_r\subset\cup_{i=1}^ng_iB_1$ and $B_1\subset\cup_{i=1}^ng_iB_r$. From this we deduce that $0<\varphi(B_r)<\infty$, for all $r>0$.
We define $\rho:(0,\infty)\rightarrow (0,\infty)$ by letting $\rho(r)=\varphi(B_r)$. Then $\rho$ is an increasing function and hence has at most countably many points of discontinuity. 

Let $s>0$ be a point at which $\rho$ is continuous. 
Next, we define $\Phi:\mathcal P(G)\rightarrow [0,1]$ by letting $\Phi(A)=\varphi(A\cap B_{s})/\varphi(B_{s})$, for every $A\in\mathcal P(G)$. Then $\Phi$ is a finitely additive measure, $\Phi(G)=1$, and we have the following

\begin{equation}\label{echi} \lim\limits_{g\in\Gamma,\; g\rightarrow e}\sup_{A\in\mathcal P(G)}|\Phi(gA)-\Phi(A)|=0. \end{equation}

To see this,  let $g\in G$ and $A\in\mathcal P(G)$. Then  we have that $g^{-1}B_{s}\setminus B_{s}\subset B_{s'}\setminus B_{s}$ and that $B_{s}\setminus g^{-1}B_{s}\subset B_{s}\setminus B_{s''}$, where $s'=s+d(g,e)$ and $s''=\max\{s-d(g,e),0\}$. Since $\rho$ is continuous at $s$, we get that $\lim\limits_{g\rightarrow e}\varphi(g^{-1}B_{s}\setminus B_{s})=0$ and $\lim\limits_{g\rightarrow e}\varphi(B_{s}\setminus g^{-1}B_{s})=0.$
Moreover, if  $g\in\Gamma$, then we have that $|\varphi(gA\cap B_{s})-\varphi(A\cap B_{s})|=|\varphi(A\cap g^{-1}B_{s})-\varphi(A\cap B_{s})|\leqslant \varphi(g^{-1}B_{s}\setminus B_{s})+\varphi(B_{s}\setminus g^{-1}B_{s})$.
The combination of these two facts clearly implies \ref{echi}.

Now, let $a$ and $b$ be two free generators of $\Gamma_0$. Let $g_n\in\Gamma_0\setminus\{e\}$ be a sequence such that $\lim\limits_{n\rightarrow\infty}g_n=e$. Note that for every $g\in\Gamma_0\setminus\{e\}$, at least one of the pairs $\{g,aga^{-1}\}$ and $\{g,bgb^{-1}\}$ generates a copy of $\mathbb F_2$. Therefore, after passing to a subsequence and after eventually swapping $a$ and $b$, we may assume that $g_n$ and $h_n=ag_na^{-1}$ generate a copy of $\mathbb F_2$, for all $n$. Note that $\lim\limits_{n\rightarrow\infty}g_n=\lim\limits_{n\rightarrow\infty}h_n=e$.

We define $X=\sqcup_{n\geqslant 1}X_n$, where $X_n=G$, to be the disjoint union of infinitely many copies of $G$. Let $c$ and $d$ be two free generators  of $\mathbb F_2$.  Let $\lim\limits_n$ be a bounded linear functional on $\ell^{\infty}(\mathbb N)$ which extends the limit.
 We consider an action $\mathbb F_2\curvearrowright X$ given by
$c\cdot x=g_nx$ and $d\cdot x=h_nx$,  for all $n\geqslant 1$ and $x\in X_n.$ Finally, we define  $\Psi:\mathcal P(X)\rightarrow [0,1]$ by 
$\Psi(A)=\lim\limits_n\Phi(A\cap X_n),$ for all $A\in\mathcal P(X).$
Then $\Psi$ is a finitely additive measure and $\Psi(X)=1$. Moreover, $\Psi$ is $\mathbb F_2$-invariant. Indeed, let $A\in\mathcal P(X)$ and write $A=\sqcup_{n\geqslant 1}A_n$, where $A_n\subset X_n$. Since $\lim\limits_{n\rightarrow\infty}g_n=e$, equation \ref{echi} implies that $$\Psi(c\cdot A)=\lim\limits_n\Phi(g_nA_n)=\lim\limits_{n}\Phi(A_n)=\Psi(A).$$

Similarly, if follows that $\Psi(d\cdot A)=\Psi(A)$. Altogether, $\Psi:\mathcal P(X)\rightarrow [0,1]$ is an $\mathbb F_2$-invariant, finitely additive measure with $\Psi(X)=1$. However, since the action $\mathbb F_2\curvearrowright X$ is free,  this contradicts the non-amenability of $\mathbb F_2$. Thus, we conclude that $B_1$ is paradoxical.

{\bf Case 2}.  There exists a compact subgroup $K<G$ such that $\Gamma_0:=\Gamma\cap K$ is non-amenable.

Let $D\in\mathcal C(G)$ be a $K$-invariant set with non-empty interior. Since the left multiplication action $\Gamma_0\curvearrowright D$ is free and $\Gamma_0$ is non-amenable, we get that $D$ is $\Gamma_0$-paradoxical and hence $\Gamma$-paradoxical.

Altogether, we have shown the existence of a subset $D\subset G$ which has non-empty interior and is paradoxical. Let $B,C\subset G$ be two bounded subsets with non-empty interior.

Since $\Gamma<G$ is dense, we can find integers $p,q\geqslant 1$ such that $B\precsim pD$ and $D\precsim qC$. Since $D$ is paradoxical,  we have $2^nD\precsim D$, for all integers $n\geqslant 1$. By combining these facts, we get that $qB\precsim qpD\precsim 2^{qp}D\precsim D\precsim qC$. The cancellation law implies that $B\precsim C$. Similarly, we get that $C\precsim B$, and \cite[Proposition 2.1.2.]{Lu94} implies that $B\sim C$, as desired. \hfill$\blacksquare$

{\bf Proof of Proposition \ref{tarski1}.} Let $\nu:\mathcal C(G)\rightarrow [0,\infty)$ be a $\Gamma$-invariant, finitely-additive measure and $B\in\mathcal C(G)$ such that $|B|=0$. We will show that $\nu(B)=0$. To this end, let $\varepsilon>0$.

Let $C\in\mathcal C(G)$ be a set with non-empty interior such that $\nu(C)<\varepsilon$.
By Lemma \ref{equidec1} we have $B\precsim C$, hence we can find a subset $D\subset C$ such that $B\sim D$. We can therefore decompose $B=\cup_{i=1}^kB_i$ and $D=\cup_{i=1}^kD_i$ such that $B_i\cap B_j=D_i\cap D_j=\emptyset$, for all $1\leqslant i<j\leqslant k$, and there exist $g_1,...,g_k\in\Gamma$ such that $g_iB_i=D_i$, for all $1\leqslant i\leqslant k$. Since $|B|=0$, we get $|D_i|=|B_i|=0$, for all $1\leqslant i\leqslant k$. In particular, $B_i,D_i$ are measurable sets with compact closures, for all $1\leqslant i\leqslant k$. Since $\nu$ is $\Gamma$-invariant, we deduce that $\nu(B)=\sum_{i=1}^k\nu(B_i)=\sum_{i=1}^k\nu(D_i)=\nu(D).$ Therefore,  since $D\subset C$, we conclude that $\nu(B)\leqslant\nu(C)<\varepsilon$ showing that $\nu(B) = 0$ because $\varepsilon > 0$ was arbitrary.
\hfill$\blacksquare$

\subsection{Proof of Theorem \ref{BRtext}} We are now ready to prove Theorem \ref{BRtext}. Let $\mathcal S\subset L^{\infty}_{\text{c}}(G,m_G)$ be the set of functions of the form $\sum_{i=1}^nc_i1_{A_i}$, where $c_1,...,c_n\in\mathbb C$ and $A_1,...,A_n\in\mathcal C(G)$.

{\bf (1) $\Longrightarrow$ (2)}. 
Assume that (1) holds true. Let $\Phi:L^{\infty}_{\text{c}}(G,m_G)\rightarrow\mathbb C$ be a $\Gamma$-invariant, positive linear functional. Then $\nu:\mathcal C(G)\rightarrow [0,\infty)$ given by $\nu(A)=\Phi(1_A)$ is a $\Gamma$-invariant, finitely additive measure. Since (1) holds,  we can find $\alpha\geqslant 0$ such that $\nu(A)=\alpha |A|$, for all $A\in\mathcal C(G)$.
Then we clearly have that $\Phi(f)=\alpha{\int_{G}f\;\text{d}m_G}$, for all $f\in\mathcal S$. 

Let $f\in L^{\infty}_{\text{c}}(G,m_G)$ be a real-valued function. Then $-\|f\|_{\infty}1_B\leqslant f\leqslant \|f\|_{\infty}1_B$, where $B$ denotes the support of $f$. Since $\Phi$ is a positive, we get that
$-\alpha|B|\|f\|_{\infty}\leqslant \Phi(f)\leqslant\alpha|B|\|f\|_{\infty}.$
Moreover, we may find a sequence $\{f_n\}$ of real-valued functions  which belong to $\mathcal S$ and whose support is included in $B$ such that $\lim\limits_{n\rightarrow\infty}\|f-f_n\|_{\infty}=0$. By applying the above inequality to $f-f_n$, it follows that ${\Phi(f)=\lim\limits_{n\rightarrow\infty}\Phi(f_n)=\lim\limits_{n\rightarrow\infty}\alpha\int_{G}f_n\;\text{d}m_G=\alpha\int_{G}f\;\text{d}m_G}.$ This proves condition  (2).
\hfill$\square$

{\bf  (2) $\Longrightarrow$ (1)}. By using Proposition \ref{tarski1}, $\Phi:\mathcal S\rightarrow\mathbb C$ given by $\Phi(f)=\sum_{i=1}^nc_i\nu(A_i)$, for every $f=\sum_{i=1}^nc_i1_{A_i}\in\mathcal S$, is a well-defined $\Gamma$-invariant positive linear functional. Moreover, if  $B$ is the support of $f$, then $|\Phi(f)|\leqslant \nu(B)\|f\|_{\infty}$. By arguing as above, we get that $\Phi$ extends to a $\Gamma$-invariant positive linear functional $\Phi:L^{\infty}_{\text{c}}(G,m_G)\rightarrow\mathbb C$. By applying  (2) to $\Phi$, the conclusion follows.
\hfill$\square$

{\bf  (2) $\Longleftrightarrow$ (3)}. Let $B\in\mathcal C(G)$ be a set with non-empty interior. Then $L^{\infty}_B(G,m_G)=L^{\infty}_{\text{c}}(G,m_G)$ by Remark \ref{CB} and the conclusion follows from Theorem \ref{mean}. \hfill$\square$

{\bf (3) $\Longrightarrow$ (4)}. This follows from the implication (4) $\Longrightarrow$ (5) from Theorem \ref{mean}. \hfill$\square$

{\bf (4) $\Longrightarrow$ (3)}. Assume that $\Gamma\curvearrowright (G,m_G)$ is strongly ergodic and let $\nu$ be a Borel probability measure on $G$ which is equivalent to $m_G$. Let $B\in\mathcal C(G)$ be an open set with compact closure.

Let $\{A_n\}$ be a sequence of measurable subsets of $G$ such that $|A_n\cap B|>0$, for all $n$, and $\lim\limits_{n\rightarrow\infty}|(gA_n\Delta A_n)\cap B|/|A_n\cap B|=0$, for all $g\in\Gamma$.
Since $L^{\infty}_B(G,m_G)=L^{\infty}_{\text{c}}(G,m_G)$,  Theorem \ref{mean} guarantees that in order to prove local spectral gap with respect to $B$, it suffices to show that $\lim\limits_{n\rightarrow\infty}|A_n\cap B|=|B|$.

After passing to a subsequence, we may assume that the limit $\ell=\lim\limits_{n\rightarrow\infty}|A_n\cap B|$ exists. By using the strong ergodicity assumption and Lemma \ref{invar} it is easy to see that if $\ell>0$, then $\ell=|B|$.  
 Thus, in order to derive the conclusion, it suffices to prove that $\ell=0$ leads to a contradiction. We will achieve this by adapting the ``averaging" argument from \cite[Lemma 14]{AJN07}.
 
Now, we let  $C=B^{-1}B$ and $D=BC$.  We claim that  $\sup_n |A_n\cap D|/|A_n\cap B|<\infty$.
Indeed, since $D\in\mathcal C(G)$ and $\cup_{g\in\Gamma}gB=G$, we can find $h_1,...,h_p\in\Gamma$ such that $D\subset\cup_{i=1}^ph_iB$. 
Then we have that $|A_n\cap D|\leqslant\sum_{i=1}^p|h_i^{-1}A_n\cap B|\leqslant p |A_n\cap B|+\sum_{i=1}^p|(h_i^{-1}A_n\Delta A_n)\cap B|$ and the claim follows. Since $\ell=0$, the claim implies that $\lim\limits_{n\rightarrow\infty}|A_n\cap D|=0$ and we can find $\kappa\in (0,1)$ such that 
\begin{equation}\label{ineq} \Big [\frac{|C|}{|A_n\cap D|}\Big]\geqslant \kappa\;\frac{|C|}{|A_n\cap B|},\;\;\text{for all $n$.}\end{equation}

We continue with the following claim:

{\bf Claim.} If $A\subset G$ is a measurable set and $p\geqslant 1$ is an integer, then

$$ \int_{C^p}|(Ag_1^{-1}\cup...\cup Ag_p^{-1})\cap B|\;\text{d}g_1...\text{d}g_p=|C|^p\int_{B}\Big(1-\Big(1-\frac{|A\cap xC|)}{|C|}\Big )^p\Big)\;\text{d}x. $$

{\it Proof of the claim.} 
 Let $q\geqslant 1$ be an integer and put $\beta_q= {\int_{C^q}|Ag_1^{-1}\cap...\cap Ag_q^{-1}\cap B|\;\text{d}g_1...\text{d}g_q}.$
By using Fubini's theorem we have that \begin{align*}\beta_q & =\int_{C^q}\;\Big(\int_B 1_{Ag_1^{-1}}(x)...1_{Ag_q^{-1}}(x)\;\text{d}x\Big)\;\text{d}g_1...\text{d}g_q \\ &=\int_B\;\Big(\int_{C^q}1_{x^{-1}A}(g_1)...1_{x^{-1}A}(g_q)\;\text{d}g_1...\text{d}g_q\Big)\;\text{d}x \\ &=\int_{B}|x^{-1}A\cap C|^q\;\text{d}x=\int_{B}|A\cap xC|^q\;\text{d}x.\end{align*}

Now, the inclusion-exclusion principle gives that the left side of the claimed identity  is equal to ${\alpha_p=\sum_{q=1}^p(-1)^{q-1}{p\choose q}\beta_q\;|C|^{p-q}}$. In combination with the above, the claim follows.\hfill$\square$ 

As one may observe, the claim deals with right translates while the group that we consider is not necessarily unimodular. We will then have to consider the constant $\tau = \sup \{ \Delta(g) \, , \, g \in C\}$, where $\Delta$ is the modular function on $G$. This constant is finite since $C$ has compact closure.

Since ${\lim\limits_{t\rightarrow 0}(1-t)^{\frac{1}{t}}=\frac{1}{e}}$, we can find $\varepsilon\in (0,1)$ such that \begin{equation}\label{e}(1-t)^{[\frac{1}{t}]}<\frac{1}{2},\;\;\;\text{for all $t\in (0,\varepsilon].$}\end{equation} 

Since $\lim\limits_{n\rightarrow\infty}|A_n\cap D|=0$, we can find $N$ such that $|A_n\cap D|/|C|\leqslant\varepsilon$, for all $n\geqslant N$.

Fix $n\geqslant N$ and let ${p_n:=\frac{\vert B \vert}{2\tau \vert C \vert}\Big[\frac{|C|}{|A_n\cap D|}\Big]}$. Also, let $x\in B$. Then we have that $B\subset xC$ and hence $|A_n\cap B|\leqslant |A_n\cap xC|$.

Defining $t_n=|A_n\cap B|/|C|$ we see that $t_n\leqslant |A_n\cap D|/|C|\leqslant\varepsilon$. By \eqref{ineq} we get $p_n\geqslant \frac{\alpha}{t_n}\geqslant \alpha[\frac{1}{t_n}]$, where $\alpha = \frac{\kappa \vert B \vert}{2\tau\vert C \vert}$.
By using \eqref{e} we deduce that \[1-\Big (1-\frac{|A_n\cap xC|}{|C|}\Big)^{p_n}\geqslant 1-(1-t_n)^{p_n}\geqslant 1-(1-t_n)^{\alpha[\frac{1}{t_n}]}>1-\frac{1}{2^{\alpha}}.\]

By combining this inequality with the claim above, we find $g_{n,1},g_{n,2},...,g_{n,p_n}\in C$ such that $\tilde A_n:=\cup_{i=1}^{p_n}A_ng_{n,i}^{-1}$ satisfies \begin{equation}\label{non}(1-\frac{1}{2^{\alpha}})|B|<|\tilde A_n\cap B|, \;\;\;\text{for all $n\geqslant N$}.\end{equation}
Moreover, we also have the upper bound
\begin{equation}\label{non2}
|\tilde A_n\cap B| \leqslant \sum_{i = 1}^{p_n} \vert A_ng_i^{-1} \cap B \vert \leqslant \sum_{i = 1}^{p_n} \vert A_n \cap Bg_i \vert \tau \leqslant  \tau p_n\vert A_n \cap D \vert \leqslant \frac{\vert B \vert}{2}, \;\;\;\text{for all $n\geqslant N$}.
\end{equation}

Finally, let $g\in\Gamma$ and $K\in\mathcal C(G)$. We claim that $\lim\limits_{n\rightarrow\infty}|(g\tilde A_n\Delta\tilde A_n)\cap K|=0$. To this end, note that $g\tilde A_n\Delta\tilde A_n\subset\cup_{i=1}^{p_n}(gA_n\Delta A_n)g_{n,i}^{-1}$. Denote $L=KC$ and notice that $Kg_{n,i}\subset L$. Also, since $B\subset D$, we get that $p_n\leqslant |C|/|A_n\cap B|$. By combining all of these facts we get that  \begin{align*}|(g\tilde A_n\Delta\tilde A_n)\cap K| & \leqslant\sum_{i=1}^{p_n}|(gA_n\Delta A_n)g_{n,i}^{-1}\cap K| \\ & \leqslant p_n\tau |(gA_n\Delta A_n)\cap L| \leqslant \tau |C|\frac{|(gA_n\Delta A_n)\cap L|}{|A_n\cap B|}.\end{align*}

Since $L\in\mathcal C(G)$, we have that $\lim\limits_{n\rightarrow\infty}|(gA_n\Delta A_n)\cap L|/|A_n\cap B|=0$. This proves our claim.

Let $\mu$ be a Borel probability measure on $G$ which is equivalent to $m_G$.
 By combining the claim and Lemma \ref{invar} we get that $\lim\limits_{n\rightarrow}\mu(g\tilde A_n\Delta\tilde A_n)=0$, for all $g\in\Gamma$. Since the action $\Gamma\curvearrowright (G,m_G)$ is strongly ergodic, we conclude that $\lim\limits_{n\rightarrow\infty}\mu(\tilde A_n)(1-\mu(\tilde A_n))=0$. By applying Lemma \ref{invar} again, we get that $\lim\limits_{n\rightarrow\infty}|\tilde A_n\cap B|(|B|-|\tilde A_n\cap B|)=0$. This however contradicts the inequalities \eqref{non} and \eqref{non2}. 
 We have therefore shown that the case $\ell=0$ leads to a contradiction, as desired. \hfill$\blacksquare$

 \section{Proofs of Corollaries \ref{monexp} and \ref{rwalk}}


\subsection{Proof of Corollary \ref{monexp}}
 
We only treat the case when $H$ is trivial. The general case can then be deduced in a similar fashion as in the proof of Corollary \ref{by}.
Assume that the conclusion is false. Then we can find a sequence of non-negligible measurable sets $A_n\subset B$ such that $|A_n|\leqslant\frac{|B|}{2}$ and \begin{equation}\label{ainv} \frac{|(gA_n\cap B)\setminus A_n|}{|A_n|}\rightarrow 0,\;\;\;\text{for every $g\in\Gamma$.}\end{equation}
Our first goal is to show that the sets $A_n$ are equidistributed, in the following sense:
\begin{equation}\label{equidib}\lim_n\frac{|C\cap A_n|}{|A_n|}=\frac{|C|}{|B|},\;\;\;\text{for every $C\in\mathcal C(B)$}.\end{equation}
To this end, we denote by $\mathcal C(B)$ the collection of all measurable subsets of $B$, and define a finitely additive measure $\nu:\mathcal C(B)\rightarrow [0,1]$ by letting $\nu(C)=\lim_{n}\frac{|C\cap A_n|}{|A_n|}$, for every $C\in\mathcal C(B)$.

Let us first show that $\nu(gC)=\nu(C)$, whenever $C\in\mathcal C(B)$ and $g\in\Gamma$ satisfy $gC\subset B$. 
Since \begin{align*} |gC\cap A_n|=|C\cap g^{-1}A_n|=|C\cap (g^{-1}A_n\cap B)|\leqslant |C\cap A_n|+|(g^{-1}A_n\cap B)\setminus A_n|,\end{align*} 
  dividing by $|A_n|$ and using \eqref{ainv} yields that $\nu(gC)\leqslant\nu(C)$. Since the same argument shows that $\nu(C)=\nu(g^{-1}(gC))\leqslant\nu(gC)$, the assertion follows.
  
Next, we claim that $\nu$ extends to a $\Gamma$-invariant finitely additive measure $\tilde\nu:\mathcal C(G)\rightarrow [0,\infty)$. 
  Let $C\in\mathcal C(G)$. Since $\Gamma$ is dense in $G$, we can find a measurable partition $\{C_i\}_{i=1}^k$ of $C$ and  
  $g_1,...,g_p\in\Gamma$ such that $g_iC_i\subset B$, for every $1\leqslant i\leqslant p$. We define $\tilde\nu(C)=\sum_{i=1}^p\nu(g_iC_i)$. To see that $\tilde\nu$ is well-defined, consider another measurable partition $\{D_j\}_{j=1}^q$ of $C$ and  
  $h_1,...,h_q\in\Gamma$ such that $h_jD_j\subset B$, for every $1\leqslant j\leqslant q$. Then for all $i,j$, we have that $g_i(C_i\cap D_j),h_j(C_i\cap D_j)\subset B$,  thus $\nu(g_i(C_i\cap D_j))=\nu(h_j(C_i\cap D_j))$. Using this fact we derive that $$\sum_{i=1}^p\nu(g_iC_i)=\sum_{i=1}^p\big(\sum_{j=1}^q\nu(g_i(C_i\cap D_j))\big)=\sum_{j=1}^q\big(\sum_{i=1}^p\nu(h_j(C_i\cap D_j))\big)=\sum_{j=1}^q\nu(h_jD_j),$$ showing that $\tilde\nu$ is well-defined.
  It is now clear that  $\tilde\nu$ is  finitely additive and extends $\nu$.

Since the action $\Gamma\curvearrowright G$ has local spectral gap and $\tilde\nu(B)=1$, Theorem \ref{BR} implies that   $\tilde\nu(C)=\frac{|C|}{|B|}$, for every $C\in\mathcal C(G)$, thus proving \eqref{equidib}.

Note that  $\alpha:=\inf_n|A_n|>0$. Otherwise, after replacing $\{A_n\}_n$ with a subsequence, we may assume that $\sum_n|A_{n}|<B$. But then $C=B\setminus(\cup_nA_{n})$ would satisfy $\tilde\nu(C)=\nu(C)=0$, while $|C|>0$.  

Define $\eta_n={\bf 1}_{A_n}-\frac{|A_n|}{|B|}{\bf 1}_B\in L^2(B)$.   Further, let $g\in\Gamma$, and note that if $A\subset B$ is a subset, then $gA\setminus A\subset ((gA\cap B)\setminus A)\cup (gB\setminus B)$. Using this, for every $n$ we get that
\begin{align*}\|g\cdot\eta_n-\eta_n\|_2&\leqslant \|1_{gA_n}-1_{A_n}\|_2+\|1_{gB}-1_B\|_2\\&=\sqrt{2|gA_n\setminus A_n|}+\sqrt{2|gB\setminus B|}\\&\leqslant\sqrt{2|(gA_n\cap B)\setminus A_n)|+2|gB\setminus B|}+\sqrt{2|gB\setminus B|}.\end{align*}
In combination with \eqref{ainv} it follows that $\lim_n\|g\cdot\eta_n-\eta_n\|_2\leqslant 2\sqrt{2|gB\setminus B|}$. 

Let $\varepsilon>0$ such that $\sup_{g\in B_{\varepsilon}(1)}|gB\setminus B|<\frac{\alpha}{128}$. By applying Theorem \ref{restricted2} we can find a probability measure $\mu$ supported on $\Gamma\cap B_{\varepsilon}(1)$ and a finite dimensional subspace $V\subset L^2(B)$ such that we have $\|\mu*F\|_2<\frac{\sqrt{\alpha}}{4\sqrt{|B|}}\|F\|_2$, for every $F\in L^2(B)\ominus V$. By using the previous paragraph we have $$\lim_n\|\mu*\eta_n-\eta_n\|_2\leqslant\sup_{g\in B_{\varepsilon}(1)}2\sqrt{2|gB\setminus B|}<\frac{\sqrt{\alpha}}{4}.$$

Moreover, since $\|\eta_n\|_2^2=|A_n|(1-\frac{|A_n|}{|B|})$ and $|A_n|\leqslant\frac{|B|}{2}$, we deduce that ${\frac{\sqrt{\alpha}}{2}}<\|\eta_n\|_2<\sqrt{|B|}$. Using this, we conclude that $\lim_n\|\mu*\eta_n\|_2>\frac{\sqrt{\alpha}}{4}>\frac{\sqrt{\alpha}}{4\sqrt{|B|}}\sup_n\|\eta_n\|_2$. On the other hand, \eqref{equidib} implies that $\eta_n\rightarrow 0$, weakly in $L^2(B)$. This gives a contradiction, and finishes the proof of the main assertion.

Let us prove the moreover assertion. Assuming that this assertion is false, we can  find
a sequence of non-negligible measurable sets $A_n\subset B$ such that $|A_n|\leqslant\frac{|B|}{2}$ and \eqref{ainv} holds for all $g\in\Gamma\cap B_{\varepsilon}(1)$. 

Let $B_0\subset B$ be a non-negligible measurable set with $B_0B_0^{-1}\subset B_{\varepsilon}(1)$. Put $A_n':=A_n\cap B_0$ and assume that $\lim_n\frac{|A_n'|}{|A_n|}>0$. We claim that $\lim_n|A_n'|=|B_0|$.  Indeed, it is easy to see that $\lim_n \frac{|gA_n'\cap B_0\setminus A_n'|}{|A_n'|}=0$, for every $g\in\Gamma\cap B_{\varepsilon}(1)$. Since $gB_0\cap B_0\not=\emptyset$ forces $g\in B_{\varepsilon}(1)$, the last identity holds for every $g\in\Gamma$, and the first part of the proof implies the claim.

Finally, choose a neighborhood $B_1$ of the identity such that $B_1B_1^{-1}\subset B_{\varepsilon}(1)$. Denote by $C$ the set of $x\in B$ such that $\lim_n\frac{|A_n\cap B_1x|}{|A_n|}>0$. By the claim, if $x\in C$, then $\lim_n|A_n\cap B_1x\cap B|=|B_1x\cap B|>0$. Since $B$ is open, it is easy to check that $C$ is both an open and closed subset. This contradicts the connectedness of $B$.
\hfill$\blacksquare$

\subsection{Proof of Corollary \ref{rwalk}} Let $S=\{g_1,...,g_k\}$ be a finite symmetric subset of $G$. 
Recall that the operator $P_S:L^2(B)\rightarrow L^2(B)$ is given by $P_S(F)=\frac{1}{k}\sum_{i=1}^k\big( {\bf 1}_{B\cap g_iB}\;g_i\cdot F+{\bf 1}_{B\setminus g_iB}\;F\big)$. We start by giving a useful formula for $\langle P_S(F),F\rangle$. 
Since ${\bf 1}_{B\setminus g_iB}F=F-{\bf 1}_{B\cap g_iB}F$, we get that $$\langle P_S(F),F\rangle=\|F\|_2^2-\frac{1}{k}\sum_{i=1}^k\big(\|F\|_{2,B\cap g_iB}^2-\langle {\bf 1}_B\;g_i\cdot F,{\bf 1}_{g_iB}\;F\rangle\big).$$
Since $\|g_i\cdot F-F\|_{2,B\cap g_iB}^2=\|F\|_{2,B\cap g_i^{-1}B}^2+\|F\|_{2,B\cap g_iB}^2-2\Re\langle {\bf 1}_B\;g_i\cdot F,{\bf 1}_{g_iB}\;F\rangle$, for all $i$, and $S$ is symmetric, we deduce that 
\begin{equation}\label{P_S}\langle P_S(F),F\rangle=\|F\|_2^2-\frac{1}{2k}\sum_{i=1}^k\|g_i\cdot F-F\|_{2,B\cap g_iB}^2.\end{equation} 
Since $\|g\cdot F-F\|_{2,B\cap gB}\leqslant 2\|F\|_2$, this calculation implies that $P_S$ is symmetric and $\|P_S\|\leqslant 1$. Moreover, if there is $1\leqslant i\leqslant k$ such that $B\cap g_iB=\emptyset$, then $\langle P_S(F),F\rangle\geqslant (-1+\frac{2}{k})\|F\|_2^2$, and therefore the spectrum of $P_S$ satisfies $\sigma(P_S)\subset [-1+\frac{2}{k},1]$. 

Assume that the conclusion of Corollary \ref{rwalk} is false. Thus, the restriction of $P_S$ to $L^2(B)\ominus\mathbb C{\bf 1}_B$ has norm $1$, for every finite symmetric set $S\subset\Gamma$. Fix $g\in\Gamma$ such that $gB\cap B=\emptyset$. It follows that $1\in\sigma({P_S}_{|L^2(B)\ominus\mathbb C{\bf 1}_B})$, for every finite symmetric set $S\subset\Gamma$ which contains $g$.  
Using \eqref{P_S} we conclude that there is a sequence $F_n\in L^2(B)\ominus\mathbb C{\bf 1}_B$ such that $\|F_n\|_2=1$, for all $n$, and $\|g\cdot F_n-F_n\|_{2,B\cap gB}\rightarrow 0$, for every $g\in\Gamma$. 

We claim that $F_n\rightarrow 0$, weakly. Indeed, any weak limit point $F$ of  $\{F_n\}$ satisfies $F(g^{-1}x)=F(x)$, for all $g\in\Gamma$ and almost every $x\in B\cap gB$. It is clear that $F$ can be extended to a $\Gamma$-invariant function $\tilde F\in L^2_{\text{loc}}(G)$. Since $\Gamma<G$ is dense, $\tilde F$ and therefore $F$ must be a constant function. 
Since $F$ has zero integral, we get that $F\equiv 0$, which proves the claim. 

Next, we define a finitely additive measure $\nu:\mathcal C(B)\rightarrow [0,1]$ by letting $\nu(C)=\lim_n\|F_n\|_{2,C}^2$. 
It is easy to check that $\nu(gC)=\nu(C)$, whenever $C\in\mathcal C(B)$ and $g\in\Gamma$ satisfy $gC\subset B$. By repeating the reasoning from the proof of Corollary \ref{monexp} it follows that $\nu(C)=\frac{|C|}{|B|}$, for every $C\in\mathcal C(B)$. 
Since $$\|g\cdot F_n-F_n\|_2\leqslant \|g\cdot F_n-F_n\|_{2,B\cap gB}+\|F_n\|_{2,B\setminus g^{-1}B}+\|F_n\|_{2,B\setminus gB},$$
we get that $\lim_n\|g\cdot F_n-F_n\|_2\leqslant 2\sqrt{\frac{|gB\setminus B|}{|B|}}$, and the proof of Corollary \ref{monexp} gives a contradiction. \hfill$\blacksquare$


\begin{appendix}
\section{Proof of Lemma \ref{BdS}}\label{appendix}

Let $G$ be a connected simple Lie group with trivial center and denote by $\frak g$ its Lie algebra.
The goal of this appendix is to prove Lemma \ref{BdS}. To this end, by relying on results from \cite{dS14,Ta06} and following closely the proof of \cite[Lemma 2.5]{BdS14}, we first prove the following $\ell^2$-flattening lemma. 

In order to do so, it will be more convenient to work with an invariant metric on $G$, rather than the $\|.\|_2$-metric used in the rest of the paper. We therefore fix an Euclidean structure on  $\frak g$,  and endow $G$ with the corresponding left-invariant Riemannian metric, denoted by $d$. 

For further reference, we note that there is a constant $C>1$ such that \begin{equation}\label{metric} C^{-1}\log(1+\|x^{-1}y-1\|_2)\leqslant d(x,y)\leqslant C\log(1+\|x^{-1}y-1\|_2),\;\;\;\text{for all $x,y\in G$}.\end{equation}

Indeed, it suffices to show that $C^{-1}\log(1+\|x-1\|_2)\leqslant d(x,1)\leqslant C\log(1+\|x-1\|_2)$, for $x\in G$. This is clear if $x$ belongs to a small enough neighborhood $V$ of the identity. On the other hand, if $x\not\in V$, then one can easily prove such an inequality by using the $KAK$ decomposition of $G$ (see \cite[Section 3]{LMR00} for details).

 Let $\delta>0$.
We denote by $B(x,\delta)=\{y\in G|d(x,y)<\delta\}$ the {\it open ball of radius $\delta$} centered at $x\in G$,  and by $A^{[\delta]}=\cup_{x\in A}B(x,\delta)$ the {\it $\delta$-neighborhood} of a set $A\subset G$, both with respect to the metric $d$.  For a probability measure $\mu$ on $G$, we denote $\mu_{\delta}=\mu*Q_{\delta}$, where $$Q_{\delta}=\frac{{\bf 1}_{B(1,\delta)}}{|B(1,\delta)|}.$$

\begin{lemma}[$\ell^2$-flattening,\cite{BdS14}]\label{BdS2} Let $G$ be a connected simple Lie group with trivial center. 
Given $\alpha,\kappa>0$, there exists $\varepsilon>0$ such that the following holds for any $\delta>0$ small enough:

Suppose that $\mu$ is a symmetric Borel probability measure on $G$ such that
\begin{enumerate}
\item supp$(\mu)\subset B(1,\varepsilon\log{\frac{1}{\delta}})$,
\item $\|\mu_{\delta}\|_2\geqslant \delta^{-\alpha}$, and
\item $(\mu*\mu)(H^{[\rho]})\leqslant\delta^{-\varepsilon}\rho^{\kappa}$, for all $\rho\geqslant\delta$ and an proper closed connected subgroup $H<G$.
\end{enumerate}
Then $\|\mu_{\delta}*\mu_{\delta}\|_2\leqslant\delta^{\varepsilon}\|\mu_{\delta}\|_2.$

\end{lemma}

{\bf Notation}. We use the notation $O(\varepsilon)$ to denote a positive quantity which is bounded by $C\varepsilon$, for some constant $C>0$ depending only on $G$. We also use the notation $\phi\ll\psi$ for functions $\phi,\psi:(0,\infty)\rightarrow (0,\infty)$ to mean the existence of a constant $C$ depending only on $G$ such that $\phi(\delta)\leqslant C\psi(\delta)$, for any small enough $\delta>0$. If $\phi\ll\psi$ and $\psi\ll\phi$, we write $\phi\simeq\psi$.

\subsection{Ingredients of the proof of Lemma \ref{BdS2}}
The proof of Lemma \ref{BdS2} relies on Bourgain and Gamburd's strategy \cite{BG05,BG06}. In order to implement this strategy, we use de Saxc\'{e}'s product theorem \cite[Theorem 3.9]{dS14}. Recall that if $A$ is a subset of $G$ and $\delta>0$, then $N(A,\delta)$ denotes the least number of open balls of radius $\delta$ needed to cover $A$.

\begin{theorem}[product theorem, \cite{dS14}]\label{dS} Let $G$ be a simple Lie group of dimension $d$. Then there exists a neighborhood $U$ of the identity in $G$ such that the following holds. 

Given $\alpha\in (0,d)$ and $\kappa>0$, there exist $\varepsilon_0=\varepsilon_0(\alpha,\kappa)$ and $\tau=\tau(\alpha,\kappa)>0$ such that, for any $\delta>0$ small enough, if $A\subset U$ is a set satisfying
\begin{enumerate}
\item $N(A,\delta)\leqslant\delta^{-d+\alpha-\varepsilon_0}$,
\item $N(A,\rho)\geqslant\rho^{-\kappa}\delta^{\varepsilon_0}$, for all $\rho\geqslant\delta$, and
\item $N(AAA,\delta)\leqslant\delta^{-\varepsilon_0}N(A,\delta)$,
\end{enumerate}
then $A$ is contained in a neighborhood of size $\delta^{\tau}$ of a proper closed connected subgroup of $G$.
\end{theorem}

To prove Lemma \ref{BdS2}, we will also need Tao's non-commutative Balog-Szemer\'{e}di-Gowers Lemma \cite[Theorem 6.10]{Ta06}. If $A,B$ are subsets of a metric group $G$ and $\delta>0$, then the {\it $\delta$-multiplicative energy} $E_{\delta}(A,B)$ is defined as $E_{\delta}(A,B)=N(\{(a,b,a',b')\in A\times B\times A\times B\;|\; d(ab,a'b')\leqslant\delta\},\delta).$
 The following inequality will be used in the proof of Lemma \ref{BdS2} 

$$E_{\delta}(A,B)\gg\delta^{-3d}\|1_A*1_B\|_2^2.$$

Indeed, define $f:A\times B\times A\rightarrow G$ by $f(a,b,a')=a'^{-1}ab$, $S=\{(a,b,a')\in A\times B\times A|f(a,b,a')\in B\}$, and $T=\{(a,b,a',b')\in A\times B\times A\times B|d(b',f(a,b,a'))\leqslant\delta\}$. Then we have $\|1_A*1_B\|_2^2=|S|$ and 
$\delta^d|S|\ll |T|\ll \delta^{4d}N(T,\delta)=\delta^{4d}E_{\delta}(A,B),$ which together prove the desired inequality.

Recall that a subset $H\subset G$ is called a {\it $K$-approximative subgroup}, for some $K>1$, if it is symmetric and there is a symmetric set $X\subset HH$ of cardinality at most $K$ such that $HH\subset XH$.

\begin{theorem}[non-commutative Balog-Szemer\'{e}di-Gowers lemma,\cite{Ta06}]\label{ta} Let $G$ be a Lie group endowed with a left-invariant Riemannian metric. 
Then there exist constants $c>0$ and $R>0$ such that the following holds for any $\delta\in (0,1)$ and $K\geqslant 2$. 

Suppose that $A,B$ are non-empty subsets of $G$ contained in $B(1,1)$ such that $$E_{\delta}(A,B)\geqslant\frac{1}{K}N(A,\delta)^{\frac{3}{2}}N(B,\delta)^{\frac{3}{2}}.$$

Then there is a $K^c$-approximate subgroup $H$ of $G$ and elements $x,y\in G$ such that
\begin{itemize}
\item $N(H,\delta)\leqslant K^c\;N(A,\delta)^{\frac{1}{2}}N(B,\delta)^{\frac{1}{2}}$,
\item $N(xH\cap A,\delta)\geqslant K^{-c}N(A,\delta)$, 
\item $N(Hy\cap B,\delta)\geqslant K^{-c}N(B,\delta)$, 
\item $x,y\in B(1,R)$ and $H\subset B(1,R)$.
\end{itemize}
\end{theorem}

A final ingredient in the proof of Lemma \ref{BdS2}  is an approximation of the measure $\mu_{\delta}$ by dyadic level sets \cite{LdS13}. A family of sets $\{A_i\}_{i\in I}$ is called  {\it essentially disjoint} if there is a constant $C$ such that the intersection of more than $C$ distinct sets $A_i$ is empty. 

\begin{lemma}\emph{\cite{LdS13}}\label{dyadic} Let $\mu$ be a Borel probability measure on $G$.
Let $\delta\in (0,1)$ and $\mathcal C$ be a maximal $\delta$-separated subset of $G$.
Define $\mathcal C_0=\{x\in\mathcal C|0<\mu_{2\delta}(x)\leqslant 1\}$ and $\mathcal C_i=\{x\in\mathcal C|2^{i-1}<\mu_{2\delta}(x)\leqslant 2^{i}\}$, for $i\geqslant 1$. For $i\geqslant 0$, let $A_i=\cup_{x\in\mathcal C_i}B(x,\delta)$. 
Then we have the following

\begin{enumerate}
\item at most $O(1)\log{\frac{1}{\delta}}$ of the sets $A_i$ are non-empty,
\item $A_i$ is an essentially disjoint union of balls of radius $\delta$, for all $i\geqslant 0$, 
\item ${\mu_{\delta}\ll\sum_{i\geqslant 0}2^i1_{A_i}}$ and ${\sum_{i>0}2^i1_{A_i}\ll\mu_{3\delta}}$.

\end{enumerate}
\end{lemma}
{\it Proof.} 
These assertions follow from \cite[Lemma 4.4]{LdS13}, but for completeness, we include a proof.
Since $|\mu_{2\delta}(x)|\leqslant\frac{1}{|B(1,2\delta)|}\ll\delta^{-d}$, for any $x\in G$,  (1) is immediate.
 Note that the balls $\{B(x,\delta)\}_{x\in\mathcal C}$ cover $G$, while the balls $\{B(x,\frac{\delta}{2})\}_{x\in\mathcal C}$ are disjoint. Since a ball of radius $3\delta$ can contain at most $O(1)$ disjoint balls of radius $\frac{\delta}{2}$,   the balls $\{B(x,\delta)\}_{x\in\mathcal C}$ are essentially disjoint. This implies (2). 

To prove (3), let $y\in G$. If  $x\in G$ is such that $y\in B(x,\delta)$, then  $$\mu_{\delta}(y)=\frac{\mu(B(y,\delta))}{|B(1,\delta)|}\leqslant\frac{\mu(B(x,2\delta))}{|B(1,\delta)|}\ll\frac{\mu(B(x,2\delta))}{|B(1,2\delta)|}=\mu_{2\delta}(x),$$ and similarly $\mu_{2\delta}(x)\ll\mu_{3\delta}(y).$ Assuming $\mu_{\delta}(y)>0$, let $x\in\mathcal C$ with $y\in B(x,\delta)$. Then  $\mu_{2\delta}(x)>0$, hence $x\in\mathcal C_i$ and $y\in A_i$, for some $i\geqslant 0$. This implies that $\mu_{\delta}(y)\ll\mu_{2\delta}(x)1_{A_i}(y)\leqslant 2^i1_{A_i}(y).$

On the other hand, if  $y\in A_i$, for  $i>0$, then there is $x\in\mathcal C$ such that $y\in B(x,\delta)$ and $\mu_{2\delta}(x)>2^{i-1}$. Hence, $2^i1_{A_i}(y)=2^i<2\mu_{2\delta}(x)\ll\mu_{3\delta}(y)$.  
Since  the balls $\{B(x,\delta)\}_{x\in\mathcal C}$ are essentially disjoint,   $y$ belongs to $O(1)$ of the sets $\{A_i\}_{i>0}$  and it follows that $\sum_{i>0}2^i1_{A_i}(y)\ll\mu_{3\delta}(y).$
This proves (3). \hfill$\blacksquare$

\begin{lemma}\label{flat} Let $a>0$ and $\mu$ be a Borel probability measure on $G$. 

Then
 ${\sup_{\delta<\delta'<1}\|\mu_{\delta'}\|_2\ll\|\mu_{\delta}\|_2}$ and $\|\mu_{\delta}\|_2\simeq\|\mu_{a\delta}\|_2$.

 \end{lemma}

{\it Proof}. 
Let $1>\delta'>\delta>0$. Then $ |B(1,\delta)| 1_{B(1,\delta')}\leqslant1_{B(1,\delta)}*1_{B(1,\delta+\delta')}$ and thus

$$Q_{\delta'}\leqslant\frac{|B(1,\delta+\delta')|}{|B(1,\delta')|}\;Q_{\delta}*Q_{\delta+\delta'}\ll Q_{\delta}*Q_{\delta+\delta'}.$$

We further get that $\|\mu*Q_{\delta'}\|_2\ll\|(\mu*Q_{\delta})*Q_{\delta+\delta'}\|_2\leqslant \|\mu*Q_{\delta}\|_2$, thus proving the first inequality. To prove the second inequality we may assume $a>1$. Then for any $x\in G$ we have

 $$\mu_{a\delta}(x)=\frac{\mu(B(x,a\delta))}{|B(1,a\delta)|}\geqslant
 \frac{|B(1,\delta)|}{|B(1,a\delta)|}\;\frac{\mu(B(x,\delta))}{|B(1,\delta)|}\gg\mu_{\delta}(x).$$ 
 
 Thus, $\|\mu_{a\delta}\|_2\gg\|\mu_{\delta}\|_2$. Since $\|\mu_{a\delta}\|_2\ll\|\mu_{\delta}\|_2$ by the above, we are done. \hfill$\blacksquare$

\subsection{Proof of Lemma \ref{BdS2}} 
Assume that $\|\mu_{\delta}*\mu_{\delta}\|_2>\delta^{\varepsilon}\|\mu_{\delta}\|_2$, for some $\varepsilon>0$. 
Following closely the proof of \cite[Lemma 2.5]{BdS14}, we will reach a contradiction for any $\varepsilon$ small enough.

Let $U$ be the neighborhood of $1\in G$ provided by Theorem \ref{dS}, and let $0<r<1$ with $B(1,r)\subset U$. 
Let $R$ be the constant given by Theorem \ref{ta}, and $\{A_i\}_{0\leqslant i\ll\log{\frac{1}{\delta}}}$ be the sets given by Lemma \ref{dyadic}. Let $C>1$ the constant appearing in inequality \eqref{metric}. By using \eqref{metric} one easily checks that \begin{equation}\label{d^{-1}}d(x^{-1},y^{-1})\leqslant Ce^{2Cd(x,1)}(e^{Cd(x,y)}-1),\;\;\;\text{for any $x,y\in G$}.\end{equation}

Assume that $3C\varepsilon<1$. Let $0\leqslant i\ll\log{\frac{1}{\delta}}$.
 Since supp$(\mu)\subset B(1,\varepsilon\log{\frac{1}{\delta}})$, we get that $A_i$ is contained in $B(1,\varepsilon\log{\frac{1}{\delta}}+3\delta)$. 
 Since  $|B(1,\varepsilon\log{\frac{1}{\delta}}+3\delta)|=\delta^{-O(\varepsilon)}$ and $|B(1,\frac{\delta^{3C\varepsilon}}{8})|=\delta^{O(\varepsilon)}$, $A_i$ can be covered by at most $\delta^{-O(\varepsilon)}$ sets of diameter at most $\frac{\delta^{3C\varepsilon}}{2}$.
Since $A_i$ is a union of balls of radius $\delta$, and $\delta\leqslant\frac{\delta^{3C\varepsilon}}{2}$, for $\delta$ small enough,
we can decompose  $A_i=\cup_{k=1}^{\delta^{-O(\varepsilon)}}A_{i,k},$  where 
  each set $A_{i,k}$ is the union of some of the balls of radius $\delta$ that make up $A_i$, and has diameter at most $\delta^{3C\varepsilon}$.
Moreover, by \eqref{d^{-1}} the diameter of $A_{i,k}^{-1}$ is at most $Ce^{2C(\varepsilon\log{\frac{1}{\delta}}+3\delta)}(e^{C\delta^{3C\varepsilon}}-1)\simeq {\delta}^{C\varepsilon}$. Thus, for $\delta$ small enough, $A_{i,k}$ has diameter at most $\min\{1,r\}$ and $A_{i,k}^{-1}$ has diameter at most $1$, for all $k$.
 
 Before continuing, let us also note that \begin{equation}\label{N} \delta^{O(\varepsilon)}N(A,\delta)\leqslant N(Ah,\delta)\leqslant\delta^{-O(\varepsilon)}N(A,\delta),\end{equation}  for every set $A\subset B(1,\varepsilon\log{\frac{1}{\delta}}+3\delta)$ and any $h\in B(1,\varepsilon\log{\frac{1}{\delta}}+3\delta+1)$. Indeed, by using \eqref{d^{-1}} it is immediate that $B(1,\delta^{1+O(\varepsilon)})\subset h^{-1}B(1,\delta)h\subset B(1,\delta^{1-O(\varepsilon)})$, which easily implies \eqref{N}.
 
Since ${\mu_{\delta}\ll\sum_{i}2^i1_{A_i}\leqslant\sum_{i,k} 2^i1_{A_{i,k}}}$, we get that $$\delta^{\varepsilon}\|\mu_{\delta}\|_2\leqslant\|\mu_{\delta}*\mu_{\delta}\|_2\leqslant\underset{1\leqslant k,l\leqslant \delta^{-O(\varepsilon)}}{\sum_{0\leqslant i,j\ll\log{\frac{1}{\delta}}}}\|2^i1_{A_{i,k}}*2^j1_{A_{j,l}}\|_2.$$

Since the sum on the right contains $\delta^{-O(\varepsilon)}$ terms,  we can find $0\leqslant i,j\ll\log\frac{1}{\delta}$ and $1\leqslant k,l\leqslant \delta^{-O(\varepsilon)}$ such that if we denote $A_i'=A_{i,k}$ and $A_j'=A_{j,l}$, then $$\|2^i1_{A_i'}*2^j1_{A_j'}\|_2\geqslant\delta^{O(\varepsilon)}\|\mu_{\delta}\|_2.$$

We claim that $i>0$. Note that $\|1_{A_0}\|_1=|A_0|\leqslant |B(1,\varepsilon\log{\frac{1}{\delta}}+3\delta)|=\delta^{-O(\varepsilon)}$ and  $\|1_{A_0}\|_2=\delta^{-O(\varepsilon)}$. Moreover, if $j>0$, then $\|2^j1_{A_j'}\|_1\leqslant \|2^j1_{A_j}\|_1\ll\|\mu_{3\delta}\|_1=1$ (by Lemma \ref{dyadic}). Young's inequality gives that $\|1_{A_0}*2^j1_{A_j'}\|_2\leqslant \|1_{A_0}\|_2\|2^j1_{A_j'}\|_1\leqslant\delta^{-O(\varepsilon)}$, for any $j\geqslant 0$. Since $\delta^{O(\varepsilon)}\|\mu_{\delta}\|_2\geqslant\delta^{O(\varepsilon)-\alpha}$, we cannot have $i=0$, provided $\varepsilon>0$ is small enough. Similarly, we must have that $j>0$.

Since $\|2^{j}1_{A_j'}\|_2\leqslant \|2^j1_{A_j}\|_2\ll\|\mu_{3\delta}\|_2\simeq\|\mu_{\delta}\|_2$ (by Lemma \ref{flat}), Young's inequality  gives that $$\delta^{O(\varepsilon)}\|\mu_{\delta}\|_2\leqslant\|2^i1_{A_i'}*2^j1_{A_j'}\|_2\leqslant \|2^i1_{A_i'}\|_1\|2^j1_{A_j'}\|_2\leqslant 2^i|A_i'\|\mu_{\delta}\|_2.$$ 

This implies that \begin{equation}\label{A1}2^i|A_i'|=\delta^{O(\varepsilon)}\;\;\;\;\text{and similarly}\;\;\;\;2^j|A_j'|=\delta^{O(\varepsilon)}.\end{equation}

Next, since $2^i|A_i'|\leqslant 2^i|A_i|\ll\|\mu_{3\delta}\|_1=1$, we deduce that \begin{align*}\delta^{O(\varepsilon)}\|\mu_{\delta}\|_2\leqslant\|2^i1_{A_i'}*2^j1_{A_j'}\|_2&\leqslant\|2^i1_{A_i'}\|_2\|2^j1_{A_j'}\|_1\ll2^i|A_i'|^{\frac{1}{2}}\leqslant|A_i'|^{-\frac{1}{2}}.\end{align*}

Using that $\|\mu_{\delta}\|_2\geqslant\delta^{-\alpha}$, it follows that $|A_i'|\leqslant\delta^{\frac{\alpha}{2}-O(\varepsilon)}.$
Since $A_i'$ is an essentially disjoint union of balls of radius $\delta$, we have  $|A_i'|\simeq\delta^dN(A_i',\delta)$. Altogether, we deduce that \begin{equation}\label{A2}N(A_i',\delta)\leqslant \delta^{-d+\frac{\alpha}{2}-O(\varepsilon)}\;\;\;\text{and similarly}\;\;\;N(A_j',\delta)\leqslant \delta^{-d+\frac{\alpha}{2}-O(\varepsilon)}.\end{equation}

By combining the inequalities $2^i|A_i'|^{\frac{1}{2}}=\|2^i1_{A_i'}\|_2\leqslant\|\mu_{3\delta}\|_2\simeq\|\mu_{\delta}\|_2$ and $2^i|A_i'|\ll 1$ with the fact that $|A_i'|\simeq\delta^{d}N(A_i',\delta)$, we derive that

\begin{align*}\|1_{A_i'}*1_{A_j'}\|_2^2  \geqslant\delta^{O(\varepsilon)}2^{-2i-2j}\|\mu_{\delta}\|_2^2 & \geqslant \delta^{O(\varepsilon)}2^{-i-j}|A_i'|^{\frac{1}{2}}|A_j'|^{\frac{1}{2}}\\&=\delta^{O(\varepsilon)}(2^i|A_i'|)^{-1}|A_i'|^{\frac{3}{2}}(2^{j}A_j')^{-1}|A_j'|^{\frac{3}{2}}\\&\geqslant\delta^{3d+O(\varepsilon)}N(A_i',\delta)^{\frac{3}{2}}N(A_j',\delta)^{\frac{3}{2}}.
\end{align*}

Since $A_i'$ and ${A_j'}^{-1}$ have diameters at most $1$, we can find $g,h\in B(1,\varepsilon\log{\frac{1}{\delta}}+3\delta+1)$ such that $gA_i'\subset B(1,1)$ and $A_j'h\subset B(1,1)$. On the other hand, combining the last inequality with \ref{N} yields

\begin{align*}E_{\delta}(gA_i',A_j'h)\gg\delta^{-3d}\|1_{gA_i'}*1_{A_j'h}\|_2^2 &=\delta^{-3d}\|1_{A_i'}*1_{A_j'}\|_2^2\\&\geqslant\delta^{O(\varepsilon)}N(A_i',\delta)^{\frac{3}{2}}N(A_j',\delta)^{\frac{3}{2}}\\&\geqslant\delta^{O(\varepsilon)}N(gA_i',\delta)^{\frac{3}{2}}N(A_j'h,\delta)^{\frac{3}{2}}\end{align*}

By applying Theorem \ref{ta} to $gA_i'$ and $A_j'h$, we deduce the existence of a $\delta^{-O(\varepsilon)}$-approximate subgroup $H\subset B(1,R)$ and  elements $z,t\in B(1,R)$ such that $N(H,\delta)\leqslant \delta^{-O(\varepsilon)}N(gA_i',\delta)^{\frac{1}{2}}N(A_j'h,\delta)^{\frac{1}{2}}$, $N(zH\cap gA_i',\delta)\geqslant\delta^{O(\varepsilon)}N(gA_i',\delta)$ and $N(Ht\cap A_j'h,\delta)\geqslant\delta^{O(\varepsilon)}N(A_j'h,\delta)$. Let $v=g^{-1}z$ and $w=th^{-1}$. By using \eqref{N} we further get that

\begin{equation}\label{A3} N(H,\delta)\leqslant \delta^{-O(\varepsilon)}N(A_i',\delta)^{\frac{1}{2}}N(A_j',\delta)^{\frac{1}{2}},\end{equation}

\begin{equation}\label{A4} N(vH\cap A_i',\delta)\geqslant\delta^{O(\varepsilon)}N(A_i',\delta)\;\;\;\text{and}\;\;\;N(Hw\cap A_j',\delta)\geqslant\delta^{O(\varepsilon)}N(A_j',\delta).\end{equation}

 The next claim allows us to replace $H$ with its $4\delta$-neighborhood:

{\bf Claim.} $\tilde H:=H^{[4\delta]}$ satisfies the following:

\begin{enumerate}
\item $\mu_{\delta}(v\tilde H\cap A_i')\geqslant\delta^{O(\varepsilon)}$.
\item $N(\tilde H^2,\delta)\leqslant N(\tilde H^6,\delta)\leqslant\delta^{-O(\varepsilon)}N(\tilde H,\delta).$
\end{enumerate}

{\it Proof of the claim}. (1) Recall that  there is a subset $\mathcal C_i'\subset\mathcal C_i$ such that $A_i'=\cup_{x\in\mathcal C_i'}B(x,\delta)$. Since  $N(vH\cap A_i',\delta)\geqslant\delta^{O(\varepsilon)}N(A_i',\delta)$ by \eqref{A4}, we get that $vH$ intersects at least $\delta^{O(\varepsilon)}N(A_i',\delta)$ of the balls $\{B(x,\delta)\}_{x\in\mathcal C_i'}$. Thus, $v\tilde H\cap A_i'=(vH)^{[4\delta]}\cap A_i'$ contains at least $\delta^{O(\varepsilon)}N(A_i',\delta)$ of the balls $\{B(x,3\delta)\}_{x\in\mathcal C_i'}$. 
Since the balls $\{B(x,\delta)\}_{x\in\mathcal C_i'}$ and hence the balls $\{B(x,3\delta)\}_{x\in\mathcal C_i'}$ are essentially disjoint,  $v\tilde H\cap A_i'$ must contain at least $\delta^{O(\varepsilon)}N(A_i',\delta)$ disjoint balls from the collection $\{B(x,3\delta)\}_{x\in\mathcal C_i'}$.

On the other hand, for every $x\in\mathcal C_i$ we have that $\mu_{2\delta}(x)>2^{i-1}$ and hence $$\mu_{\delta}(B(x,3\delta))\geqslant\mu(B(x,2\delta))=|B(1,2\delta)|\mu_{2\delta}(x)\gg\delta^d2^i.$$

Since $\delta^dN(A_i',\delta)\simeq |A_i'|$, and $2^i|A_i'|=\delta^{O(\varepsilon)}$ by \ref{A1}, we finally derive that $$\mu_{\delta}(v\tilde H\cap A_i)\geqslant\delta^{O(\varepsilon)}N(A_i',\delta)\delta^d2^i\simeq\delta^{O(\varepsilon)}2^i|A_i'|=\delta^{O(\varepsilon)}.$$

(2) 
Let $X\subset G$ be a set of cardinality $\delta^{-O(\varepsilon)}$ such that $HH\subset HX$. 
Since $H\subset B(1,R)$ and $R$ is an absolute constant, by using \eqref{d^{-1}} it follows that $\tilde H\tilde H\subset (HH)^{[O(1)\delta]}$, and thus $\tilde H\tilde H\subset H^{[D\delta]}X$, for some constant $D>1$. Let $S\subset H$ be a maximal set such that the balls $\{B(x,\delta)\}_{x\in S}$ are disjoint. Then $\cup_{x\in S}B(x,\delta)\subset H^{[\delta]}$ and $H\subset\cup_{x\in S}B(x,2\delta)$, therefore $H^{[D\delta]}\subset\cup_{x\in S}B(x,(D+2)\delta)$. Let $Y$ be a set of cardinality $O(1)$ with $B(1,(D+2)\delta)\subset B(1,\delta)Y$. Then $B(x,(D+2)\delta)\subset B(x,\delta)Y$, for any $x\in G$, implying that $H^{[D\delta]}\subset H^{[\delta]}Y$. Altogether, we get that $\tilde H\tilde H\subset H^{[\delta]}Z\subset \tilde HZ$, where $Z=YX$. Since $\tilde H$ is symmetric, we have $\tilde H\tilde H\subset Z^{-1}\tilde H$. Since the cardinality of $Z$ is $\delta^{-O(\varepsilon)}$,  (2) follows.
\hfill$\square$

Let us show that for $\varepsilon>0$ small enough, the set $\tilde H^2\cap U$ satisfies the assumptions of Theorem \ref{dS}.
Firstly, using (2) above in combination with \eqref{A2} and \eqref{A3} we get that \begin{equation}\label{A5}N(\tilde H^2\cap U,\delta)\leqslant N(\tilde H^2,\delta)\leqslant\delta^{-d+\frac{\alpha}{2}-O(\varepsilon)}.\end{equation}

Secondly, since $A_i'$ has diameter at most $r$, we get that ${A_i'}^{-1}{A_i'}\subset B(1,r)\subset U$. This implies that $(v\tilde H\cap A_i')^{-1}(v\tilde H\cap A_i')\subset \tilde H^2\cap U$. In combination with (1) from the claim, it follows that \begin{equation}\label{mumu}\check{\mu}_{\delta}*\mu_{\delta}(\tilde H^2\cap U)\geqslant\check{\mu}_{\delta}((v\tilde H\cap A_i')^{-1})\mu_{\delta}(v\tilde H\cap A_i')=\mu_{\delta}(v\tilde H\cap A_i')^2\geqslant\delta^{O(\varepsilon)}.\end{equation}

Let $1\geqslant\rho\geqslant\delta$ and $x\in G$ such that $B(x,\rho)\cap U\not=\emptyset$. Then $B(1,\delta)B(x,\rho)B(1,\delta)\subset B(x,\rho+O(1)\delta)$.
Since $B(x,\rho+O(1)\delta)$ is contained in the $(\rho+O(1)\delta)$-neighborhood of a proper closed connected subgroup of $G$,  the hypothesis implies that $\mu*\mu(B(x,\rho+O(1)\delta))\ll\delta^{-\varepsilon}\rho^{\kappa}$. By using the fact that $\mu$ and $Q_{\delta}$ are symmetric, we get that 

\begin{align*} \check{\mu}_{\delta}*\mu_{\delta}(B(x,\rho))&=(Q_{\delta}*\mu*\mu*Q_{\delta})(B(x,\rho))\\&\leqslant \mu*\mu(B(1,\delta)B(x,\rho)B(1,\delta))\ll\delta^{-\varepsilon}\rho^{\kappa}.\end{align*}

Since this holds for any ball of radius $\rho$ that intersects $U$, in combination with \eqref{mumu} it gives that 

\begin{equation}\label{A6} N(\tilde H^2\cap U,\rho)\geqslant\delta^{O(\varepsilon)}\rho^{-\kappa},\;\;\;\text{for all $1\geqslant\rho\geqslant\delta$}.\end{equation}

Thirdly, let $C\subset\tilde H$ be a set of diameter at most $r$ and fix $x\in C$. Since $B(1,r)\subset U$, we get that $x^{-1}C\subset C^{-1}C\subset\tilde H^2\cap U$. Thus, $N(C,\delta)=N(x^{-1}C,\delta)\leqslant N(\tilde H^2\cap U,\delta)$. On the other hand, since $H\subset B(1,R)$, we can cover $\tilde H$ by $O(1)$  sets of diameter at most $r$. Altogether, we conclude that $N(\tilde H,\delta)\ll N(\tilde H^2\cap U,\delta)$. Using (2), it follows that

\begin{equation}\label{A7} N((\tilde H^2\cap U)^3,\delta)\leqslant N(\tilde H^6,\delta)\leqslant\delta^{-O(\varepsilon)}N(\tilde H,\delta)\leqslant\delta^{-O(\varepsilon)}N(\tilde H^2\cap U,\delta).  \end{equation}

Equations \eqref{A5}, \eqref{A6}, \eqref{A7} together guarantee that we are in position to apply Theorem \ref{dS} to $\tilde H^2\cap U$. Thus, there is a proper closed connected subgroup $L<G$ such that $\tilde H^2\cap U\subset L^{[\delta^{\tau}]}$. 
On the other hand, by using \eqref{mumu} and reasoning similarly to the above we conclude that \begin{align*}\delta^{O(\varepsilon)}\leqslant\check{\mu}_{\delta}*\mu_{\delta}(\tilde H^2\cap U)&\leqslant\check{\mu}_{\delta}*\mu_{\delta}(L^{[\delta^{\tau}]}\cap U)\\&\leqslant \mu*\mu(B(1,\delta)(L^{[\delta^{\tau}]}\cap U)B(1,\delta))\leqslant\mu*\mu(L^{[\delta^{\tau}+O(1)\delta]}).\end{align*}

Since the hypothesis implies that $\mu*\mu(L^{[\delta^{\tau}+O(1)\delta]})\leqslant\delta^{-\varepsilon}(\delta^{\tau}+O(1)\delta)^{\kappa}$,
it is now clear that choosing $\varepsilon>0$ small enough yields a contradiction. \hfill$\blacksquare$

\subsection{Proof of Lemma \ref{BdS}} Let $\alpha,\kappa>0$. Let $\varepsilon>0$ such that the conclusion of Lemma \ref{BdS2} holds.
 Let $\mu$  be a symmetric Borel probability measure on $G$. Then we have
 
 \begin{itemize}
\item[(a)] $\|\mu*P_{\delta}\|_2\simeq \|\mu*Q_{\delta}\|_2.$ 
\end{itemize}

Moreover, assuming that supp$(\mu)\subset B_{\delta^{-\beta}}(1)$, for some $\beta>0$, we have that

\begin{itemize}
\item[(b)] $\mu*Q_{\delta}\ll\delta^{-O(\beta)}\;Q_{\delta^{1-O(\beta)}}*\mu*Q_{\delta^{1-O(\beta)}}$.
\item[(c)] $A^{[\rho]}\cap$ supp$(\mu*\mu)\subset A^{(\delta^{-O(\beta)}\rho)}$, for any set $A\subset G$ and every $1>\rho>\delta$. 

\end{itemize}

Indeed,  \ref{metric} implies that $B(1,\frac{\delta}{2C})\subset B_{\delta}(1)\subset B(1,C\delta)$, for small $\delta>0$.
 This readily gives that $$P_{\delta}=\frac{1_{B_{\delta}(1)}}{|B_{\delta}(1)|}\leqslant\frac{|B(1,C\delta)|}{|B(1,\frac{\delta}{2C})|} Q_{C\delta}\ll Q_{C\delta}$$ and similarly ${Q_{\frac{\delta}{2C}}\ll P_{\delta}}.$ Therefore,  $\|\mu*Q_{\frac{\delta}{2C}}\|_2\ll \|\mu*P_{\delta}\|\ll \|\mu*Q_{C\delta}\|_2$. On the other hand,  Lemma \ref{flat} implies that $\|\mu*Q_{\frac{\delta}{2C}}\|_2\simeq\|\mu*Q_{\delta}\|_2\simeq\|\mu*Q_{C\delta}\|_2$. This altogether proves (a).

By \eqref{metric} we have that $B_{\delta^{-\beta}}(1)\subset B(1,2C\beta\log{\frac{1}{\delta}})$, hence $\mu$ and $\mu*\mu$ are supported on $B(1,4C\beta\log{\frac{1}{\delta}})$.
To prove (b), we may therefore assume that $\mu=\delta_{x}$, for some $x\in B(1,4C\beta\log{\frac{1}{\delta}})$. Using \eqref{d^{-1}} we get that  $B(1,4C\beta\log{\frac{1}{\delta}})B(1,\delta)B(1,4C\beta\log{\frac{1}{\delta}})\subset B(1,\delta^{1-O(\beta)})$, hence $\delta_x*1_{B(1,\delta)}*\delta_{x^{-1}}\leqslant 1_{B(1,\delta^{1-O(\beta)})}.$
Since $|B(1,\delta^{1-O(\beta)}|\ll\delta^{-O(\beta)}|B(1,\delta)|$, this implies that $\mu*Q_{\delta}\ll\delta^{-O(\beta)}Q_{\delta^{1-O(\beta)}}*\mu$. 
Similarly, we have that $Q_{\delta}\ll \delta^{-O(\beta)}Q_{\delta}*Q_{\delta^{1-O(\beta)}}$, and combining the last two inequalities implies (b).
Finally, \eqref{metric} gives that  $\|x-y\|_2\leqslant\delta^{-O(\beta)}\rho$, for any $x\in G$ and $y\in B(1,4C\beta\log{\frac{1}{\delta}})$ satisfying $d(x,y)<\rho$. This clearly implies (c).

 To finish the proof, assume that $\mu$ additionally satisfies
 $\|\mu*P_{\delta}\|_2\geqslant \delta^{-\alpha}$, and
 $(\mu*\mu)(H^{(\rho)})\leqslant\delta^{-\gamma}\rho^{\kappa}$, for all $\rho\geqslant\delta$ and any proper closed connected subgroup $H<G$,
 for some $\gamma>0$.

By using (a), (c) and Lemma \ref{flat} we get that 

\begin{itemize}
\item supp$(\mu)\subset B(1,2C\beta\log{\frac{1}{\delta}}).$
\item $\|\mu*Q_{\delta^{1-O(\beta)}}\|_2\ll\|\mu*Q_{\delta}\|_2\ll\delta^{-\alpha}$. 
\item $(\mu*\mu)(H^{[\rho]})\leqslant(\mu*\mu)(H^{(\delta^{-O(\beta)}\rho)})\leqslant\delta^{-(\gamma+\kappa O(\beta))}\rho^{\kappa}$,  for all $\rho\geqslant\delta$ and any  proper closed connected subgroup $H<G$.
\end{itemize}

If  $\beta,\gamma>0$ are chosen small enough, then Theorem \ref{BdS2} implies that $\|\mu_{\delta}*\mu_{\delta}\|_2<\delta^{\varepsilon}\|\mu_{\delta}\|_2$. Moreover, if $\beta,\gamma$ are small enough, then  by combining this inequality with (a) and (b) we derive that 
\begin{align*}
\|\mu*\mu*P_{\delta}\|_2\ll\|\mu*\mu*Q_{\delta}\|_2&\leqslant\delta^{-O(\beta)}\|\mu*Q_{\delta^{1-O(\beta)}}*\mu*Q_{\delta^{1-O(\beta)}}\|_2\\&\leqslant\delta^{-O(\beta)+(1-O(\beta))\varepsilon}\|\mu*Q_{\delta^{1-O(\beta)}}\|_2\\&\leqslant\delta^{\varepsilon-O(\beta)}\|\mu*Q_{\delta}\|_2\\&\leqslant\delta^{\varepsilon-O(\beta)}\|\mu*P_{\delta}\|_2\\&\leqslant\delta^{\gamma}\|\mu*P_{\delta}\|_2.
\end{align*}

This concludes the proof of Lemma \ref{BdS}. \hfill$\blacksquare$

\end{appendix}

\bibliographystyle{alpha1}

\begin{thebibliography}{AA}
\bibitem [AE10]{AE10} M. Ab\'{e}rt, G. Elek: {\it Dynamical properties of profinite actions}, Erg. Th. Dynam. Sys. {\bf 32} (2012),  1805-1835.
\bibitem [AJN07]{AJN07} M. Ab\'{e}rt, A. Jaikin-Zapirain, N. Nikolov: {\it The rank gradient from a combinatorial viewpoint} Groups Geom. Dyn. {\bf 5} (2011), 213-230.
\bibitem [BNP08]{BNP08} L. Babai, N. Nikolov, L. Pyber: {\it Product Growth and Mixing in Finite Groups}, 19th ACM-SIAM Symposium on Discrete Algorithms, SIAM, 2008, 248-257.
\bibitem [Ba23]{Ba23} S. Banach: {\it Sur le probl\`{e}me de la mesure}, Fund. Math. {\bf 4} (1923), 7-33.
\bibitem [BdHV08]{BdHV08} M. Bekka, P. de la Harpe, A. Valette: {\it  Kazhdan's property (T)}, New Mathematical Monographs, 11. Cambridge University Press, Cambridge, 2008. xiv+472 pp.
\bibitem [BdS14]{BdS14}  Y. Benoist, N. de Saxc\`{e}: {\it A spectral gap theorem in simple Lie groups}, preprint arXiv:1405.1808.
\bibitem [Bo09]{Bo09} J. Bourgain: {\it Expanders and dimensional expansion}, Comptes Rendus Mathematique {\bf 347} (2009), 357-362.
\bibitem [BG05]{BG05} J. Bourgain, A. Gamburd: {\it Uniform expansion bounds for Cayley graphs of $SL_2(\mathbb F_p)$}, Ann. of Math. (2) {\bf 167} (2008),  625-642.
\bibitem [BG06]{BG06} J. Bourgain, A. Gamburd: {\it On the spectral gap for finitely-generated subgroups of $SU(2)$}, Invent. Math. {\bf 171} (2008), 83-121.
\bibitem [BG10]{BG10} J. Bourgain, A. Gamburd: {\it A spectral gap theorem in $SU(d)$}, J. Eur. Math. Soc. (JEMS) {\bf 14} (2012), 1455-1511.
\bibitem [BY11]{BY11} J. Bourgain, A. Yehudayoff: {\it Expansion in $SL_2(\mathbb R)$ and monotone expanders}, Geom. Funct. Anal. {\bf 23} (2013), 1-41.
\bibitem [Br08]{Br08} E. Breuillard: {\it A strong Tits alternative}, preprint arXiv:0804.1395.
\bibitem [BrG02]{BrG02} E. Breuillard, T. Gelander: {\it On dense free subgroups of Lie groups}, J. Algebra {\bf 261} (2003), 448-467.
\bibitem[BG04]{BG04} E. Breuillard, T. Gelander: {\it A topological Tits alternative},  Ann. of Math. {\bf 166} (2007), 427-474. 
\bibitem[CFW81]{CFW81} A. Connes, J. Feldman, B.Weiss: {\it An amenable equivalence relations is generated by a single transformation},
Ergodic Th. Dynam. Sys. {\bf 1}(1981), 431-450.
\bibitem[CW80]{CW80} A. Connes, B. Weiss: {\it Property T and asymptotically invariant sequences}, Israel J. Math. {\bf 37} (1980), 209-210.
\bibitem[CdH14]{CdH14} Y. de Cornulier, de la Harpe: {\it Metric geometry of locally compact groups}, book in progress, preprint arXiv:1403.3796.
\bibitem[Dr84]{Dr84} V. Drinfeld: {\it Finitely-additive measures on $S^2$ and $S^3$, invariant with respect to rotations}, Funct. Anal. Appl. {\bf 18} (1984), 245-246.
\bibitem[EMO05]{EMO05} A. Eskin, S. Mozes, H. Oh: {\it On Uniform exponential growth for linear groups,}
 Invent. Math. {\bf 160} (2005), 1-30.
 \bibitem[Fu09]{Fu09} A. Furman: {\it A survey of Measured Group Theory}, Geometry, Rigidity, and Group Actions, 296-374, The University of Chicago Press, Chicago and London, 2011.
 \bibitem[FS99]{FS99} A. Furman, Y. Shalom: \textit{Sharp ergodic theorems for group actions and strong ergodicity},
Ergodic Theory Dynam. Systems {\bf 19} (1999), no. 4, 1037-1061.
\bibitem[Ga10]{Ga10} D. Gaboriau: \textit{Orbit equivalence and measured group theory}, In Proceedings of the ICM (Hyderabad, India, 2010), Vol. III, Hindustan Book Agency, 2010, pp. 1501-1527.
\bibitem [GJS99]{GJS99} A. Gamburd, D. Jakobson, P. Sarnak: {\it Spectra of elements in the group ring of $SU(2)$}, J. Eur. Math. Soc. (JEMS) {\bf 1} (1999),  51-85.
\bibitem [He05]{He05} H. Helfgott: {\it Growth and generation in $SL_2(\mathbb Z/p\mathbb Z)$}, Ann. Math.  {\bf 167} (2008), 601-623.
\bibitem [HV12]{HV12} C. Houdayer, S. Vaes: {\it Type III factors with unique Cartan decomposition}, J. Math. Pures Appl. {\bf 100} (2013),  564-590.
\bibitem [Io13]{Io13} A. Ioana: {\it Orbit equivalence and Borel reducibility rigidity for profinite actions with spectral gap}, preprint arXiv:1309.3026.
\bibitem [Io14]{Io14} A. Ioana: {\it Strong ergodicity, property (T), and orbit equivalence rigidity for translation actions}, preprint arXiv:1406.6628, to appear in J. Reine Angew. Math.
\bibitem [dJR79]{dJR79} A. del Junco, J. Rosenblatt: {\it Counterexamples in ergodic theory and number theory}, Math. Ann. {\bf 245} (1979), 185-197.
\bibitem [Ke59]{Ke59} H. Kesten: {\it Symmetric random walks on groups}, Trans. Amer. Math. Soc. {\bf 92} (1959), 336-354.
\bibitem [KS71]{KS71} A. W. Knapp and E. M. Stein: {\it Interwining Operators for Semisimple Groups}, Ann. of Math.
{\bf 93} (1971), 489-578.
\bibitem[LdS13]{LdS13} E. Lindenstrauss, N. de Saxc\'{e}: {\it Hausdorff dimension and subgroups of $SU(2)$}, preprint, to appear in Israel J. Math.
 \bibitem [Lu94]{Lu94} A. Lubotzky: {\it Discrete Groups, Expanding Graphs and Invariant Measures. With an appendix by Jonathan D. Rogawski}, Progress in Mathematics, 125. Birkh\"{a}user Verlag, Basel, 1994. xii+195 pp.
\bibitem[LMR00]{LMR00} A. Lubotzky, S. Mozes, M.S. Raghunathan: {\it The word and Riemannian metrics on lattices of semisimple groups}, Inst. Hautes \'{E}tudes Sci. Publ. Math.  {\bf 91} (2000), 5-53.
\bibitem[LPS86]{LPS86}  A. Lubotzky, R. Phillips, P. Sarnak: {\it Hecke operators and distributing points on the sphere} I, Commun. Pure Appl. Math. {\bf 39} (1986), 149-186.
\bibitem[LPS87]{LPS87} A. Lubotzky, R. Phillips, P. Sarnak: {\it Hecke operators and distributing points on $S^2$} II, Commun. Pure Appl. Math. {\bf 40} (1987), 401-420.
\bibitem[Ma80]{Ma80} G. Margulis: {\it Some remarks on invariant means}, Monatsh. Math. {\bf 90} (1980), 233-235. 
 \bibitem [Ma82]{Ma82} G. Margulis: {\it Finitely-additive invariant measures on Euclidian spaces}, Ergodic Theory Dynam. Systems {\bf 2} (1982), 383-396.
 \bibitem [Ma91]{Ma91} G. Margulis: {\it Discrete subgroups of semisimple Lie groups}, Ergebnisse der Mathematik und ihrer Grenzgebiete (3) [Results in Mathematics and Related Areas (3)], 17. Springer-Verlag, Berlin, 1991. x+388 pp. 
 \bibitem[MW83]{MW} D.W. Masser and G. W\"{u}stholz, 
{\it Fields of large transcendence degree generated by values of elliptic functions,} 
Invent. Math. {\bf 72} (1983), no. 3, 407-464.
 \bibitem[MZ55]{MZ55} D. Montgomery, L. Zippin: {\it Topological transformation groups}, Interscience Publishers, New York-London, 1955. xi+282 pp.
 \bibitem[Oh05]{Oh05} H. Oh: {\it The Ruziewicz problem and distributing points on homogeneous spaces of a compact Lie group},
Israel J. Math. (Furstenberg volume), {\bf 149} (2005), 301-316.
\bibitem[OP07]{OP07} N. Ozawa, S. Popa: {\it On a class of II$_1$ factors with at most one Cartan subalgebra}, Ann. of Math. (2), {\bf 172}
(2010), 713-749.
\bibitem[OW80]{OW80} D. Ornstein, B.Weiss: {\it Ergodic theory of amenable groups. I. The Rokhlin lemma.,}
Bull. Amer. Math. Soc. (N.S.) {\bf 1} (1980), 161-164.
 \bibitem[Po07]{Po07} S. Popa: \textit{Deformation and rigidity for group actions and von Neumann algebras}, In Proceedings of the ICM (Madrid, 2006), Vol. I, European Mathematical Society Publishing House, 2007, 445-477.
\bibitem[PV11]{PV11} S. Popa, S. Vaes: {\it Unique Cartan decomposition for II$_1$ factors arising from arbitrary actions of free
groups}, Acta Math. {\bf 212} (2014), no. 1, 141-198.
 \bibitem[Ra72]{Ra72} M.S. Raghunathan: {\it Discrete subgroups of Lie groups}, Ergebnisse der Mathematik und ihrer Grenzgebiete, Band 68. Springer-Verlag, New York-Heidelberg, 1972. ix+227 pp. 
\bibitem[RS88]{RS88} F. Ricci, E. Stein: {\it Harmonic analysis on nilpotent groups and singular integrals. II. Singular kernels supported on sub manifolds.} J. Funct. Anal., {\bf 78} (1988), 56-84.
\bibitem[Ro81]{Ro81} J. Rosenblatt: {\it Uniqueness of invariant means for measure-preserving transformations}, Trans. Amer. Math. Soc. {\bf 265} (1981).
\bibitem[SGV11]{SGV11} A. Salehi Golsefidy, P. Varj\'{u}: {\it Expansion in perfect groups}, Geom. Funct. Anal. {\bf 22} (2012), no. 6, 1832-1891.
\bibitem [SX91]{SX91} P. Sarnak, X. Xue: {\it Bounds for multiplicities of automorphic representations}, Duke Math. J., {\bf 64} (1991), 207-227.
\bibitem [dS14]{dS14} N. de Saxc\'{e}: {\it A product theorem in simple Lie groups}, preprint  arXiv:1405.2003, to appear in Geom. Funct. Anal.
\bibitem[Sc80]{Sc80} K. Schmidt: {\it Asymptotically invariant sequences and an action of $SL(2;\mathbb Z)$ on the $2$-sphere}, Israel J. Math.
{\bf 37} (1980), 193-208.
\bibitem [Sc81]{Sc81} K. Schmidt: {\it Amenability, Kazhdan's property T, strong ergodicity and invariant means for ergodic group-actions}, Ergodic Theory Dynamical Systems {\bf 1} (1981),  223-236.
\bibitem [St93]{St93} E. Stein: {\it Harmonic Analysis: Real-variable Methods, Orthogonality and Oscillatory Integrals}, Princeton University Press, 1993.

\bibitem[Su81]{Su81} D. Sullivan: {\it For $n>3$ there is only one finitely additive rotationally invariant measure on the $n$-sphere on all Lebesgue measurable sets}, Bull. Am. Math. Soc. {\bf 4} (1981), 121-123.
\bibitem[Va10]{Va10} P. Varj\'{u}: {\it Expansion in $SL_d(\mathcal O_K/I)$, $I$ square-free}, J. Eur. Math. Soc. (JEMS) {\bf 14} (2012),  273-305.
\bibitem[Ta06]{Ta06} T. Tao: {\it Product set estimates for non-commutative groups}, Combinatorica {\bf 28} (2008), no. 5, 547-594.
\bibitem[Ta15]{Ta15} T. Tao: {\it Expansion in finite simple groups of Lie type},  Graduate Studies in Mathematics, 164. American Mathematical Society, Providence, RI, 2015. xiv+303 pp. 
\bibitem [Wa93]{Wa93} S. Wagon: {\it The Banach-Tarski paradox. With a foreword by Jan Mycielski.} Encyclopedia of Mathematics and its Applications, 24. Cambridge University Press, Cambridge, 1985. xvi+251 pp.
\bibitem  [Zi78]{Zi78} R. J. Zimmer: {\it Amenable ergodic group actions and an application to Poisson boundaries of random walks}, J. Funct. Anal. {\bf 27} (1978),  350-372.
\end{thebibliography}

\end{document}